%% file: pf_article.tex
\documentclass[10pt,a4paper]{amsart}
\pdfoutput=1

\usepackage[austrian,english]{babel}
\usepackage[latin1]{inputenc}
\usepackage[a4paper]{geometry}

\usepackage{hyperref}
\usepackage{color}
\usepackage{tikz}
\usetikzlibrary{arrows,decorations.pathreplacing,patterns}
\usepackage[font=footnotesize]{caption}

\usepackage{amsmath}
\usepackage{amssymb}
\usepackage{faktor}
\usepackage{mathtools}

\def\Q{\check{Q}}
\def\wa{\smash{\widetilde{W}}}

\def\ac{A_{\circ}}

\def\wdom{\wshi^{\mathsf{dom}}}

\def\wshi{W_{\mathsf{Shi}}}
\def\A{\mathcal{A}}

\newcommand{\nnpark}{\op{Park}(\Phi)}
\newcommand{\walls}{J(\lambda)}

\newcommand{\LS}{\mathrlap{\overleftarrow{S}}\rule{.8em}{0pt}}
\newcommand{\RS}{\mathrlap{\overrightarrow{S}}\rule{.8em}{0pt}}
\newcommand{\ld}{\mathcal L^{\bullet}}
\newcommand{\bd}{\mathcal B^{\bullet}}

\newcommand{\waf}{\widetilde{w}}
\newcommand{\wf}{\waf_f}

\def\n{K}
\def\area{\op{area}}
\def\dinv{\op{dinv}}
\def\bounce{\op{bounce}}

\newcommand{\N}{\mathbb N}
\newcommand{\R}{\mathbb R}
\newcommand{\Z}{\mathbb Z}
\renewcommand{\S}{\mathfrak S}
\newcommand{\sgn}{\operatorname{sgn}}

\newcommand{\abs}[1]{\left| #1 \right|}

\newcommand{\skal}[1]{\langle #1 \rangle}
\newcommand{\op}[1]{\operatorname{#1}}
\newcommand{\affS}{\smash{\stackrel{\sim}{\smash{\mathfrak{S}}\rule{0pt}{1.05ex}}}}

\newcommand{\myfontA}[1]{\textnormal{\textbf{#1}}}
\newcommand{\myfontB}[2][]{\myfontA{#2}\textnormal{#1}}

\newcommand{\sk}{

\smallskip}

\newcommand{\myi}[1]{

\smallskip\myfontB{#1.}}

\newcounter{mythrm}\numberwithin{mythrm}{section}
\numberwithin{equation}{section}

\newenvironment{mymathenvironmentA}[3]{
\refstepcounter{mythrm}\label{#3:#2}
\sk\myfontA{#3 \arabic{section}.\arabic{mythrm}#1.}\begin{itshape}}
{\end{itshape}}

\newenvironment{mymathenvironmentB}[3]{
\refstepcounter{mythrm}\label{#3:#2}
\smallskip\myfontB{#3 \arabic{section}.\arabic{mythrm}#1.}}
{\hfill$\Box$}

\newenvironment{mydef}[2][]{
\begin{mymathenvironmentB}{#1}{#2}{Definition}}
{\end{mymathenvironmentB}}

\newenvironment{myex}[2][]{
\begin{mymathenvironmentB}{#1}{#2}{Example}}
{\end{mymathenvironmentB}}

\newenvironment{mylem}[2][]{
\begin{mymathenvironmentA}{#1}{#2}{Lemma}}
{\end{mymathenvironmentA}}

\newenvironment{myprob}[2][]{
\begin{mymathenvironmentB}{#1}{#2}{Problem}}
{\end{mymathenvironmentB}}

\newenvironment{myproof}{
\smallskip\myfontA{Proof.}}
{\hfill$\blacksquare$}

\newenvironment{myproof*}[1]{
\smallskip\myfontA{Proof #1.}}
{\hfill$\blacksquare$}

\newcommand{\noproof}[1]{#1\hfill$\blacksquare$\sk}

\newenvironment{myprop}[2][]{
\begin{mymathenvironmentA}{#1}{#2}{Proposition}}
{\end{mymathenvironmentA}}

\newenvironment{myrem}[2][]{
\begin{mymathenvironmentB}{#1}{#2}{Remark}}
{\end{mymathenvironmentB}}

\newenvironment{mythrm}[2][]{
\begin{mymathenvironmentA}{#1}{#2}{Theorem}}
{\end{mymathenvironmentA}}

\newcommand{\refx}[2]{#2~\ref{#2:#1}}

\newcommand{\refd}[1]{\refx{#1}{Definition}}
\newcommand{\refe}[1]{\refx{#1}{Example}}
\newcommand{\reff}[1]{\refx{#1}{Figure}}
\newcommand{\refl}[1]{\refx{#1}{Lemma}}
\newcommand{\refp}[1]{\refx{#1}{Proposition}}
\newcommand{\refq}[1]{\eqref{eq:#1}}

\newcommand{\refs}[1]{\refx{#1}{Section}}
\newcommand{\reft}[1]{\refx{#1}{Theorem}}

\newcommand{\ballot}[1]{
	\pgfmathsetmacro{\y}{2*(#1)}
	\draw[gray] (0,0)--(0,\y);
	\foreach \x in {1,...,#1}{
		\pgfmathsetmacro{\y}{2*(#1)-\x+1}
		\draw[gray] (0,\x)--(\x,\x)--(\x,\y)--(0,\y);
	}
}
\newcommand{\ballotodd}[1]{
	\pgfmathsetmacro{\y}{2*(#1)-1}
	\draw[gray] (0,0)--(0,\y);
	\pgfmathsetmacro{\y}{(#1)-1}
	\draw[gray] (0,#1)--(\y,#1);
	\foreach \x in {1,...,\y}{
		\pgfmathsetmacro{\z}{2*(#1)-\x}
		\draw[gray] (0,\x)--(\x,\x)--(\x,\z)--(0,\z);
	}
}

\newcommand{\ballotsign}[3]{
\draw[xshift=#1cm,yshift=#2cm](.5,0)node[circle,fill=white,inner sep=-1.3pt]{$#3$};
}

\begin{document}

\title{On parking functions and the zeta map in types $B,C$ and $D$}
\date{2016}

\author{Robin Sulzgruber}
\thanks{Research supported by the Austrian Science Fund (FWF), grant S50-N15 in the framework of the Special Research Program ``Algorithmic and Enumerative Combinatorics'' (SFB F50).}
\address{Fakult\"at f\"ur Mathematik, Universit\"at Wien, Oskar-Morgenstern-Platz 1, 1090 Wien, Austria}
\email{robin.sulzgruber@univie.ac.at}

\author{Marko Thiel}
\address{Institut f\"ur Mathematik, Universit\"at Z\"urich, Winterthurerstrasse 190, 8057 Z\"urich, Switzerland}
\email{marko.thiel@math.uzh.ch}

\begin{abstract} 
Let $\Phi$ be an irreducible crystallographic root system with Weyl group $W$, coroot lattice $\Q$ and Coxeter number $h$. Recently the second named author defined a uniform $W$-isomorphism $\zeta$ between the finite torus $\Q/(mh+1)\Q$ and the set of non-nesting parking fuctions $\op{Park}^{(m)}(\Phi)$. If $\Phi$ is of type $A_{n-1}$ and $m=1$ this map is equivalent to a map defined on labelled Dyck paths that arises in the study of the Hilbert series of the space of diagonal harmonics. 

In this paper we investigate the case $m=1$ for the other infinite families of root systems ($B_n$, $C_n$ and $D_n$). In each type we define models for the finite torus and for the set of non-nesting parking functions in terms of labelled lattice paths. The map $\zeta$ can then be viewed as a map between these combinatorial objects. Our work entails new bijections between (square) lattice paths and ballot paths.
\end{abstract}

\thispagestyle{empty}
\maketitle






\pagenumbering{arabic}
\pagestyle{headings}

\input{uniform}
\input{preps}
\input{typeC}
\input{typeD}

\input{typeB}

\input{conclusion}


\bibliographystyle{alpha}
\bibliography{robinbib}



\end{document}

%% file: uniform.tex
\section{Introduction}

\subsection{The space of diagonal harmonics}
One of the most well-studied objects in algebraic combinatorics is the space of diagonal harmonics of the symmetric group $\S_n$.
Its Hilbert series has two combinatorial interpretations \cite[Conjecture 5.2]{Haglund2008} \cite{CarMel2016}: 
\begin{align*}
\mathcal{DH}(n;q,t)=\sum_{(P,v)\in \op{Vert}(A_{n-1})}q^{\dinv'(P,v)}t^{\area(P,v)}=\sum_{(D,w)\in \op{Diag}(A_{n-1})}q^{\area'(D,w)}t^{\bounce(D,w)},
\end{align*}
where $ \op{Vert}(A_{n-1})$ is the set of \emph{parking functions} of length $n$, viewed as vertically labelled Dyck paths, and $ \op{Diag}(A_{n-1})$ is the set of \emph{diagonally labelled} Dyck paths with $2n$ steps.
There is a bijection $\zeta_A$ due to 
Haglund and Loehr 
\cite{HagLoe2005} that maps $ \op{Vert}(A_{n-1})$ to $ \op{Diag}(A_{n-1})$ and sends the bistatistic $(\dinv',\area)$ to $(\area',\bounce)$, demonstrating the second equality.

\subsection{The zeta map beyond type \emph{A}}
The combinatorial objects $ \op{Vert}(A_{n-1})$ and $ \op{Diag}(A_{n-1})$ may be viewed as corresponding to the type $A_{n-1}$ cases of more general algebraic objects associated to any irreducible crystallographic root system $\Phi$.
These are, respectively, the \emph{finite torus} \mbox{$\Q/(h+1)\Q$} and the set of \emph{non-nesting parking functions} $\nnpark$ of $\Phi$. Here $\Q$ is the coroot lattice and $h$ is the Coxeter number of $\Phi$.
\sk
Recently, the second named author has defined a uniform bijection $\zeta$ from \mbox{$\Q/(h+1)\Q$} to $\nnpark$ \cite[Thm 15.4]{Thiel2015}.
Its most salient property is that it commutes with the action of the \emph{Weyl group} $W$ of $\Phi$ on both sets.
In type $A_{n-1}$, the uniform map $\zeta$ corresponds to the combinatorial map $\zeta_A$ \cite[Thm 16.3]{Thiel2015}. This is illustrated by the following commutative diagram of bijections:
\begin{center}
\begin{tikzpicture}
\draw (0,0) node(q){$\Q/(n+1)\Q$};
\draw (4,0) node(park){$\mathrm{Park}(\Phi)$};
\draw (0,1.4) node(v){$\op{Vert}(A_{n-1})$};
\draw (4,1.4) node(d){$\op{Diag}(A_{n-1})$};
\draw[->] (q)--node[below]{$\zeta$}(park);
\draw[->] (v)--node[above]{$\zeta_A$}(d);
\draw[->] (v)--node[left]{$\psi$}(q);
\draw[->] (d)--node[right]{$\varphi$}(park);
\end{tikzpicture}
\end{center}
In this paper, for each of the other classical types $B_n,C_n$ and $D_n$, we define combinatorial models for $\Q/(h+1)\Q$ and $\nnpark$, in terms of vertically labelled paths and diagonally labelled paths respectively.
Furthermore, we define the combinatorial bijections $\zeta_B$, $\zeta_C$, and $\zeta_D$ which correspond to the uniform bijection $\zeta$ in those types.

Besides their meaning in the theory of reflection groups these maps can also be appreciated from a purely combinatorial point of view.
As such $\zeta_B$ and $\zeta_C$ are new bijections between lattice paths in an $n\times n$-square and ballot paths with $2n$ steps, both of which are well-known to be counted by the central binomial coefficients $\binom{2n}{n}$.
Furthermore $\zeta_D$ is a new bijection between lattice paths in an $(n-1)\times n$ rectangle and ballot paths with $2n-1$ steps.

\subsection{Outline} This paper is structured as follows. In \refs{defs} we present the needed facts about Weyl groups and root systems. In particular we recall the relations between objects such as the Shi arrangement, the finite torus and the non-nesting parking functions, and discuss in detail the definition of the uniform zeta map.
In \refs{BCD} we turn our attention to the root systems of types $B_n, C_n$ and $D_n$. In particular we recall the combinatorial models for their affine Weyl groups in terms of affine permutations and show how to obtain the decomposition of an affine permutation into a product of an element of the Weyl group and a translation by an element in the coroot lattice combinatorially.

We then develop the combinatorics of the zeta map of type $C_n$ in \refs{typeC}, type $D_n$ in \refs{typeD} and type $B_n$ in \refs{typeB}.
Each of these sections is structured similarly.
First we define combinatorial models for the finite torus in terms of (labelled) lattice paths and the non-nesting parking functions in terms of (labelled) ballot paths.
Secondly we define an area vector, a diagonal reading word and two versions of a combinatorial zeta map, one between labelled objects and one between unlabelled objects.
We then prove that the combinatorial zeta map coincides with the uniform zeta map of the given type, and is thus a bijection.
\refs{typeC} includes an alternative description of the zeta map of type $C_n$ akin to the sweep map and an additional section on lattice path statistics.

Finally, in \refs{end} we offer some possible directions for further research.

\sk
An extended abstract~\cite{SulThi2015} of this paper containing mostly the results of \refs{typeC} has appeared in the conference proceedings of FPSAC 2015 in Daejeon.

\section{Definitions and Preliminaries}\label{Section:defs}

Given $n\in\N$ we set $[n]=\{1,2,\dots,n\}$ and $[\pm n]=[n]\cup\{-i:i\in[n]\}$.

\subsection{Root systems and Weyl groups}

Let $\Phi$ be an irreducible crystallographic root system of rank $r$, with simple system $\Delta=\{\alpha_1,\alpha_2,\ldots,\alpha_r\}$, positive system $\Phi^+$ and ambient space $V$.
For background on root systems and reflection groups see \cite{Humph1990}.
For $\alpha\in\Phi$, let $s_{\alpha}$ be the reflection in the hyperplane 
\begin{align*}
H_{\alpha}=\{x\in V:\langle x,\alpha\rangle=0\}.
\end{align*}
Then the \emph{Weyl group} $W$ of $\Phi$ is the group of linear automorphisms of $V$ generated by all the $s_{\alpha}$ with $\alpha\in\Phi$.
Define the \emph{Coxeter arrangement} of $\Phi$ as the central hyperplane arrangement in $V$ consisting of all the hyperplanes $H_{\alpha}$ for $\alpha\in\Phi$.
The connected components of the complement of the union of these hyperplanes are called \emph{chambers}.
The Weyl group $W$ acts simply transitively on the chambers, so we define the \emph{dominant chamber} as
\begin{align*}
C=\{x\in V:\langle x,\alpha\rangle>0\text{ for all }\alpha\in\Delta\}
\end{align*}
and write every chamber as $wC$ for a unique $w\in W$.

The Weyl group $W$ is a \emph{Coxeter group} generated by $S=\{s_{\alpha}:\alpha\in\Delta\}$. If $I\subseteq \Delta$, we call $W_I=\langle\{s_{\alpha}:\alpha\in I\}\rangle$ a \emph{standard parabolic subgroup} of $W$.
A \emph{parabolic subgroup} is any subgroup conjugate to a standard parabolic subgroup.
\subsection{The affine Weyl group}
Define a partial order on $\Phi^+$ by letting
$\alpha\leq\beta$ if and only if $\beta-\alpha$ can be written as a sum of positive roots.
The set $\Phi^+$ with this partial order is called the \emph{root poset} of $\Phi$. It has a unique maximal element, called the \emph{highest root} $\tilde{\alpha}$.
Write $\tilde{\alpha}=\sum_{i=1}^rc_i\alpha_i$ and define the \emph{Coxeter number} of $\Phi$ as $h=1+\sum_{i=1}^rc_i$.

For $\alpha\in\Phi$ and $k\in\Z$, let $s_{\alpha,k}$ be the reflection in the affine hyperplane 
\begin{align*}
H_{\alpha,k}=\{x\in V:\langle x,\alpha\rangle=k\}.
\end{align*}
Then the \emph{affine Weyl} group $\wa$ of $\Phi$ is the group of affine transformations of $V$ generated by all the $s_{\alpha,k}$ for $\alpha\in\Phi$ and $k\in\Z$.
Define the \emph{affine Coxeter arrangement} as the affine hyperplane arrangement in $V$ consisting of all the $H_{\alpha,k}$ for $\alpha\in\Phi$ and $k\in\Z$.
The connected components of the complement of the union of these hyperplanes are called \emph{alcoves}.
The affine Weyl group $\wa$ acts simply transitively on the alcoves, so we define the \emph{fundamental alcove} as
\begin{align*}
\ac=\{x\in V:\langle x,\alpha\rangle>0\text{ for all }\alpha\in\Delta\text{ and }\langle x,\tilde{\alpha}\rangle<1\}
\end{align*}
and write every alcove as $\widetilde w\ac$ for a unique $\widetilde w\in\wa$.

For each root $\alpha\in\Phi$ we define its coroot as $\alpha^{\vee}=2\alpha/\skal{\alpha,\alpha}$. The coroot lattice $\Q$ of $\Phi$ is the integer span of the coroots. The simple coroots are a natural lattice basis for it: $\Q=\bigoplus_{\alpha\in\Delta}\Z\alpha^{\vee}$. The affine Weyl group $\wa$ acts faithfully on the coroot lattice $\Q$. To each coroot lattice point $\mu\in\Q$ corresponds a translation $t_{\mu}:V\to V$ given by $t_{\mu}(x)=x+\mu$. If we identify $\Q$ with its translation group we may write $\wa=W\ltimes\Q$ as a semidirect product.

The affine Weyl group is also a Coxeter group with generating set $\widetilde{S}=S\cup\{s_{\tilde{\alpha},1}\}$.
The Weyl group $W=\langle S\rangle$ is a standard parabolic subgroup of the Coxeter group $\wa$. Each (left) coset in the quotient $\wa/W$ contains a unique element of minimal length. The minimal length representatives are called \emph{(left) Gra{\ss}mannian}. An element $\widetilde w\in\wa$ is Gra{\ss}mannian if and only if the alcove $\widetilde w^{-1}\ac$ lies in the dominant chamber.

\subsection{The Shi arrangement}\label{Section:shi}
Define the \emph{Shi arrangement} as the hyperplane arrangement consisting of the hyperplanes $H_{\alpha,k}$ for $\alpha\in\Phi^+$ and $k=0,1$. The complement of the union of these hyperplanes falls apart into connected components called the \emph{regions} of the arrangement. 
The hyperplanes that support facets of a region $R$ are called the \emph{walls} of $R$. Those walls of $R$ that do not contain the origin and separate $R$ from the origin are called the \emph{floors} of $R$.
Define the walls and floors of an alcove similarly. Notice that every wall of a region is a hyperplane of the Shi arrangement, but the walls of an alcove need not be. We call a region or alcove \emph{dominant} if it is contained in the dominant chamber.
\sk
\begin{mythrm}{thrm:min} {\textnormal{\cite[Prop~7.1]{Shi1987}}}
Every region $R$ of the Shi arrangement has a unique \emph{minimal alcove} $\waf_R\ac\subseteq R$, which is the alcove in $R$ closest to the origin.
\end{mythrm}
\begin{figure}[h]
\begin{center}
 \resizebox*{6cm}{!}{\includegraphics{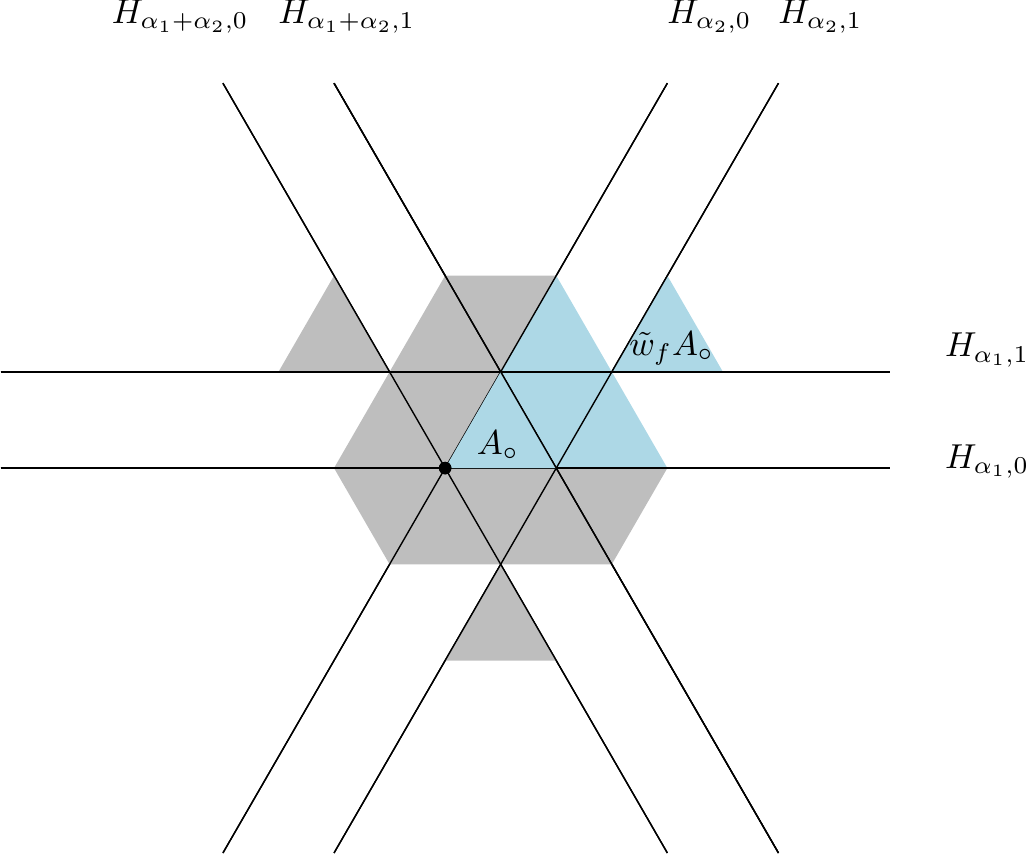}}\\
\end{center}
\caption{The 16 regions of the Shi arrangement of type $A_2$ with their minimal alcoves. The minimal alcoves of the 5 dominant regions are coloured blue.}
\end{figure}
\sk
Furthermore the floors of the minimal alcove $\waf_R\ac$ of $R$ are exactly the floors of $R$ \cite[Thm 12.4]{Thiel2015}.
We define $\wshi=\{\waf_R: R\text{ is a Shi region}\}$. The corresponding alcoves $\waf_R\ac$ we call \emph{Shi alcoves}. That is, we call an alcove a Shi alcove if it is the minimal alcove of the Shi region containing it.
Define $\wdom=\{\waf_R: R\text{ is a dominant Shi region}\}$.

An \emph{order filter} in the root poset is a set $J\subseteq\Phi^+$ such that whenever $\alpha\in J$ and $\beta\geq\alpha$ then $\beta\in J$.
For a dominant Shi region $R$, we consider 
\[J=\{\alpha\in\Phi^+:\langle x,\alpha\rangle>1\text{ for all }x\in R\}.\]
Then $J$ is an order filter in the root poset. In fact, the map $R\mapsto J$ is a bijection from the set of dominant Shi regions to the set of order filters in the root poset \cite[Thm~1.4]{Shi1997}.
The set of minimal elements of $J$ is
\[A=\{\alpha\in\Phi^+:H_{\alpha,1}\text{ is a floor of }R\},\]
an antichain in the root poset. Since the correspondence between order filters and the antichains of their minimal elements is bijective, the map $R\mapsto A$ is a bijection from the set of dominant Shi regions to the set of antichains in the root poset.

We conclude this section by introducing an important affine transformation considered by Athanasiadis \cite[Section 4]{Athan2005}. Denote by $\wf\in\wa$ the element of the affine Weyl group such that $\wf\ac$ is the minimal alcove of the dominant Shi region furthest from the origin, corresponding to the order filter $J=\Phi^+$ and the antichain $A=\Delta$:
\begin{align*}
\wf\ac=\big\{x\in V:\skal{x,\alpha}>1\text{ for all }\alpha\in\Delta\text{ and }\skal{x,\tilde\alpha}<h\big\},
\end{align*}
where $h$ denotes the Coxeter number of $\Phi$. 

\subsection{Affine roots}

Let $\delta$ denote a formal variable and set $\widetilde V=V\oplus\R\delta$. We define the set of affine roots as
\begin{align*}
\widetilde\Phi=\big\{\alpha+k\delta:\alpha\in\Phi,k\in\Z\big\}\subseteq\widetilde V.
\end{align*}
To an affine root $\alpha+k\delta$ we associate the half space
\begin{align*}
\mathcal H_{\alpha+k\delta}=\big\{x\in V:\skal{x,\alpha}>-k\big\}.
\end{align*}
The natural action of the affine Weyl group on half spaces carries over to affine roots as follows. Suppose $\widetilde w\in\wa$ has the unique decomposition $\widetilde w=t_{\mu}w$ with $\mu\in\Q$ and $w\in W$. Then
\begin{align*}
\widetilde w\cdot(\alpha+k\delta)
=w\cdot\alpha + (k-\skal{\mu,w\cdot\alpha})\delta,
\end{align*}
where $w\cdot\alpha$ denotes the usual action of $W$ on $\Phi$.
Define the positive affine roots as
\begin{align*}
\widetilde\Phi^+=\Phi^+\cup\big\{\alpha+k\delta:\alpha\in\Phi,k>0\big\},
\end{align*}
and the simple affine roots as
\begin{align*}
\widetilde\Delta=\Delta\cup\{-\tilde\alpha+\delta\}.
\end{align*}
Thereby positive affine roots correspond to those half spaces that contain the fundamental alcove. Simple affine roots correspond to half spaces that contain the fundamental alcove and share one of its walls. As a consequence we obtain the following lemma that characterises the floors and separating hyperplanes of an alcove.
\begin{mylem}{affroots} Let $\widetilde w\in\wa$ and $\alpha+k\delta\in\widetilde\Phi^+$.
\myi{(i)} The hyperplane $H_{\alpha,-k}$ separates $\widetilde w\ac$ from $\ac$ if and only if $\widetilde w^{-1}\cdot(\alpha+k\delta)\in-\widetilde\Phi^+$.
\myi{(ii)} The hyperplane $H_{\alpha,-k}$ is a floor of $\widetilde w\ac$ if and only if $\widetilde w^{-1}\cdot(\alpha+k\delta)\in-\widetilde\Delta$.
\end{mylem}\noproof{}
\sk
Translating the walls of $\wf\ac$ into affine roots we obtain the following lemma.
\begin{mylem}{wf} The element $\wf\in\wa$ is uniquely determined by the property
\begin{align*}
\wf(\widetilde\Delta)=(\Delta-\delta)\cup\{-\tilde\alpha+h\delta\}.
\end{align*}
\end{mylem}\noproof{}

\subsection{The finite torus}\label{Section:torus}
The finite torus is defined as the quotient of the coroot lattice
\begin{align*}
T=\Q/(h+1)\Q.
\end{align*}
The action of the Weyl group on the coroot lattice induces an action on the finite torus. The following result due to Haiman provides a very useful description of the orbits under this action.

\begin{mythrm}{stab} {\textnormal{\cite[Lemma~7.4.1]{Haiman1994}}} The set $\Q\cap(h+1)\overline{\ac}$ forms a system of representatives for the $W$-orbits of $\Q/(h+1)\Q$. The stabiliser subgroup $H=\op{Stab}(\lambda)\leq W$ of an element $\lambda+(h+1)\Q$ of the finite torus, where $\lambda\in\Q\cap(h+1)\overline{\ac}$, is generated by $\{s_{\alpha}:\alpha\in\walls\}$, where 
\begin{align*}
\walls=\big\{\alpha\in\Delta\cup\{-\tilde\alpha\}:\lambda\text{ is contained in the wall of $(h+1)\ac$ perpendicular to }\alpha\big\}.
\end{align*}
\end{mythrm}
By a result of Sommers \cite[Thm~6.4]{Sommers2005}, there exists a $w\in W$ such that $w\cdot\walls\subseteq\Delta$.
So in particular the stabiliser subgroup $H$ of $\lambda+(h+1)\Q$ in $W$ is a parabolic subgroup of $W$,
and any coset in $W/H$ contains a unique element $u$ such that $u\cdot\walls\subseteq\Phi^+$. We call this the canonical representative of the coset with respect to $\walls$.
We deduce the following lemma.

\begin{mylem}{canon}
We may represent any element of the finite torus uniquely as $u\cdot\lambda$ where $\lambda\in\Q\cap(h+1)\overline{\ac}$ and $u\in W$ is such that $u\cdot\walls\subseteq\Phi^+$.
\end{mylem}
\sk
Athanasiadis has defined a bijection from dominant Shi alcoves to $W$-orbits of $\Q/(h+1)\Q$.
\begin{mythrm}{rho} {\textnormal{\cite[Thm~4.2]{Athan2005}}} The map $\mathcal \rho:\wdom\to\Q\cap(h+1)\overline{\ac}$ given by
\begin{align*}
\waf_R\mapsto \wf\waf_R^{-1}\cdot0
\end{align*}
is a bijection.
\end{mythrm}
\sk
The second named author has recently extended this to a bijection called Anderson map between the minimal Shi alcoves and the finite torus.
\begin{mythrm}{anderson} {\textnormal{\cite[Thms~8.4 and 12.6]{Thiel2015}}} The map $\mathcal A:\wshi\to\Q/(h+1)\Q$ given by
\begin{align*}
\waf_R\mapsto -\waf_R\wf^{-1}\cdot0+(h+1)\Q
\end{align*}
is a bijection.
\end{mythrm}
\sk
The way in which $\A$ extends $\rho$ is that the following diagram commutes:
\begin{center}
\begin{tikzpicture}
\draw (0,0) node(q){$\wdom$};
\draw (4,0) node(park){$\Q\cap(h+1)\overline{\ac}$};
\draw (0,1.4) node(v){$\wshi$};
\draw (4,1.4) node(d){$\Q/(h+1)\Q$};
\draw[->] (q)--node[below]{$\rho$}(park);
\draw[->] (v)--node[above]{$\A$}(d);
\draw[<-] (v)--node[left]{$\iota$}(q);
\draw[->] (d)--node[right]{$\pi$}(park);
\end{tikzpicture}
\end{center}
Here $\iota$ is the natural inclusion and $\pi$ is the map that sends an element of $\Q/(h+1)\Q$ to the element of $\Q\cap(h+1)\overline{\ac}$ that represents its $W$-orbit.
\subsection{Non-nesting parking functions}\label{nnsec}
The non-nesting parking functions were introduced by Armstrong, Reiner and Rhoades~\cite{ARR2015} as a model for the regions of the Shi arrangement that carries an action of $W$.
Define an equivalence relation on the set of pairs of an element of the Weyl group and an antichain in the root poset $\Phi^+$ by $(w,A)\sim(w',A')$ if and only if $A=A'$ and $wH=w'H$ where $H\leq W$ denotes the subgroup generated by $\{s_{\alpha}:\alpha\in A\}$. In other words, $H$ is the stabiliser subgroup of $\bigcap_{\alpha\in A}H_{\alpha}$ in $W$. We define the set of \emph{non-nesting parking functions} of $\Phi$ as the set of equivalence classes
\begin{align*}
\op{Park}(\Phi)=\big\{[w,A]:w\in W,A\subseteq\Phi^+\text{ is an antichain}\big\}.
\end{align*}
Using \cite[Thm~6.4]{Sommers2005} the same reasoning that lead to \refl{canon} in the previous section allows us to chose a canonical representative in each class.
\begin{mylem}{nnpark} Every non-nesting parking function $[w,A]$ contains a unique representative $(u,A)$ such that $u\cdot A\subseteq\Phi^+$.
\end{mylem}
\sk
We will also call $u$ the canonical coset representative of $wH$. Note that this depends implicitly not just on $H$, but also on the antichain $A$.
\sk
The Weyl group $W$ acts on $\nnpark$ by
$$u\cdot[w,A]=[uw,A]$$
for $u\in W$. So the orbits of this action are indexed by the antichains in the root poset.
Non-nesting parking functions are naturally in bijection with the regions of the Shi arrangement.
\begin{mythrm}{shitheta} {\textnormal{\cite[Thm~15.2]{Thiel2015}}} The map $\Theta^{-1}:\wshi\to\op{Park}(\Phi)$ given by $\waf_R\mapsto [w,A]$,
where $\waf_R\ac\subseteq wC$ and $A\subseteq\Phi^+$ consists of the roots $\alpha$ such that $H_{\alpha,1}$ is a floor of $w^{-1}\waf_R\ac$, is a bijection.
\end{mythrm}

\sk
Note that a similar bijection using ceilings instead of floors is given in \cite[Prop~10.3]{ARR2015}. We remark that the map $\Theta^{-1}$ extends the natural bijection between the dominant regions of the Shi arrangement and the antichains in the root poset.

\subsection{The uniform zeta map}\label{Section:uniform}
The \emph{zeta map} is defined as the composition of the two bijections introduced in the previous sections
\begin{align}\label{eq:HLzetauni}
\zeta:\Q/(h+1)\Q&\to\op{Park}(\Phi),\\\notag
u\cdot\lambda+(h+1)\Q&\mapsto[w,A]=\left(\Theta^{-1}\circ\mathcal A^{-1}\right)\big(u\cdot\lambda+(h+1)\Q\big).
\end{align}
By a slight abuse of notation we also write $\zeta(u,\lambda)=[w,A]$, where $u\cdot\lambda$ is always understood to be the unique representative as in \refl{canon}.
Since the zeta map commutes with the action of $W$ on both $\Q/(h+1)\Q$ and $\nnpark$, it restricts to a bijection between the $W$-orbits on both sets:
\begin{align}\label{eq:zetauni}
\overline{\zeta}:\Q\cap(h+1)\overline{\ac}&\to\op{Antichains}(\Phi^+),\\\notag
\lambda&\mapsto A.
\end{align}
This can also be seen as an inverse of the bijection $\rho$. To be precise, we have a commutative diagram:
\begin{center}
\begin{tikzpicture}
\draw (0,0) node(q){$\op{Antichains}(\Phi^+)$};
\draw (4,0.7) node(park){$\Q\cap(h+1)\overline{\ac}$};
\draw (0,1.4) node(v){$\wdom$};
\draw[<-] (q)--node[below]{$\overline{\zeta}$}(park);
\draw[->] (v)--node[above]{$\rho$}(park);
\draw[<->] (v)--node[left]{}(q);
\end{tikzpicture}
\end{center}
Here the vertical map is the natural bijection mentioned in \refs{shi}.
\sk

In the remainder of this section we explore $\zeta$ in more detail to obtain a more explicit description of $[w,A]$ in terms of the pair $(u,\lambda)$. To this end we recall some arguments that can be found in similar form for example in \cite[Sec~4]{Athan2005} or \cite[Sec~8.5]{Thiel2015}.

Suppose that $R$ is a Shi region with minimal alcove $\waf_R\ac$.
Say $\A(\waf_R)=u\cdot\lambda+(h+1)\Q$, and $\Theta^{-1}(\waf_R)=[w,A]$.
By definition of $\Theta^{-1}$ we have $R\subseteq wC$ and for each $\alpha\in\Phi^+$ the hyperplane $H_{\alpha,1}$ is a floor of the dominant Shi alcove $w^{-1}\waf_R\ac$ if and only if $\alpha\in A$.

Choose $\mu,\nu\in\Q$ and $\tau,\sigma\in W$ such that $\waf_R^{-1}w=t_{\mu}\sigma$ and $\wf=t_{\nu}\tau$.

Then $\lambda=\wf\waf_R^{-1}\cdot0=\nu+\tau\mu$ and $\A(\waf_R)=-\waf_R\wf^{-1}\cdot0=w(\tau\sigma)^{-1}\cdot\lambda$.
It follows from the first identity that $\mu$ can be computed from $\lambda$ independently of $u$.
Since $\sigma$ is the unique element of $W$ such that $(t_{\mu}\sigma)^{-1}\ac$ lies in the dominant chamber, also $\sigma$ does not depend on $u$.

\begin{myprop}{antichain} We have
\begin{align*}
A=(\tau\sigma)^{-1}\cdot\walls
\end{align*}
\end{myprop}

\begin{myproof} Let $\beta\in\Phi^+$ then
\begin{align*}
\beta\in A
&\Leftrightarrow H_{\beta,1}\text{ is a floor of }w^{-1}\waf_R\ac\\
\text{(\refl{affroots})}&\Leftrightarrow \waf_R^{-1}w\cdot(-\beta+\delta)\in-\widetilde\Delta\\
\text{(\refl{wf})}&\Leftrightarrow \wf\waf_R^{-1}w\cdot(\beta-\delta)\in \wf(\widetilde\Delta)=(\Delta-\delta)\cup\{-\tilde\alpha+h\delta\}\\
&\Leftrightarrow t_{\lambda}\tau\sigma\cdot\beta\in\Delta\cup\{-\tilde\alpha+(h+1)\delta\}\\
&\Leftrightarrow
\begin{cases}
\tau\sigma\cdot\beta\in\Delta\text{ and }\skal{\lambda,\tau\sigma\cdot\beta}=0\text{, or}\\
\tau\sigma\cdot\beta=-\tilde\alpha\text{ and }\skal{\lambda,\tau\sigma\cdot\beta}=h+1.
\end{cases}
\end{align*}
Hence $\beta\in A$ if and only if $\tau\sigma\cdot\beta\in\walls$.
\end{myproof}

\sk
\refp{antichain} provides an explicit description of the map $\overline{\zeta}$ on $W$-orbits. The final result of this section explains how the Weyl group element $w$ can be obtained when $u$ is taken into account.

\begin{myprop}{utausig} We have $w=u\tau\sigma$.
\end{myprop}

\begin{myproof} Recall that
\begin{align*}
u\cdot\lambda
&=\A(\waf_R)
=-\waf_R\wf^{-1}\cdot0
=w(\tau\sigma)^{-1}\cdot(\nu+\tau\mu)
=w(\tau\sigma)^{-1}\cdot\lambda.
\end{align*}
Using \reft{stab} and \refp{antichain} it follows that
\begin{align*}
w^{-1}u\tau\sigma&\in(\tau\sigma)^{-1}\op{Stab}(\lambda)\tau\sigma\\
&\quad=\skal{s_{(\tau\sigma)^{-1}\cdot\alpha}:\alpha\in\walls}\\
&\quad=\skal{s_{\beta}:\beta\in A}.
\end{align*}
In other words $wH=u\tau\sigma H$ where $H=\skal{s_{\beta}:\beta\in A}$. On the one hand $w\cdot A\subseteq\Phi^+$ due to the fact that $\waf_R\ac$ is a Shi alcove contained in the chamber $wC$ \cite[Sec~12.3]{Thiel2015}.
On the other hand $u\tau\sigma\cdot A=u\cdot\walls\subseteq\Phi^+$ using \refp{antichain}.
Since there is only one element $v$ of $wH$ with $v\cdot A\subseteq\Phi^+$, we must have $w=u\tau\sigma$.
\end{myproof}

\subsection{Lattice paths}\label{Section:paths}

We call $N=(0,1)$ a \emph{North step} and $E=(1,0)$ an \emph{East step}. Given $a,b\in\N$, denote by $\mathcal L_{a,b}$ the set of lattice paths in the plane starting at the origin, consisting only of North and East steps and ending at $(a,b)$.
That is, $\mathcal L_{a,b}$ is the set of sequences $s_1s_2\dots s_{a+b}$ with $s_i\in\{N,E\}$ consisting of $a$ East steps and $b$ North steps.
A lattice path in $\mathcal L_{n,n}$ is called \emph{Dyck path} if it never goes below the main diagonal, that is, each initial segment $s_1s_2\cdots s_k$ contains at least as many North steps as East steps.
We denote by $\mathcal D_n\subseteq\mathcal L_{n,n}$ the set of all Dyck paths with $2n$ steps.
A ballot path of length $n$ is a lattice path starting at the origin, consisting of $n$ North and/or East steps that never goes below the main diagonal.
We denote the set of ballot paths with $n$ steps by $\mathcal B_n$.

A pattern of the form $NN$ is called \emph{rise}. A pattern $EN$ is called \emph{valley}. More precisely, let $\pi$ be any lattice path with steps $s_i\in\{N,E\}$. We say $i$ is a rise of $\pi$ if its $i$-th North step is followed by a North step. We say $(i,j)$ is a valley of $\pi$ if its $i$-th East step is followed by its $j$-th North step. If $\beta$ is a ballot path consisting of $i$ East steps and $j$ North steps, and $\beta$ ends with an East step, then we adopt the convention of counting $(i,j+1)$ as a valley of $\beta$. See Figures \ref{Figure:ballotpath} and \ref{Figure:dyckpath}.

\begin{figure}[ht]
\begin{minipage}[t]{.49\linewidth}
\begin{center}
\begin{tikzpicture}[scale=.5]
\draw[gray] (0,1)--(1,1)--(1,6)--(0,6)--(0,0)
(0,2)--(2,2)--(2,5)--(0,5)
(0,3)--(3,3)--(3,4)--(0,4);
\draw[very thick] (0,0)--(0,1)--(1,1)--(1,4)--(2,4);
	\draw[xshift=5mm,yshift=5mm]
		(0,1) node{\large{$\bullet$}}
		(1,4) node{\large{$\circ$}};
\end{tikzpicture}
\caption{A ballot path $\beta\in\mathcal B_6$ with two valleys $(1,2),(2,5)$ and two rises $2,3$.}
\label{Figure:ballotpath}
\end{center}
\end{minipage}
\begin{minipage}[t]{.48\linewidth}
\begin{center}
\begin{tikzpicture}[scale=.5]
\draw[gray] (0,0) grid (6,6)
	(0,0)--(6,6);
\draw[very thick] (0,0)--(0,4)--(2,4)--(2,5)--(3,5)--(3,6)--(6,6);
\draw (1.5,4.5) node{\large{$\bullet$}}
	(2.5,5.5) node{\large{$\bullet$}};
\end{tikzpicture}
\caption{A Dyck path $\pi\in\mathcal D_6$ with valleys $(2,5),(3,6)$ and rises $1,2,3$.}
\label{Figure:dyckpath}
\end{center}
\end{minipage}
\end{figure}

In order to define the combinatorial zeta maps, which are bijections between sets of (labelled) lattice paths, we need the following building blocks.
Let $a=(a_1,a_2,\dots,a_n)$ be a vector with integer entries.
Define a word $\RS_j^+(a)$ in the alphabet $\{E,N\}$ as follows:
Read $a$ from left to right.
Whenever you encounter an entry $a_i=j$ write down $N$.
Whenever you encounter an entry $a_i=j+1$ write down $E$.
Similarly we define $\LS_j^+(a)$ except now one reads $a$ from right to left.
Moreover we define $\RS_j^-(a)$ as follows:
Read $a$ from left to right and write down $N$ for each $-j$ you encounter and $E$ for each $-j-1$ you encounter.
Define $\LS_j^-(a)$ analogously.

We conclude this section by recalling the definitions of the two versions of the zeta map in type $A_{n-1}$, that is, one on labelled paths and one on unlabelled paths.
See \reff{zetaA} for an accompanying example.

Given a Dyck path $\pi\in\mathcal{D}_n$ we define its \emph{area vector} $(\mu_1,\dots,\mu_n)$ by letting $\mu_i$ be the number of $1\times1$-squares in the $i$-th row that lie between the path and the main diagonal.
For example the Dyck path in \reff{dyckpath} has area vector $(0,1,2,3,2,2)$.
It is not hard to see that a vector with non-negative integer entries is the area vector of a Dyck path if and only if $\mu_1=0$ and $\mu_{i+1}\leq\mu_{i}+1$.

The original \emph{zeta map} $\zeta_A:\mathcal{D}_n\to\mathcal{D}_n$ first appears in a paper of Andrews, Krattenthaler, Orsina and Papi~\cite{AKOP2002}.
It was rediscovered and popularised by Haglund and Haiman, and an explicit treatment of its compatibility with the statistics on Dyck paths mentioned in the introduction is found in~\cite{Haglund2008}.
The image of the Dyck path $\pi$ under map $\zeta_A$ can be defined as the concatenation
\begin{align*}
\zeta_A(\pi)&=\RS_0^-(\mu)\RS_{-1}^-(\mu)\cdots\RS_{-n}^-(\mu)
\end{align*}
where $\mu$ is the area vector of $\pi$.

A \emph{vertically labelled Dyck path} is a pair $(\pi,v)$ of a Dyck path $\pi\in\mathcal{D}_n$ and a permutation $v\in\S_n$ such that $v(i)<v(i+1)$ for each rise $i$ of $\pi$.
Denote the set of all vertically labelled Dyck paths by $\operatorname{Vert}(A_{n-1})$.

The \emph{diagonal reading word} $d_A(\pi,v)$ of a vertically labelled Dyck path is the permutation in $\S_n$ obtained as follows:
Let $(\mu_1,\dots,\mu_n)$ be the area vector of $\pi$.
For each $i=0,1,\dots,n-1$ read $\mu$ from left to right and write down $v(j)$ for each encountered entry $\mu_j=i$.

A \emph{diagonally labelled Dyck path} is a pair $(\pi,w)$ of a Dyck path $\pi\in\mathcal{D}_n$ and a permutation $w\in\S_n$ such that $w(i)<w(j)$ for each valley $(i,j)$ of $\pi$.
Let $\operatorname{Diag}(A_{n-1})$ denote the set of all diagonally labelled Dyck paths.

A generalisation of the original zeta map, which we call \emph{Haglund--Loehr-zeta map}, appears in~\cite{HagLoe2005} and sends vertically labelled Dyck paths to diagonally labelled Dyck paths.
It can be defined as
\begin{align*}
\zeta_A:\operatorname{Vert}(A_{n-1})&\to\operatorname{Diag}(A_{n-1}),\\
(\pi,v)&\mapsto(\zeta_A(\pi),d_A(\pi,v)).
\end{align*}
Note that the rises of $(\pi,v)$ correspond precisely to the valleys of $\zeta_A(\pi,v)$ in the following sense:
For each rise $i$ of $\pi$ there exists a valley $(j,k)$ of $\zeta_A(\pi)$ such that $(v(i),v(i+1))=(d_A(\pi,v)(j),d_A(\pi,v)(k))$ and vice versa \cite[Section 5.2]{ALW2014} \cite[Theorem 14.1]{Thiel2015}.
\begin{figure}[t]
\begin{center}
\begin{tikzpicture}[scale=.5]
\draw[gray]
(0,0)grid(6,6)
;
\draw[very thick]
(0,0)--(0,4)--(2,4)--(2,5)--(3,5)--(3,6)--(6,6)
;
\draw[xshift=5mm,yshift=5mm]
(0,0)node{\large{$2$}}
(0,1)node{\large{$3$}}
(0,2)node{\large{$4$}}
(0,3)node{\large{$6$}}
(2,4)node{\large{$1$}}
(3,5)node{\large{$5$}}
;
\draw[xshift=8cm,yshift=25mm,thick,->]
(0,0)--node[above]{$\zeta_A$}(1,0)
;
\begin{scope}[xshift=11cm]
\draw[gray]
(0,0)grid(6,6)
;
\draw[very thick]
(0,0)--(0,1)--(1,1)--(1,2)--(2,2)--(2,5)--(3,5)--(3,6)--(6,6)
;
\draw[xshift=5mm,yshift=5mm]
(0,0)node{\large{$2$}}
(1,1)node{\large{$3$}}
(2,2)node{\large{$4$}}
(3,3)node{\large{$1$}}
(4,4)node{\large{$5$}}
(5,5)node{\large{$6$}}
(0,1)node{\large{$\bullet$}}
(1,2)node{\large{$\bullet$}}
(2,5)node{\large{$\bullet$}}
;
\end{scope}
\end{tikzpicture}
\caption{A vertically labelled Dyck path $(\pi,v)$ and the diagonally labelled Dyck path $\zeta_A(\pi,v)$.}
\label{Figure:zetaA}
\end{center}
\end{figure}
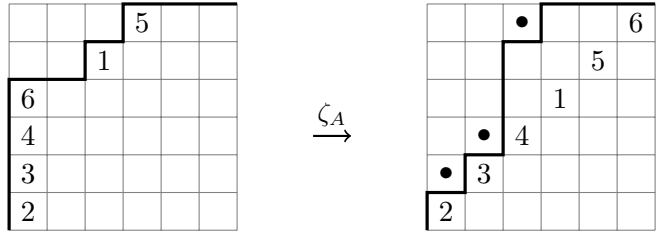

%% file: preps.tex
\section{Types \emph{B}, \emph{C} and \emph{D}}\label{Section:BCD}

In this section we fix further notation and recall some facts specific to the root systems and Weyl groups of types $B_n, C_n$ and $D_n$.

\subsection{The root systems of types \emph{B}, \emph{C} and \emph{D}}

The roots, positive roots and simple roots of type $B_n$ are given by
\begin{align*}
\Phi&=\{\pm e_i\pm e_j:1\leq i<j\leq n\}\cup\{\pm e_i:i\in[n]\}, \\
\Phi^+&=\{\pm e_i+e_j:1\leq i<j\leq n\}\cup\{e_i:i\in[n]\}, \text{and} \\
\Delta&=\{-e_i+e_{i+1}:i\in[n-1]\}\cup\{e_1\}.
\end{align*}
We denote the simple roots by $\alpha_0^B=e_1$ and $\alpha_i^B=e_{i+1}-e_i$ for $i\in[n-1]$, and the highest root by $\tilde\alpha^B=e_{n-1}+e_n$.
The root lattice and coroot lattice are
\begin{align*}
Q
&=\bigoplus_{\alpha\in\Delta}\Z\alpha
=\Z^n,&&
\Q
=\bigoplus_{\alpha\in\Delta}\check\Z\alpha
=\Big\{(x_1,\dots,x_n)\in\Z^n:\sum_i{x_i}\in2\Z\Big\}.
\end{align*}
The roots, positive roots and simple roots of type $C_n$ are given by
\begin{align*}
\Phi&=\{\pm e_i\pm e_j:1\leq i<j\leq n\}\cup\{\pm2e_i:i\in[n]\}, \\
\Phi^+&=\{\pm e_i+e_j:1\leq i<j\leq n\}\cup\{2e_i:i\in[n]\}, \text{and} \\
\Delta&=\{-e_i+e_{i+1}:i\in[n-1]\}\cup\{2e_1\}.
\end{align*}
We denote the simple roots by $\alpha_0^C=2e_1$ and $\alpha_i^C=e_{i+1}-e_i$ for $i\in[n-1]$, and the highest root by $\tilde\alpha^C=2e_n$.
The root lattice and coroot lattice are dual to the type $B_n$ case
\begin{align*}
Q
=\Big\{(x_1,\dots,x_n)\in\Z^n:\sum_i{x_i}\in2\Z\Big\},&&
\Q
=\Z^n.
\end{align*}
The roots, positive roots and simple roots of type $D_n$ are given by
\begin{align*}
\Phi&=\{\pm e_i\pm e_j:1\leq i<j\leq n\}, \\
\Phi^+&=\{\pm e_i+e_j:1\leq i<j\leq n\}, \text{and} \\
\Delta&=\{-e_i+e_{i+1}:i\in[n-1]\}\cup\{e_1+e_2\}.
\end{align*}
We denote the simple roots by $\alpha_0^D=e_1+e_2$ and $\alpha_i^D=e_{i+1}-e_i$ for $i\in[n-1]$, and the highest root by $\tilde\alpha^D=e_{n-1}+e_n$.
The root lattice and coroot lattice are both isomorphic to the coroot lattice of type $B_n$.
\begin{align*}
Q
&=\Q
=\Big\{(x_1,\dots,x_n)\in\Z^n:\sum{x_i}\in2\Z\Big\}
\end{align*}

\subsection{Affine permutation groups}

The affine Weyl groups of types $A_{n-1}$, $B_n$, $C_n$ and $D_n$ can all be realised as groups of certain bijections on integers also called affine permutations. In this section we recall the descriptions of the affine Weyl groups of types $B_n$, $C_n$ and $D_n$ in terms of affine permutations given in \cite{BjoBre2005}.

Set $\n=2n+1$. A bijection $\omega:\Z\to\Z$ is called \emph{affine permutation} if $\omega(i+\n)=\omega(i)+\n$ and $\omega(-i)=-\omega(i)$ for all $i\in\Z$. The set of all such bijections forms a group under composition which we denote by $\affS_n^C$. Each affine permutation is fully determined by its \emph{window}
\begin{align*}
[\omega(1),\omega(2),\dots,\omega(n)].
\end{align*}
The group $\affS_n^C$ is generated by the following $n+1$ \emph{simple transpositions of type $C_n$}, defined as
\begin{align*}
s_0^C&=[-1,2,\dots,n],\\
s_i^C&=[1,\dots,i-1,i+1,i,i+2,\dots,n]&&\text{for }i\in[n-1]\text{ and}\\
s_n^C&=[1,\dots,n-1,n+1].
\end{align*}
The group $\affS_n^C$ contains a subgroup $\affS_n^B$ that consists of all affine permutations $\omega$ such that the finite set $\{i\in\Z:i\leq n,\omega(i)>n\}$ has even cardinality. This subgroup is generated by the \emph{simple transpositions of type $B_n$}, given by
\begin{align*}
s_0^B&=[-1,2,\dots,n],\\
s_i^B&=[1,\dots,i-1,i+1,i,i+2,\dots,n]&&\text{for }i\in[n-1]\text{ and}\\
s_n^B&=[1,\dots,n-2,n+1,n+2].
\end{align*}
The group $\affS_n^B$ contains a subgroup $\affS_n^D$ that consists of all affine permutations $\omega$ such that both finite sets $\{i\in\Z:i\leq n,\omega(i)>n\}$ and $\{i\in\Z:i\geq0,\omega(i)<0\}$ have even cardinality. This subgroup is generated by the \emph{simple transpositions of type $D_n$}, that is,
\begin{align*}
s_0^D&=[-1,-2,3,\dots,n],\\
s_i^D&=[1,\dots,i-1,i+1,i,i+2,\dots,n]&&\text{for }i\in[n-1]\text{ and}\\
s_n^D&=[1,\dots,n-2,n+1,n+2].
\end{align*}
Let $\Phi$ be a root system of type $B_n$, $C_n$ or $D_n$. The group $\affS_n^{\Phi}$ is isomorphic to the affine Weyl group $\wa$ of $\Phi$. An explicit isomorphism in terms of the generators is obtained by mapping $s_i^{\Phi}$ to $s_{\alpha_i^{\Phi}}$ for $0\leq i\leq n-1$, and $s_n^{\Phi}$ to $s_{\tilde\alpha^{\Phi},1}$. Let $\S_n^{\Phi}$ denote the subgroup of $\affS_n^{\Phi}$ corresponding to the Weyl group $W$ under this isomorphism. Then an affine permutation $\omega\in\affS_n^{\Phi}$ lies in $\S_n^{\Phi}$ if and only if its window is a subset of $[\pm n]$. More precisely, $\S_n^B=\S_n^C$ is the group of signed permutations while $\S_n^D\leq\S_n^B$ consists of the signed permutations with an even number of sign changes. Furthermore, this isomorphism affords an action of the affine permutations on the coroot lattice $\Q$ of $\Phi$, which is made explicit by the following rules
\begin{align*}
s_0^C\cdot(x_1\dots,x_n)&=(-x_1,x_2\dots,x_n),\\
s_i^C\cdot(x_1\dots,x_n)&=(x_1,\dots,x_{i+1},x_i,\dots,x_n),\quad\text{for }i\in[n-1],\\
s_n^C\cdot(x_1\dots,x_n)&=(x_1,\dots,x_{n-1},-x_n+1),\\
s_n^B\cdot(x_1\dots,x_n)&=(x_1,\dots,x_{n-2},-x_n+1,-x_{n-1}+1)\quad\text{and}\\
s_0^D\cdot(x_1\dots,x_n)&=(-x_1,-x_2,x_3,\dots,x_n).
\end{align*}
The first tool we need is a combinatorial description of the Gra{\ss}mannian affine permutations in $\affS_n^{\Phi}$. That is, given an affine permutation $\omega\in\affS_n^{\Phi}$ we want to decide whether $\omega$ is a minimal length coset representative just by looking at its window. This question was answered in \cite[Prop~8.4.4,~8.5.4, and~8.6.4]{BjoBre2005}.

\begin{myprop}{grassmannian}
An affine permutation $\omega\in\affS_n^{\Phi}$ is the minimal length representative of the coset $\omega\S_n^{\Phi}\in\affS_n^{\Phi}/\S_n^{\Phi}$ if and only if
\begin{align*}
\begin{cases}
0<\omega(1)<\omega(2)<\dots<\omega(n)&\quad\text{if }\Phi\text{ is of type }B_n\text{ or }C_n,\\
0<\abs{\omega(1)}<\omega(2)<\dots<\omega(n)&\quad\text{if }\Phi\text{ is of type }D_n.
\end{cases}
\end{align*}
\end{myprop}

\sk
Next we want to obtain a combinatorial description of the decomposition of $\omega$ into a product of a translation by an element of the coroot lattice and an element of the Weyl group, that is, a signed permutation. To this end, for $\omega\in\affS_n^{\Phi}$ and each $i\in\Z$ write $\omega(i)=a_i\n+b_i$ such that $a_i\in\Z$ and $b_i\in[\pm n]$. Define $\sigma(\omega,i)=b_i$, $\mu(\omega,b_i)=-a_i$ and $\nu(\omega,i)=a_i$.

\begin{mylem}{munusig} Let $\omega\in\affS_n^{\Phi}$ be an affine permutation.
\myi{(i)} The map $i\mapsto\sigma(\omega,i)$, where $i\in[n]$, defines a signed permutation $\sigma(\omega)\in\S_n^{\Phi}$.
\myi{(ii)} The vectors $\mu(\omega)=(\mu(\omega,1),\dots,\mu(\omega,n))$ and $\nu(\omega)=(\nu(\omega,1),\dots,\nu(\omega,n))$ lie in the coroot lattice $\Q$ of $\Phi$, and for all $i$ with $0\leq i\leq n$ we have
\begin{align*}
\mu(s_i^{\Phi}\omega)=s_i^{\Phi}\cdot\mu(\omega),&&
\nu(\omega s_i^{\Phi})=s_i^{\Phi}\cdot\nu(\omega).
\end{align*}
\myi{(iii)} We have $\mu(\omega^{-1})=\nu(\omega)$ and $\sigma(\omega^{-1})=\sigma(\omega)^{-1}$.
\myi{(iv)} We have $\omega\cdot(0,\dots,0)=\mu(\omega)$ and $\omega^{-1}\cdot(0,\dots,0)=\nu(\omega)$
\myi{(v)} We have $\mu(\omega)=-\sigma\cdot\nu(\omega)$.
\end{mylem}

\begin{myproof} Claims (i) and (v) are immediate from the definitions while (iii) and (iv) follow directly from (ii). Thus it only remains to show (ii) which is done for each type using induction on the length of $\omega$.
\end{myproof}

\sk
For $q\in\Q$ define an affine permutation $t_q\in\affS_n^{\Phi}$ by setting $t_q(i)=-q_i\n+i$ for $i\in[n]$. Set $T_{\Q}=\{t_q:q\in\Q\}$. We call an affine permutation $\omega\in\affS_n^{\Phi}$ translation by $q\in\Q$ if $\omega\cdot x=x+q$ for all $x\in\Q$. Thus by definition the translations in $\affS_n^{\Phi}$ correspond to translations in $\wa$.

\begin{myprop}[ (i)]{translations} Let $\omega\in\affS_n^{\Phi}$ be an affine permutation, and set $\sigma=\sigma(\omega)$, $\mu=\mu(\omega)$ and $\nu=\nu(\omega)$. Then $\omega=t_{\mu}\sigma=\sigma t_{-\nu}$.
\myi{(ii)} Let $x,y\in\Q$ then $t_xt_y=t_{x+y}$ and $(t_x)^{-1}=t_{-x}$. Thus $T_{\Q}$ is a subgroup of $\affS_n^{\Phi}$ that is isomorphic to the coroot lattice.
\myi{(iii)} An affine permutation $\omega\in\affS_n^{\Phi}$ is a translation if and only if $\omega\in T_{\Q}$.
\end{myprop}

\begin{myproof} To prove claim (i) let $i\in[n]$ and $\omega(i)=a_i\n+b_i$ such that $a_i\in\Z$ and $b_i\in[\pm n]$. Then
\begin{align*}
t_{\mu}\sigma(i)
&=t_{\mu}(\sigma(\omega,i))
=t_{\mu}(b_i)
=-\mu(\omega,b_i)\n+b_i
=a_i\n+b_i,\\
\sigma t_{-\nu}(i)
&=\sigma(\nu(\omega,i)\n+i)
=\sigma(a_i\n+i)
=a_i\n+\sigma(i)
=a_i\n+b_i.
\end{align*}
To see claim (ii), note that for each $i\in[n]$ we have
\begin{align*}
t_xt_y(i)
&=t_x(-y_i\n+i)
=-y_i\n+t_x(i)
=-(y_i+x_i)\n+i
=t_{x+y}(i).
\end{align*}
Finally, let $q,x\in\Q$. Then using~(i),~(ii) and \refl{munusig}~(iv) we see that $t_q$ is indeed a translation as
\begin{align*}
t_q\cdot x
&=t_q\cdot t_x\cdot(0,\dots,0)
=(t_qt_x)\cdot(0,\dots,0)
=t_{q+x}\cdot(0,\dots,0)
=q+x.
\end{align*}
Conversely, if $\omega$ is a translation by $q$ then $\omega=t_{q}$ because $\affS_n^{\Phi}$ acts faithfully on the coroot lattice.

\end{myproof}

\sk
In fact from \refp{translations}~(i) and~(ii) it follows directly that $\affS_n^{\Phi}$ is the semi direct product of $T_{\Q}$ and $\S_n^{\Phi}$ (the Weyl group acting on the coroot lattice). Hence the decomposition $\omega=t_{\mu}\sigma$ is the one we were looking for, and it can easily be obtained from the window of $\omega$.

The following lemma is an easy consequence of the two propositions above but it shall serve us as a reference in the sections to come.

\begin{mylem}{sigma} Let $\omega\in\affS_n^{\Phi}$ be the minimal length coset representative of $\omega\S_n^{\Phi}$, $\mu=\mu(\omega)\in\Q$ and $\sigma=\sigma(\omega)\in\S_n^{\Phi}$. Then for each $i\in[n]$
\begin{align*}
\abs{\sigma^{-1}(i)}
&=\#\big\{k\in[n]:\abs{\mu_k\n-k}\leq\abs{\mu_i\n-i}\big\}.
\end{align*}
If $\Phi$ is of type $B_n$ or $C_n$, or if $\Phi$ is of type $D_n$ and $\abs{\sigma^{-1}(i)}\neq1$, then $\sigma^{-1}(i)>0$ if and only if $\mu_i\leq0$. 
However, if $\Phi$ is of type $D_n$ and $\abs{\sigma^{-1}(i)}=1$, then $\sigma^{-1}(i)>0$ if and only if either $\mu_i\leq0$ and the number of positive entries of $\mu$ is even, or $\mu_i>0$ and the number of positive entries of $\mu$ is odd.
\end{mylem}

\begin{myproof} Since $\omega$ is a minimal length coset representative, the absolute values of the entries of the window of $\omega$ must be increasing, that is, $0<\abs{\omega(1)}<\abs{\omega(2)}<\dots<\abs{\omega(n)}$. On the other hand, $\abs{\sigma(j)}=i$ if and only if $\omega(j)=\abs{\mu_i\n-i}$. Hence $\abs{\sigma(j)}=i$ is equivalent to
\begin{align*}
j=\#\big\{k\in[n]:\abs{\mu_k\n-k}\leq\abs{\mu_i\n-i}\big\}.
\end{align*}
Furthermore, if $\Phi$ is of type $B_n$ or $C_n$, or if $\Phi$ is of type $D_n$ and $\abs{\sigma^{-1}(i)}\neq1$, then $\omega(i)>0$. Hence $\sigma(i)>0$ if and only if $-\n\mu_i+i>0$, which is the case if and only if $\mu_i\leq0$.

If $\Phi$ is of type $D_n$ then the sign of $\omega(1)$ possibly has to be changed such that there is an even number of integers $j\in\Z$ with $j\geq0$ and $\omega(j)<0$.
\end{myproof}

%% file: typeC.tex
\section{The zeta map of type C}\label{Section:typeC}

In \refs{vertC} we define a combinatorial model for the elements of the finite torus in terms of vertically labelled lattice paths in an $n\times n$-square.
The underlying lattice paths correspond to the orbits of the finite torus under the action of the Weyl group $\S_n^C$.
That is, acting by an element of $\S_n^C$ only changes the labels but leaves the path intact.
In \refs{diagC} we encode the non-nesting parking functions in term of diagonally labelled ballot paths.
Again acting by an element of the Weyl group only changes the labels, such that the ballot paths correspond to the antichains in the root poset.

We then demonstrate how the dominant Shi region corresponding to a square lattice path can be recovered using its area vector in \refs{areaC}.
Moreover we obtain the Shi region corresponding to a vertically labelled lattice path using the diagonal reading word in \refs{readC}.
Finally, in \refs{zetaC} we provide a combinatorial construction of the image of a vertically labelled lattice path under the zeta map of type $C_n$.
We prove that this construction is invertible and that it coincides with the uniform zeta map.

\refs{sweepC} offers an alternative description of the type $C_n$ zeta map in the spirit of the sweep map.
We conclude by defining type $C_n$ analogues of the $\area$ and $\dinv$ statistics, and discuss their connection via the zeta map in \refs{statC}


\subsection{The finite torus}\label{Section:vertC}

The Coxeter number of type $C_n$ is given by $h=2n$ hence the finite torus equals $\Z^n/(2n+1)\Z^n$. As Athanasiadis~\cite[Sec~5.2]{Athan2005} pointed out, a system of representatives for the orbits of the finite torus under the action of the Weyl group $\S_n^C$ is given by
\begin{align*}
\Q\cap (2n+1)\overline{\ac}
=\Big\{(\lambda_1,\lambda_2,\dots,\lambda_n)\in\Q: 0\leq\lambda_1\leq\lambda_2\leq\dots\leq\lambda_n\leq n\Big\},
\end{align*}
and the stabiliser subgroup $\op{Stab}(\lambda)\leq\S_n^C$ of such a $\lambda$ is generated by the simple transpositions $s_i$ for each $i\in[n-1]$ such that $\lambda_i=\lambda_{i+1}$, and the simple transposition $s_0$ if $\lambda_1=0$.


\begin{mydef}{vertC} A \emph{vertically labelled lattice path} is a pair $(\pi,u)$ of a lattice path $\pi\in\mathcal L_{n,n}$ and a signed permutation $u\in\S_n^C$ such that $u(i)<u(i+1)$ for each rise $i$ of $\pi$ and such that $u(1)>0$ if $\pi$ begins with a North step. We denote the set of all vertically labelled lattice paths by $\op{Vert}(C_n)$.
\end{mydef}
\sk
Given a lattice path $\pi\in\mathcal L_{n,n}$ and a signed permutation $u\in\S_n^C$ we picture $(\pi,u)$ by placing the label $\sigma(i)$ to the left of the $i$-th North step of $\pi$ as is shown in \reff{vertC}. Thus we obtain $(\pi,u)\in\op{Vert}(C_n)$ if the labels increase along columns from bottom to top, and if all labels left of the starting point of $\pi$ are positive.

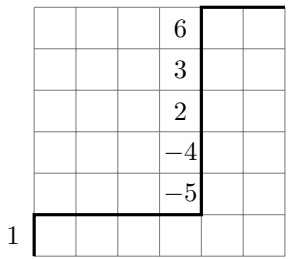
\begin{figure}[ht]
\begin{center}
\begin{tikzpicture}[scale=.55]
	\draw[gray] (0,0) grid (6,6);
	\draw[very thick] (0,0)--(0,1)--(4,1)--(4,6)--(6,6);
	\draw[xshift=5mm,yshift=5mm]
		(-1,0) node{$1$}
		(3,1) node{$-5$}
		(3,2) node{$-4$}
		(3,3) node{$2$}
		(3,4) node{$3$}
		(3,5) node{$6$};
\end{tikzpicture}
\caption{The vertically labelled lattice path $(NEEEENNNNNEE,[1,-5,-4,2,3,6])\in\op{Vert}(C_6)$.}
\label{Figure:vertC}
\end{center}
\end{figure}

It is not difficult to see that vertically labelled lattice paths encode the elements of the finite torus.

\begin{myprop}{vertC} The map $\psi:\op{Vert}(C_n)\to\Q/(2n+1)\Q$ given by $(\pi,u)\mapsto u\cdot\lambda+(2n+1)\Q$, where $\lambda\subseteq(n^n)$ is the partition\footnote{For our purposes here a partition is an increasing sequence $\lambda_1\leq\lambda_2\leq\dots\leq\lambda_n$ of non-negative integers.} with South-East boundary $\pi$, is a bijection.
\end{myprop}

\begin{myproof} First note that the set of partitions that fit inside the square $(n^n)$ coincides with the set $\Q\cap(2n+1)\overline{\ac}$.
Let $\pi\in\mathcal L_{n,n}$ be a lattice path, $u\in\S_n^C$ be a signed permutation, and $\lambda\subseteq(n^n)$ be the partition with South-East boundary $\pi$.
By \refl{canon} it suffices to show that $(\pi,u)\in\op{Vert}(C_n)$ if and only if $u\cdot\walls\subseteq\Phi^+$.
Suppose $i\in[n-1]$ and $(\pi,u)\in\op{Vert}(C_n)$ then
\begin{align*}
\alpha_i\in\walls
&\Leftrightarrow \skal{\lambda,e_{i+1}-e_i}=0
\Leftrightarrow \lambda_i=\lambda_{i+1}\\
&\Rightarrow u(i)<u(i+1)
\Leftrightarrow u\cdot\alpha_i\in\Phi^+,\\
\alpha_0\in\walls
&\Leftrightarrow \skal{\lambda,2e_1}=0
\Leftrightarrow \lambda_1=0\\
&\Rightarrow 0<u(1)
\Leftrightarrow u\cdot\alpha_0\in\Phi^+.
\end{align*}
Conversely if $i\in[n-1]$ and $u\cdot\walls\subseteq\Phi^+$ then
\begin{align*}
\lambda_i=\lambda_{i+1}
&\Leftrightarrow\alpha_i\in\walls\\
&\Rightarrow u\cdot\alpha_i\in\Phi^+
\Leftrightarrow u(i)<u(i+1),\\
\lambda_1=0
&\Leftrightarrow\alpha_0\in\walls\\
&\Rightarrow u\cdot\alpha_0\in\Phi^+
\Leftrightarrow 0<u(1).
\end{align*}
To complete the proof, note that $\skal{\lambda,\tilde\alpha}\leq2n$ implies $\tilde\alpha\notin\walls$.
\end{myproof}

\begin{myex}{vertC} The vertically labelled lattice path in \reff{vertC} corresponds to the partition $\lambda=(0,4,4,4,4,4)$ and the element $(0,4,4,-4,-4,4)+13\Z^6$ of the finite torus.
\end{myex}

\subsection{Non-nesting parking functions}\label{Section:diagC} In this section we encode the non-nesting parking functions of type $C_n$ in terms of diagonally labelled ballot paths.

\begin{mydef}{DiagC} A \emph{diagonally labelled ballot path} is a pair $(\beta,w)$ of a ballot path $\beta\in\mathcal B_{2n}$ and a signed permutation $w\in\S_n^C$ such that for each valley $(i,j)$ of $\beta$ we have
\begin{align*}
w(n+1-i)>
\begin{cases}
w(n+1-j)&\quad\text{if }j\leq n,\\
w(j-n)&\quad\text{if }j>n.
\end{cases}
\end{align*}
We denote the set of all diagonally labelled ballot paths by $\op{Diag}(C_n)$.
\end{mydef}
\sk
Suppose $(\beta,w)\in\op{Diag}(C_n)$. If we place the labels $w(i)$, where $i=n,n-1,\dots,1,-1,\dots,-n$, in the diagonal as in \reff{diagC}, then for each valley the label to its right will be smaller than the label below it. In particular, if the path ends with an East step then the label below will be positive.

\begin{figure}[t]
\begin{center}
\begin{tikzpicture}[scale=.55]
	\ballot{6}
	\draw[very thick] (0,0)--(0,2)--(1,2)--(1,3)--(2,3)--(2,5)--(3,5)--(3,6)--(4,6)--(4,7)--(5,7);
	\draw[xshift=5mm,yshift=5mm]
		(0,0) node{$5$}
		(1,1) node{$6$}
		(2,2) node{$4$}
		(3,3) node{$3$}
		(4,4) node{$1$}
		(5,5) node{$-2$}
		(6,6) node{$2$}
		(5,7) node{$-1$}
		(4,8) node{$-3$}
		(3,9) node{$-4$}
		(2,10) node{$-6$}
		(1,11) node{$-5$}
		(0,2) node{\large{$\bullet$}}
		(1,3) node{\large{$\bullet$}}
		(2,5) node{\large{$\bullet$}}
		(3,6) node{\large{$\bullet$}}
		(4,7) node{\large{$\circ$}};
\end{tikzpicture}
\caption{The diagonally labelled ballot path $(NNENENNENENE,[-2,1,3,4,6,5])\in\op{Diag}(C_6)$.}
\label{Figure:diagC}
\end{center}
\end{figure}
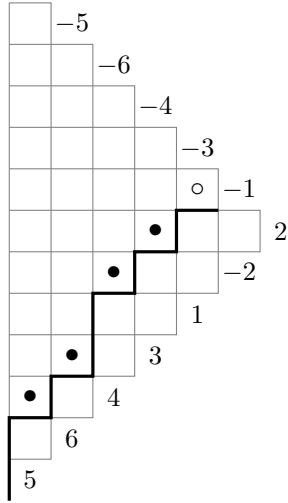

\sk
Let $\beta\in\mathcal B_{2n}$ be a ballot path and $(i,j)$ a valley of $\beta$. We define the positive root
\begin{align*}
\alpha_{i,j}&=
\begin{cases}
e_{n+1-i}-e_{n+1-j}&\quad\text{if }j\leq n,\\
e_{n+1-i}+e_{j-n}&\quad\text{if }j>n.
\end{cases}
\end{align*}
Furthermore set
\begin{align*}
A_{\beta}=\big\{\alpha_{i,j}:(i,j)\text{ is a valley of }\beta\big\}.
\end{align*}
We will use the following well-known result.

\begin{mylem}{AbetaC} The map $\varphi:\mathcal B_{2n}\to\op{Antichains}(\Phi^+)$ given by $\beta\mapsto A_{\beta}$ is a bijection between ballot paths of length $2n$ and the set of antichains in the root poset $\Phi^+$ of type $C_n$.
\end{mylem}

\begin{myex}{diagC} The diagonally labelled ballot path $(\beta,w)$ in \reff{diagC} has five valleys: $(1,3)$, $(2,4),(3,6),(4,7)$ and $(5,8)$. Note that for each valley the number to its right is less than the number below (for example $w(6+1-1)=5>4=w(6+1-3)$) and that the label below the terminal East step is positive: $w(6+1-5)=1>0$. The corresponding antichain $A_{\beta}$ contains the roots $e_6-e_4,e_5-e_3,e_4-e_1,e_3+e_1$ and $2e_2$.
\end{myex}
\sk
Taking the labels into account, we extend \refl{AbetaC} to a bijection between diagonally labelled ballot paths and non-nesting parking functions.

\begin{myprop}{diagC} The map $\varphi:\op{Diag}(C_n)\to\op{Park}(C_n)$ given by $(\beta,w)\mapsto[w,A_{\beta}]$ is a bijection.
\end{myprop}

\begin{myproof} Let $\beta\in\mathcal B_{2n}$ be a ballot path and $w\in\S_n^C$ a signed permutation.
Due to \refl{AbetaC} and \refl{nnpark} it suffices to show that $(\beta,w)\in\op{Diag}(C_n)$ if and only if $w\cdot A_{\beta}\subseteq\Phi^+$. Assume $j\leq n$. If $(\beta,w)\in\op{Diag}(C_n)$ then
\begin{align*}
\alpha_{i,j}\in A_{\beta}
&\Leftrightarrow (i,j)\text{ is a valley of }\beta\\
&\Rightarrow w(n+1-i)>w(n+1-j)
\Leftrightarrow w\cdot\alpha_{i,j}
=w\cdot(e_{n+1-i}-e_{n+1-j})\in\Phi^+.
\end{align*}
Conversely, if $w\cdot A_{\beta}\subseteq\Phi^+$ then
\begin{align*}
(i,j)\text{ is a valley of }\beta
&\Leftrightarrow \alpha_{i,j}\in A_{\beta}\\
&\Rightarrow w\cdot\alpha_{i,j}\in\Phi^+
\Leftrightarrow w(n+1-i)>w(n+1-j).
\end{align*}
The case $j>n$ is treated similarly.
\end{myproof}

\subsection{The area vector}\label{Section:areaC}

Let $\pi\in\mathcal L_{n,n}$ be a lattice path and $\lambda\subseteq(n^n)$ the partition with South-East boundary $\pi$.
We have seen above that $\lambda\in\Q\cap(2n+1)\overline{\ac}$.
Thus $\lambda$ corresponds to a dominant Shi region with minimal alcove $\waf_D\ac$ by means of the Anderson map of \refs{torus}.
Our goal in this section is to obtain the element of the affine Weyl group $\waf_D$ from the lattice path $\pi$.

To this end and write $\waf_D^{-1}=t_{\mu}\sigma$ where $\mu\in\Q$ and $\sigma\in\S_n^C$. Since $\waf_D^{-1}$ is Gra{\ss}mannian, the signed permutation $\sigma$ is determined by \refl{sigma} once we know $\mu$.
Note that $\waf_D$ corresponds to the affine transformation $w^{-1}\waf_R$ from \refs{uniform}. The vector $\mu$ therefore appears implicitly in the identity $\lambda=\waf_f\waf_D^{-1}\cdot0$. Hence all we need is an explicit description of $\waf_f$. This is provided by the next lemma.

\begin{mylem}{wfC} Let $\nu\in\Q$ and $\tau\in W$ be such that $\waf_f=t_{\nu}\tau$. Then
\begin{align*}
\nu=
(1,2,\dots,n),&&
\tau=
[-n,\dots,-2,-1].
\end{align*}
\end{mylem}
\begin{myproof} By \refl{wf} it suffices to show that $t_{\nu}\tau(\widetilde\Delta)=(\Delta-\delta)\cup\{-\tilde\alpha+h\delta\}$. We have
\begin{align*}
t_{\nu}\tau(\alpha_i)
&=\alpha_{n+1-i}-\skal{\nu,\alpha_{n+1-i}}\delta
=\alpha_{n+1-i}-\delta
&&\text{for all }i\in[n-1],\\
t_{\nu}\tau(\alpha_{0})
&=-2e_n-\skal{\tau,-2e_n}\delta
=-\tilde\alpha+2n\delta
&&\text{and}\\
t_{\nu}\tau(-\tilde\alpha+\delta)
&=2e_1+(1-\skal{\tau,2e_1})\delta
=\alpha_0-\delta.
\end{align*}
%
\end{myproof}

\sk
For example if $n=5$ then
\begin{align*}
\waf_f=[50,40,30,20,10].
\end{align*}
\begin{mydef}{areaC} Let $\pi\in\mathcal L_{n,n}$ and $\lambda\subseteq(n^n)$ be the partition with South-East boundary $\pi$. Moreover define $\nu$ and $\tau$ as in \refl{wfC} above. We define the type $C_n$ \emph{area vector} of $\pi$ as
\begin{align*}
\mu
=\tau(\lambda-\nu)
=\tau(\lambda_1-1,\lambda_2-2\dots,\lambda_n-n)
=(n-\lambda_n,\dots,2-\lambda_2,1-\lambda_1).
\end{align*}
\end{mydef}

Indeed $\mu_i$ counts the number of boxes in the $(n+1-i)$-th row between the path $\pi$ and the path $(EN)^n\in\mathcal L_{n,n}$ consisting of alternating North and East steps. In this regard $\mu$ is quite similar to the type $A_{n-1}$ area vector of a Dyck path. Note that the entries of the type $C_n$ area vector are negative as long as $\pi$ is to the right of $(EN)^n$ in the respective row. See \reff{areaC}.

\begin{myex}{areaC} The lattice path $\pi=NEEEENNNNNEE$ in \reff{vertC} has area vector
\begin{align*}
\mu
=(6,5,4,3,2,1)-(4,4,4,4,4,0)
=(2,1,0,-1,-2,1).
\end{align*}
The translation by $\mu$ is given by $t_{\mu}=[-25,-11,3,17,31,-7]\in\affS_n^C$. Moreover $\sigma=[3,-6,-2,4,-1,5]$ is the unique element of $\S_n^C$ such that $\waf_D^{-1}=t_{\mu}\sigma=[3,7,11,17,25,31]$ is a Gra{\ss}mannian affine permutation. Finally $\waf_D=[21,10,1,-9,-20,11]$ and $\waf_D\ac$ is the minimal alcove of the dominant Shi region corresponding to the lattice path $\pi$.
\end{myex}

\begin{figure}[ht]
\begin{center}
\begin{tikzpicture}[scale=.55]
\draw[gray] (0,0) grid (4,4);
\draw[very thick] (0,0)--(1,0)--(1,1)--(2,1)--(2,2)--(3,2)--(3,3)--(4,3)--(4,4);
\begin{scope}[xshift=7cm]
\fill[black!12] (0,0)--(1,0)--(1,2)--(2,2)--(2,1)--(0,1)--cycle
	(3,2)--(4,2)--(4,3)--(3,3)--cycle;
\draw[gray] (0,0) grid (4,4);
\draw[very thick] (0,0)--(0,1)--(1,1)--(1,2)--(4,2)--(4,4);
\end{scope}
\begin{scope}[xshift=14cm]
\fill[black!12] (1,0)--(3,0)--(3,2)--(2,2)--(2,1)--(1,1)--cycle
	(3,3)--(4,3)--(4,4)--(3,4)--cycle;
\draw[gray] (0,0) grid (4,4);
\draw[very thick] (0,0)--(3,0)--(3,4)--(4,4);
\end{scope}
\end{tikzpicture}
\caption{The lattice paths with type $C_n$ area vectors $(0,0,0,0),(0,-1,1,1)$ and $(1,0,-1,-2)$.}
\label{Figure:areaC}
\end{center}
\end{figure}
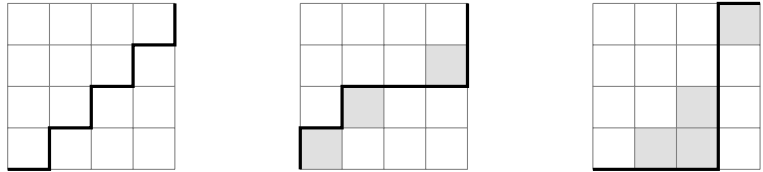

\sk
We conclude this section by proving some auxiliary results on area vectors for later use.

\begin{mylem}{areaC} Let $\pi\in\mathcal L_{n,n}$ be a lattice path with type $C_n$ area vector $\mu$.
\myi{(i)} Let $i,j\in[n]$ with $i<j$ such that $\mu_j=\mu_i-1$ and $\mu_{\ell}\notin\{\mu_i-1,\mu_i\}$ for all $\ell$ with $i<\ell<j$. Then $j=i+1$.
\myi{(ii)} For all $i\in[n-1]$  we have $\mu_i\leq\mu_{i+1}+1$.
\myi{(iii)} Let $j\in[n]$ such that $\mu_j<0$. Then there exist $i\in[j-1]$ with $\mu_i=\mu_j+1$.
\myi{(iv)} Let $i\in[n]$ such that $\mu_i>1$. Then there exists $j\in[n]$ with $i<j$ such that $\mu_j=\mu_i-1$.
\myi{(v)} Let $i\in[n]$ such that $\mu_i=1$ and $\mu_{\ell}\notin\{0,1\}$ for all $\ell\in[n]$ with $i<\ell$. Then $i=n$.
\end{mylem}

\begin{myproof} We start by proving claim~(i). From $\mu=\tau\cdot(\lambda-\nu)$ we obtain $\lambda_{n+1-\ell}=n+1-\ell-\mu_{\ell}\notin\{n+1-\ell-\mu_i,n+2-\ell-\mu_i\}$ for all $\ell$ with $i<\ell<j$. Since $\lambda_{n+1-j}=n+2-j-\mu_i$ and $\lambda_{n-\ell}\leq\lambda_{n+1-\ell}$ it follows inductively that $n+2-\ell-\mu_i<\lambda_{n+1-\ell}$ for all $\ell$ with $i<\ell<j$. But this yields a contradiction for $\ell=i+1$, namely $\lambda_{n-i}>n+1-i-\mu_i=\lambda_{n+1-i}$. Thus $j=i+1$.

Claim~(ii) is an immediate consequence of $\lambda_{n-i}\leq\lambda_{n+1-i}$ for all $i\in[n-1]$. Claim~(iii) follows from~(ii) and $\mu_1\geq0$. Similarly claims~(iv) and~(v) follow from~(ii) and $\mu_n\leq1$.
\end{myproof}

\subsection{The diagonal reading word}\label{Section:readC}

Let $(\pi,u)\in\op{Vert}(C_n)$ be a vertically labelled lattice path. Recall that $(\pi,u)$ corresponds to an element $u\cdot\lambda+(2n+1)\Q$ of the finite torus by \refp{vertC} and hence to a region of the Shi arrangement with minimal alcove $\waf_R\ac$ by \reft{anderson}. That is, $u\cdot\lambda=\A(\waf_R)=-\waf_R\waf_f^{-1}\cdot0$. Our goal for this section is to recover the element $\waf_R$ of the affine Weyl group from the labelled path $(\pi,u)$.

In order to do so write $\waf_R=w\waf_D$ where $w\in\S_n^C$ and $\waf_D\ac$ is the minimal alcove of the dominant Shi region corresponding to $\lambda$ discussed in the previous section (\refs{areaC}). Recall that $w=u\tau\sigma$ by \refp{utausig}. Since we already constructed $\waf_D$ from $\pi$ in the previous section, we are in principle also able to compute the signed permutation $w$. However, it turns out that this permutation can simply be read off the vertically labelled path (diagonally!).

\begin{mydef}{readC} Let $(\pi,u)\in\op{Vert}(C_n)$ be a vertically labelled lattice path and $\mu$ the area vector of $\pi$. Define the type $C_n$ \emph{diagonal reading word} $d_C(\pi,u)$ as follows: For each $i=0,1,\dots,n$ first write down the negative labels $-u(j)$ of the rows with $\mu_{n+1-j}=-i$ from top to bottom, then write down the labels $u(j)$ of rows with $\mu_{n+1-j}=i+1$ from bottom to top.
\end{mydef}

\begin{figure}[ht]
\begin{center}
\begin{tikzpicture}[scale=.55]
\begin{scope}
\draw[gray] (0,0) grid (6,6);
\foreach \x in {0,...,5}{
\pgfmathsetmacro{\y}{\x-.1}
\pgfmathsetmacro{\z}{5.1-\x}
\draw[->,very thick,lightgray](-.9,\z)--(\y,5.9);
\pgfmathsetmacro{\y}{0.1+ \x}
\pgfmathsetmacro{\z}{5.9-\x}
\draw[->,very thick,lightgray](5.9,\z)--(\y,.1);
}
\draw[xshift=5mm,yshift=5mm]
(5,5)node{1}
(4,4)node{2}
(3,3)node{3}
(2,2)node{4}
(1,1)node{5}
(0,0)node{6}
(-1,0)node{7}
(0,1)node{8}
(1,2)node{9}
(2,3)node{10}
(3,4)node{11}
(4,5)node{12}
(5,4)node{13}
(4,3)node{14}
(3,2)node{15}
(2,1)node{16}
(1,0)node{17}
(-1,1)node{18}
(0,2)node{19}
(1,3)node{20}
(2,4)node{21}
(3,5)node{22}
(5,3)node{23}
(4,2)node{24}
(3,1)node{25}
(2,0)node{26}
(-1,2)node{27}
(0,3)node{28}
(1,4)node{29}
(2,5)node{30}
(5,2)node{31}
(4,1)node{32}
(3,0)node{33}
(-1,3)node{34}
(0,4)node{35}
(1,5)node{36}
(5,1)node{37}
(4,0)node{38}
(-1,4)node{39}
(0,5)node{40}
(5,0)node{41}
(-1,5)node{42};
\draw[very thick,red] (0,0)--(0,1)--(1,1)--(1,2)--(2,2)--(2,3)--(3,3)--(3,4)--(4,4)--(4,5)--(5,5)--(5,6)--(6,6)
(6.5,.5)--(7,.5)
(-2,5.5)--(-1.5,5.5)
(-1.75,5.75)--(-1.75,5.25);
\end{scope}
\begin{scope}[xshift=12cm]
	\draw[gray] (0,0) grid (6,6);
	\draw[very thick] (0,0)--(0,1)--(4,1)--(4,6)--(6,6);
	\draw[xshift=5mm,yshift=5mm]
		(-1,0) node{$1$}
		(3,1) node{$-5$}
		(3,2) node{$-4$}
		(3,3) node{$2$}
		(3,4) node{$3$}
		(3,5) node{$6$};
\end{scope}
\end{tikzpicture}
\caption{The diagonal reading order of type $C_6$ indicated left, and a vertically labelled lattice path with $d_C(\pi,u)=[-2,1,3,4,6,5]$ on the right.}
\label{Figure:readC}
\end{center}
\end{figure}
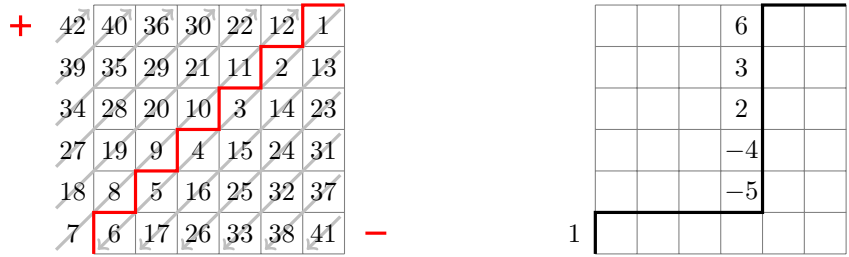

\begin{myex}{readC} Consider the vertically labelled lattice path $(\pi,u)$ from our running example, which is depicted again in \reff{readC} (right). We first deduce the diagonal reading word $d_C(\pi,u)$ following \refd{readC}. In \refe{areaC} we already determined the area vector $\mu=(2,1,0,-1,-2,1)$ of $\pi$.

For $i=0$ we have $-i=0=\mu_3=\mu_{6+1-4}$ hence we write down $-u(4)=-2$. Moreover $i+1=1=\mu_6=\mu_{6+1-1}$ and $i+1=\mu_2=\mu_{6+1-5}$. We write down the labels $u(1)=1$ and $u(5)=3$ from bottom to top, that is, first $1$ and then $3$. Next set $i=1$. We have $-i=-1=\mu_4=\mu_{6+1-3}$ and write down $-u(3)=4$. As $i+1=2=\mu_1=\mu_{6+1-6}$ we write down $u(6)=6$. Finally set $i=2$ to obtain $-i=-2=\mu_5=\mu_{6+1-2}$ and write down $-u(2)=5$. The complete diagonal reading word is thus given by
\begin{align*}
d_C(\pi,u)=[-2,1,3,4,6,5].
\end{align*}

The diagonal reading word of type $C_n$ can also be read off quickly by scanning all boxes that may contain labels according to the diagonal reading order, which is indicated in \reff{readC} (left). In our example the first box with respect to the diagonal reading order to contain a label has the number $3$. The encountered label is $2$. Since this box lies below the alternating path $(NE)^n$ we write down its negative $-2$. The next two boxes that contain labels have numbers $7$ and $11$. Hence we write down $1$ and then $3$. It is not difficult to verify that this procedure always yields $d_C(\pi,u)$.

In light of the following proposition the reader may wish to check that
\begin{align*}
d_C(\pi,u)
=u\tau\sigma
=[1,-5,-4,2,3,6]\cdot[-6,-5,-4,-3,-2,-1]\cdot[3,-6,-2,4,-1,5].
\end{align*}
As a consequence we obtain the affine permutation $\waf_R=d_C(\pi,u)\,\waf_D=[20,10,-2,-9,-21,12]$ that takes the fundamental alcove to the minimal alcove of the Shi region $R$ corresponding to the vertically labelled lattice path $(\pi,u)$. Indeed, we have $\waf_R\waf_f^{-1}=[1,47,48,54,55,58]$ and
\begin{align*}
\A(\waf_R)
=-\waf_R\waf_f^{-1}\cdot0+13\Z^6
=(0,4,4,-4,-4,4)+13\Z^6,
\end{align*}
which is in accordance with \refe{vertC}.
\end{myex}

\sk
The next result confirms that the diagonal reading word is the correct signed permutation.

\begin{myprop}{readC} Let $(\pi,u)\in\op{Vert}(C_n)$ be a vertically labelled lattice path with area vector $\mu$, let $\tau$ be defined as in \refl{wfC}, and define $\sigma\in\S_n^C$ as in \refl{sigma} such that $t_{\mu}\sigma$ is Gra{\ss}mannian. Then $d_C(\pi,u)=u\tau\sigma$.
\end{myprop}

\begin{myproof} Let $i,j\in[n]$. By definition we have $\abs{d_C(\pi,u)(i)}=\abs{u(j)}$ if and only if
\begin{align*}
i
&=\#\big\{r\in[n],\abs{\mu_r}<\abs{\mu_{n+1-j}}\big\}
+\#\big\{r\in[n+1-j]:\mu_r=\mu_{n+1-j}>0\big\} \notag\\
&\quad+\#\big\{r\in[n]:r\geq n+1-j,\mu_r=\mu_{n+1-j}\leq0\big\}
+\#\big\{r\in[n],\mu_r=-\mu_{n+1-j}>0\big\}\\
&=\big\{r\in[n]:\abs{\mu_r\n-r}\leq\abs{\mu_{n+1-j}\n-(n+1-j)}\big\}.
\end{align*}
By \refl{sigma} we obtain $\abs{\sigma(i)}=n+1-j$, hence
\begin{align*}
\abs{d_C(\pi,u)(i)}
=\abs{u(j)}
=\abs{u(n+1-\abs{\sigma(i)})}
=\abs{u\tau\sigma(i)}.
\end{align*}
Moreover $d_C(\pi,u)(i)=u(j)$ if and only if $\mu_{n+1-j}>0$. On the other hand $\tau\sigma(i)=j$ is equivalent to $\sigma(i)=-(n+1-j)<0$, which is the case if and only if $\mu_{n+1-j}>0$ by \refl{sigma}. We can therefore omit the absolute values in the identity above and the proof is complete.
\end{myproof}

\subsection{The zeta map}\label{Section:zetaC}

We are now in a position to define the combinatorial zeta map of type $C_n$, which is made up of the building blocks defined in \refs{paths}.
Note that there are two versions of the map: one defined on square lattice paths that corresponds to the uniform map defined on $W$-orbits, and one defined on vertically labelled lattice paths that corresponds to the uniform map defined on the finite torus.
To avoid confusion we call the zeta map defined on vertically labelled paths the Haglund--Loehr-zeta map.


\begin{mydef}{zetaC} Given a lattice path $\pi\in\mathcal L_{n,n}$ with type $C_n$ area vector $\mu$ we define its image under the \emph{zeta map} of type $C_n$ as
\begin{align*}
\zeta_C(\pi)
=\LS_n^-(\mu)\RS_n^+(\mu)
\LS_{n-1}^-(\mu)\RS_{n-1}^+(\mu)\cdots
\LS_1^-(\mu)\RS_1^+(\mu)
\LS_0^-(\mu)\RS_0^+(\mu).
\end{align*}
Moreover, given a vertically labelled lattice path $(\pi,u)\in\op{Vert}(C_n)$, we define its image under the \emph{Haglund--Loehr-zeta map} of type $C_n$ as
\begin{align*}
\zeta_C(\pi,u)=\big(\zeta_C(\pi),d_C(\pi,u)\big).
\end{align*}
\end{mydef}

Note that by definition $\zeta_C(\pi)$ is a ballot path, that is, $\zeta_C:\mathcal L_{n,n}\to\mathcal B_{2n}$. It will soon turn out that the Haglund--Loehr-zeta map sends vertically labelled lattice paths to diagonally labelled ballot paths.

\begin{figure}[ht]
\begin{center}
\begin{tikzpicture}[scale=.55]
\begin{scope}
\draw[gray](0,0)grid(6,6);
\draw[very thick](0,0)--(0,1)--(4,1)--(4,6)--(6,6);
\draw[xshift=5mm,yshift=5mm]
(-1,0)node{$1$}
(3,1)node{$-5$}
(3,2)node{$-4$}
(3,3)node{$2$}
(3,4)node{$3$}
(3,5)node{$6$};
\end{scope}
\draw[xshift=8cm,yshift=2.5cm,->,thick](0,0)--node[above]{$\zeta_C$}(1,0);
\begin{scope}[xshift=11cm,yshift=-3cm]
\ballot{6}
\draw[very thick](0,0)--(0,2)--(1,2)--(1,3)--(2,3)--(2,5)--(3,5)--(3,6)--(4,6)--(4,7)--(5,7);
\draw[very thick,red]
(0,0)--(0,1)
(0,2)--(1,2)--(1,3)
(2,5)--(3,5)--(3,6);
\draw[xshift=5mm,yshift=5mm]
(0,0)node{$5$}
(1,1)node{$6$}
(2,2)node{$4$}
(3,3)node{$3$}
(4,4)node{$1$}
(5,5)node{$-2$}
(6,6)node{$2$}
(5,7)node{$-1$}
(4,8)node{$-3$}
(3,9)node{$-4$}
(2,10)node{$-6$}
(1,11)node{$-5$}
(0,2)node{\large{$\bullet$}}
(1,3)node{\large{$\bullet$}}
(2,5)node{\large{$\bullet$}}
(3,6)node{\large{$\bullet$}}
(4,7)node{\large{$\circ$}};
\end{scope}
\end{tikzpicture}
\caption{A vertically labelled lattice path and its image under the Haglund--Loehr-zeta map.}
\label{Figure:zetaC}
\end{center}
\end{figure}
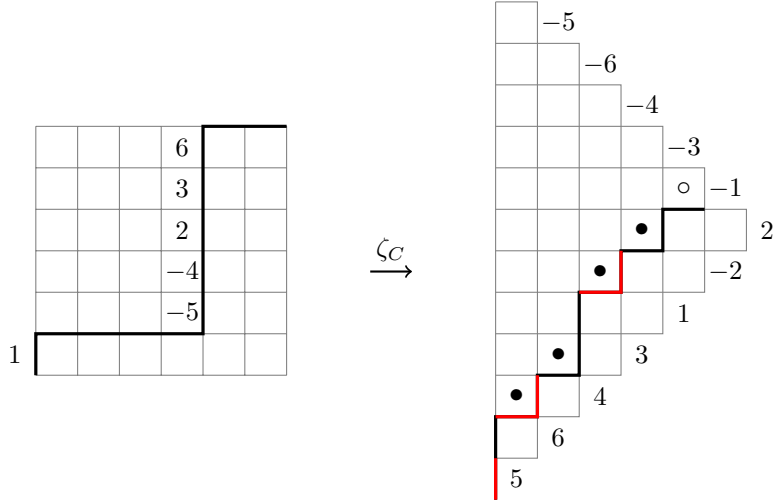

\begin{myex}{zetaC} Recall the area vector $\mu=(2,1,0,-1,-2,1)$ and the diagonal reading word $d=[-2,1,3,4,6,5]$ of the vertically labelled lattice path shown in \reff{zetaC}, which were computed in earlier examples. We construct the individual segments of $\zeta_C(\pi)$ as follows
\begin{align*}
\LS_2^-(\mu)=N, 
\RS_2^+(\mu)=N, 
\LS_1^-(\mu)=EN, 
\RS_1^+(\mu)=ENN, 
\LS_0^-(\mu)=EN,
\RS_0^+(\mu)=ENE.
\end{align*}
\end{myex}

Our first main result of this section is the fact that the zeta map can be inverted using a construction reminiscent of the bounce path of a Dyck path in type $A_{n-1}$.

\begin{mythrm}{zetaC} The zeta map $\zeta_C:\mathcal L_{n,n}\to\mathcal B_{2n}$ is a bijection.
\end{mythrm}

\begin{myproof} Let $\pi\in\mathcal L_{n,n}$ be a lattice path with type $C_n$ area vector $\mu$. For each $k$ with $0\leq k\leq n$ let $\alpha_k$ denote the number of indices $i\in[n]$ such that $\abs{\mu_i}=k$.

Define the \emph{bounce path} of a ballot path $\beta\in\mathcal B_{2n}$ as follows: start at the end point of $\beta$ and go South until you hit the diagonal. Bounce  off it and travel to the West until you reach the upper end of a North step of $\beta$. Bounce off the path $\beta$ to the South until you hit the diagonal again, and repeat until you arrive at $(0,0)$.

Now suppose that $\beta=\zeta_C(\pi)$. By definition of the zeta map the end point of $\beta$ is $(n-\alpha_0,n+\alpha_0)$. The bounce path of $\beta$ meets the diagonal for the first time in the point $(n-\alpha_0,n-\alpha_0)$, and then travels West to the point $(n-\alpha_0-\alpha_1,n-\alpha_0)$, which is the starting point of the segment $\LS_0^-(\mu)\RS_0^+(\mu)$. We claim that this point is a peak of the bounce path.

To see this note that \refl{areaC}~(iii) and~(iv) imply that each non-empty segment $\LS_j^-(\mu)$ or $\RS_j^+(\mu)$ ends with a North step, except possibly $\RS_0^+(\mu)$. In particular, the starting point of any segment $\LS_j^-(\mu)\RS_j(\mu)$ is either $(0,0)$ or the endpoint of a North step of $\beta$.

Inductively the peaks of the bounce path therefore encode the numbers $\alpha_0,\alpha_1,\dots,\alpha_n$.

Knowing $\alpha_0$ and $\alpha_1$ we can recover the number and relative order of zeroes, ones and minus ones in $\mu$ from the segment $\LS_0^-(\mu)\RS_0^+(\mu)$. Since $\LS_0^-(\mu)$ ends with a North step, we first obtain the sequences $\LS_0^-(\mu)$ and $\RS_0^+(\mu)$ and thus the number occurrences of ones and minus ones in $\mu$. Moreover these paths encode the relative order of zeroes and minus ones, respectively the relative order of zeroes and ones. The relative order of ones and minus ones is implied by the following observation: If $\mu_i=1$ and $\mu_j=-1$ for some $i,j\in[n]$ with $i<j$ then there exists $\ell$ with $i<\ell<j$ and $\mu_{\ell}=0$ due to \refl{areaC}~(ii).

Similarly one now reconstructs the numbers of twos and minus twos, as well as the relative order of zeroes, ones, minus ones, twos and minus twos, using the segment $\LS_1^-(\mu)\RS_1^+(\mu)$. Continuing in this fashion one recovers the entire area vector $\mu$. Thus $\zeta_C$ is injective. Since $\mathcal L_{n,n}$ and $\mathcal B_{2n}$ both have cardinality $\binom{2n}{n}$ it is also bijective.
\end{myproof}

\sk
Next we state the respective result for the Haglund--Loehr-zeta map of type $C_n$.

\begin{mythrm}{HLzetaC} The Haglund--Loehr-zeta map $\zeta_C:\op{Vert}(C_n)\to\op{Diag}(C_n)$ is a bijection.
\end{mythrm}
\sk
Note that it follows from \reft{zetaC} that the Haglund--Loehr-zeta map is injective, since the signed permutation $u$ can be recovered from $d_C(\pi,u)$ once $\pi$ is known. The missing fact that $(\pi,u)\in\op{Vert}(C_n)$ if and only if $(\zeta_C(\pi),d_C(\pi,u))\in\op{Diag}(C_n)$ is implied by the next theorem, which relates the rises of a vertically labelled path to the valleys of the corresponding diagonally labelled path.

\sk
Let $i$ be a rise of a vertically labelled lattice path $(\pi,u)$. We say $i$ is labelled $(u(i),u(i+1))$. Let $(i,j)$ be a valley of a diagonally labelled ballot path $(\beta,w)$. We say
\begin{align*}
(i,j)\text{ is labelled by }
\begin{cases}
\big(w(n+1-i),w(n+1-j)\big)&\quad\text{if }j\leq n,\\
\big(w(n+1-i),w(n-j)\big)&\quad\text{if }j>n.
\end{cases}
\end{align*}
Note that with our usual way of picturing diagonally labelled ballot paths, each valley is labelled by the number below it and the number to its right.

\begin{myex}{risevalleyC} The vertically labelled lattice path in \reff{zetaC} has rises $2,3,4$ and $5$, which are labelled $(-5,-4),(-4,2),(2,3)$ and $(3,6)$. Moreover, its initial North step has label $1$. The image under the Haglund--Loehr-zeta map has valleys $(1,3),(2,4),(3,6),(4,7)$ and $(5,8)$, which are labelled by $(5,4),(6,3),(4,-2),(3,2)$ and $(1,-1)$.
\end{myex}

\begin{mythrm}{risevalleyC} Let $(\pi,u)\in\op{Vert}(C_n)$ be a vertically labelled lattice path and $a,b\in u([n])$.
Then $(\pi,u)$ has a rise labelled by $(a,b)$ if and only if $\zeta_C(\pi,u)$ has a valley labelled by $(b,a)$ or $(-a,-b)$.
Moreover $\pi$ begins with a North step labelled by $a$ if and only if $\zeta_C(\pi,u)$ has a valley labelled by $(a,-a)$.
\end{mythrm}

\begin{myproof*}{(Part 1)} We first assume that we are given a valley of $\zeta_C(\pi,u)$ and show that there exists a fitting rise in $(\pi,u)$. Let $\mu$ be the area vector of $\pi$. A valley of $\zeta_C(\pi)$ can either occur within a sequence $\LS_k^-(\mu)$ or $\RS_k^+(\mu)$, or if $\zeta_C(\pi)$ ends with an East step. We treat these three cases independently. No valley may arise at the join of two such sequences because of \refl{areaC}~(iii) and~(iv).

\myi{(1.1)} There is a valley within the sequence $\LS_k^-(\mu)$. Then there exist indices $i,j\in[n]$ with $i<j$ such that $\mu_i=-k$, $\mu_j=-k-1$ and $\mu_{\ell}\notin\{-k-1,-k\}$ for all $\ell$ with $i<\ell<j$. By \refl{areaC}~(i) we have $j=i+1$. Hence $\lambda_{n-i}=n-i-\mu_{i+1}=n+1-i-\mu_i=\lambda_{n+1-i}$ and $n-i$ is a rise of $\pi$. We claim that the labels of our valley are compatible with the labels $(u(n-1),u(n-i+1))$ of this rise.

Suppose $(x,y)$ is our valley. Then $x$ equals the number of East steps in the sequence
\begin{align*}
\LS_n^-(\mu)\RS_n^+(\mu)\cdots\LS_{k+1}^-(\mu)\RS_{k+1}^+(\mu)\LS_k^-(\mu_{i+1},\dots,\mu_n).
\end{align*}
In other words
\begin{align*}
x
&=\#\big\{r\in[n]:\abs{\mu_r}>k+1\big\}+\#\big\{r\in[n]:r\geq i+1,\mu_r=-k-1\big\}\\
&=\#\big\{r\in[n]:\abs{\mu_r\n-r}\geq\abs{(-k-1)\n-(i+1)}\big\}.
\end{align*}
By \refl{sigma} we obtain
\begin{align*}
n+1-x&=\#\big\{r\in[n]:\abs{\mu_r\n-r}\leq\abs{\mu_{i+1}\n-(i+1)}\big\}
=\abs{\sigma^{-1}(i+1)}.
\end{align*}
Since $\mu_{i+1}\leq0$ we have $\sigma^{-1}(i+1)>0$ and thus
\begin{align*}
d_C(\pi,u)(n+1-x)
=u\tau\sigma(n+1-x)
=u\tau(i+1)
=-u(n-i).
\end{align*}
Similarly $y$ equals the number of North steps in the sequence
\begin{align*}
\LS_n^-(\mu)\RS_n^+(\mu)\cdots\LS_{k+1}^-(\mu)\RS_{k+1}^+(\mu)\LS_k^-(\mu_i,\dots,\mu_n).
\end{align*}
We may rewrite this as
\begin{align*}
y
&=\#\big\{r\in[n]:\abs{\mu_r}>k\big\}+\#\big\{r\in[n]:r\geq i,\mu_r=-k\big\}\\
&=\#\big\{r\in[n]:\abs{\mu_r\n-r}\geq\abs{-k\n-i}\big\}.
\end{align*}
As before \refl{sigma} provides
\begin{align*}
n+1-y
&=\#\big\{r\in[n]:\abs{\mu_r\n-r}\leq\abs{\mu_i\n-i}\big\}
=\abs{\sigma^{-1}(i)}
=\sigma^{-1}(i)
\end{align*}
and we compute
\begin{align*}
d_C(\pi,u)(n+1-y)
=u\tau\sigma(n+1-y)
=u\tau(i)
=-u(n+1-i).
\end{align*}
\myi{(1.2)} The valley appears within the sequence $\RS_k^+(\mu)$. Then there exist indices $i,j\in[n]$ with $i<j$ such that $\mu_i=k+1$, $\mu_j=k$ and $\mu_{\ell}\notin\{k,k+1\}$ for all $\ell$ with $i<\ell<j$. We obtain $j=i+1$ and $n-i$ is a rise of $\pi$ just as in~(1.1).

Let $(x,y)$ be our valley. Then $x$ equals the number of East steps in
\begin{align*}
\LS_n^-(\mu)\RS_n^+(\mu)\cdots\LS_k^-(\mu)\RS_k^+(\mu_1,\dots,\mu_i).
\end{align*}
Equivalently
\begin{align*}
x
&=\#\big\{r\in[n]:\abs{\mu_r}>k+1\big\}+\#\big\{r\in[n]:\mu_r=-k-1\big\}+\#\big\{r\in[i]:\mu_r=k+1\big\}\\
&=\#\big\{r\in[n]:\abs{\mu_r\n-r}\geq\abs{(k+1)\n-i}\big\}.
\end{align*}
Since $\mu_i>0$ \refl{sigma} implies
\begin{align*}
n+1-x
&=\#\big\{r\in[n]:\abs{\mu_r\n-r}\leq\abs{\mu_i\n-i}\big\}
=\abs{\sigma^{-1}(i)}
=-\sigma^{-1}(i)
\end{align*}
and we compute
\begin{align*}
d_C(\pi,u)(n+1-x)
=u\tau\sigma(n+1-x)
=-u\tau(i)
=u(n+1-i).
\end{align*}
On the other hand $y$ equals the number of North steps in the sequence
\begin{align*}
\LS_n^-(\mu)\RS_n^+(\mu)\cdots\LS_k^-(\mu)\RS_k^+(\mu_1,\dots,\mu_{i+1}).
\end{align*}
If $k>0$ then
\begin{align*}
y
&=\#\big\{r\in[n]:\abs{\mu_r}>k\big\}+\#\big\{r\in[n]:\mu_r=-k\big\}+\#\big\{r\in[i+1]:\mu_r=k\big\}\\
&=\#\big\{r\in[n]:\abs{\mu_r\n-r}\geq\abs{k\n-(i+1)}\big\},
\end{align*}
and \refl{sigma} provides
\begin{align*}
n+1-y
&=\#\big\{r\in[n]:\abs{\mu_r\n-r}\leq\abs{\mu_{i+1}\n-(i+1)}\big\}
=\abs{\sigma^{-1}(i+1)}
=-\sigma^{-1}(i+1)
\end{align*}
In particular $y\leq n$ and we compute
\begin{align*}
d_C(\pi,u)(n+1-y)
=u\tau\sigma(n+1-y)
=-u\tau(i+1)
=u(n-i).
\end{align*}
Otherwise $k=0$ and 
\begin{align*}
y
=n+\#\big\{r\in[i+1]:\mu_r=0\big\}
=n+\#\big\{r\in[n]:\abs{\mu_r\n-r}\leq i+1\big\}
\end{align*}
From \refl{sigma} we obtain
\begin{align*}
y-n
=\#\big\{r\in[n]:\abs{\mu_r\n-r}\leq\abs{\mu_{i+1}\n-(i+1)}\big\}
=\abs{\sigma^{-1}(i+1)}
=\sigma^{-1}(i+1).
\end{align*}
Since $y>n$ the second label of the valley is given by
\begin{align*}
d_C(\pi,u)(n-y)
=-u\tau\sigma(y-n)
=-u\tau(i+1)
=u(n-i).
\end{align*}
\myi{(1.3)} The path $\zeta_C(\pi)$ ends with an East step. Then there exists $i\in[n]$ such that $\mu_i=1$ and $\mu_{\ell}\notin\{0,1\}$ for all $\ell$ with $i<\ell$. From \refl{areaC}~(v) we know that $i=n$. Consequently $\lambda_1=1-\mu_n=0$ and $\pi$ begins with a North step.

Let $(x,y)$ be the valley above the final East step. Then $x$ equals the number of East steps in $\zeta_C(\pi)$, that is,
\begin{align*}
x
&=\#\big\{r\in[n]:\abs{\mu_r}>0\big\}
=\#\big\{r\in[n]:\abs{\mu_r\n-r}\geq n+1\big\}\\
n+1-x
&=\#\big\{r\in[n]:\abs{\mu_r\n-r}\leq\abs{\mu_n\n-n}\big\}
=\abs{\sigma^{-1}(n)}
=-\sigma^{-1}(n)
\end{align*}
The valley's first label is
\begin{align*}
d_C(\pi,u)(n+1-x)
=u\tau\sigma(n+1-x)
=-u\tau(n)
=u(1).
\end{align*}
On the other hand, $y$ equals the number of North steps of $\zeta_C(\pi)$ plus one. We have
\begin{align*}
y
&=n+1+\#\big\{r\in[n]:\abs{\mu_r\n-r}\leq n\big\}\\
y-n
&=\#\big\{r\in[n]:\abs{\mu_r\n-r}\geq n+1\big\}
=-\sigma^{-1}(n).
\end{align*}
Hence the valley's second label is
\begin{align*}
d_C(\pi,u)(n-y)
=-u\tau\sigma(y-n)
=u\tau(n)
=-u(1).
\end{align*}
\myi{(Part 2)} To complete the proof we need to demonstrate the reverse implication. Thus assume that $i\in[n-1]$ is a rise of $\pi$, that is, $\lambda_i=\lambda_{i+1}$. Then $\mu_{n+1-i}=i-\lambda_i$ and $\mu_{n-i}=i+1-\lambda_{i+1}=\mu_{n+1-i}+1$.

If $\mu_{n-i}=-k\leq0$ then $\mu_{n+1-i}=-k-1$ and we are in the situation of~(1.1). If $\mu_{n-i}=k+1>0$ then $\mu_{n+1-i}=k$ and we are in the situation of~(1.2).

Finally assume that $\pi$ begins with a North step. Then $\lambda_1=0$ and $\mu_n=1-\lambda_1=1$. Hence we are in the situation of~(1.3).
\end{myproof*}

\sk
By now we have assembled all tools needed to prove that the Haglund--Loehr-zeta map $\zeta_C$ is indeed the type $C_n$ special case of the uniform bijection discussed in the \refs{uniform}.

\begin{mythrm}{uniformC} Let $\Phi$ be the root system of type $C_n$ with coroot lattice $\Q$ and zeta map $\zeta$, and let $\psi$ and $\varphi$ be defined as in \refp{vertC} and \refp{diagC}. Then the following diagram commutes.
\begin{center}
\begin{tikzpicture}
\draw (0,0) node(q){$\Q/(2n+1)\Q$};
\draw (4,0) node(park){$\mathrm{Park}(\Phi)$};
\draw (0,1.4) node(v){$\op{Vert}(C_n)$};
\draw (4,1.4) node(d){$\op{Diag}(C_n)$};
\draw[->] (q)--node[below]{$\zeta$}(park);
\draw[->] (v)--node[above]{$\zeta_C$}(d);
\draw[->] (v)--node[left]{$\psi$}(q);
\draw[->] (d)--node[right]{$\varphi$}(park);
\end{tikzpicture}
\end{center}
\end{mythrm}

\begin{myproof} Let $(\pi,u)\in\op{Vert}(C_n)$ and $(\beta,w)=\zeta_C(\pi,u)\in\op{Diag}(C_n)$. Moreover let $\lambda$ be the partition with South East boundary $\pi$, define $\nu$ and $\tau$ as in \refl{wfC}, and let $\mu$ be the area vector of $\pi$. Chose $\sigma\in\S_n^C$ such that $t_{\mu}\sigma$ is a Gra{\ss}mannian affine permutation.

By \refp{readC} we have $w=u\tau\sigma$, which agrees with \refp{utausig}. Thus, by \refp{antichain} it suffices to show that
\begin{align}\label{eq:antichainC}
A_{\beta}=(\tau\sigma)^{-1}\cdot\walls
\end{align}
We start out by proving that $(\tau\sigma)^{-1}\cdot\walls\subseteq A_{\beta}$. To this end let $i\in[n-1]$. Then
\begin{align*}
\alpha_i\in\walls
&\Leftrightarrow\lambda_i=\lambda_{i+1}
\Leftrightarrow\pi\text{ has a rise }i.
\end{align*}
By \reft{risevalleyC} this is the case if and only if $(\beta,w)$ has a valley $(x,y)$ labelled by $(u(i+1),u(i))$ or $(-u(i),-u(i+1))$. Moreover, a closer look at the proof of \reft{risevalleyC} reveals that the second case only occurs if $y\leq n$. In other words
\begin{align*}
u(i+1)&=w(n+1-x)
&\text{or}
&&-u(i)&=w(n+1-x),\\
u(i)&=\begin{cases}w(n+1-y)&\quad\text{if }y\leq n,\\w(n-y)&\quad\text{if }y>n,\end{cases}
&&&-u(i+1)&=w(n+1-y)\quad\text{and }y\leq n.
\end{align*}
Applying $u^{-1}$ to the above identities yields
\begin{align*}
i+1&=\tau\sigma(n+1-x)
&\text{or}
&&-i&=\tau\sigma(n+1-x)\\
i&=
\begin{cases}
\tau\sigma(n+1-y)&\quad\text{if }y\leq n,\\
\tau\sigma(n-y)&\quad\text{if }y>n,
\end{cases}
&&&-i-1&=
\tau\sigma(n+1-y)\quad\text{and }y\leq n.
\end{align*}
We obtain
\begin{align*}
(\tau\sigma)^{-1}\cdot\alpha_i=\alpha_{x,y}\in A_{\beta}.
\end{align*}
Furthermore $\alpha_0\in\walls$ if and only if $\lambda_1=0$, that is, $\pi$ begins with a North step. By \reft{risevalleyC} this is equivalent to $(\beta,w)$ ending with an East step such that the corresponding valley $(x,y)$ has labels $(w(n+1-x),w(n-y))=(u(1),-u(1))$. Thus $\tau\sigma(n+1-x)=1=\tau\sigma(y-n)$ and we obtain
\begin{align*}
(\tau\sigma)^{-1}\cdot\alpha_0=\alpha_{x,y}\in A_{\beta}.
\end{align*}
Finally we have already observed that $\tilde\alpha\notin\walls$, and the first inclusion in \refq{antichainC} follows.

Conversely, let $\alpha_{x,y}\in A_{\beta}$. By similar reasoning as above it follows that $\tau\sigma\cdot\alpha_{x,y}$ is of the form $e_{i+1}-e_i$ for a rise $i$ of $\pi$, unless the valley comes from a terminal East step of $\beta$, in which case $\pi$ begins with a North step and $\tau\sigma\cdot\alpha_{x,y}=2e_1=\alpha_0$. Hence, we also have $A_{\beta}\subseteq(\tau\sigma)^{-1}\cdot\walls$ and the proof is complete.
\end{myproof}

\subsection{The sweep map}\label{Section:sweepC}

In type $A_{n-1}$ there is a generalisation of the zeta map to rational Dyck paths called the sweep map \cite{ALW2014SweepMaps}. The concept of the sweep map is as follows. Given a path one assigns to each step a label, the labels being distinct integers. To obtain the image of a path under the sweep map, one rearranges the steps such that the labels are in increasing order.

We now give a description of the zeta map of type $C_n$ similar to the sweep map on Dyck paths. Given a path $\pi=s_1s_2,\dots,s_{2n}\in\mathcal L_{n,n}$, where $s_i\in\{E,N\}$, assign a label $\ell_i$ to each step $s_i$ by setting $\ell_1=0$, and $\ell_{i+1}=\ell_i+2n+1$ if $s_i=N$, and $\ell_{i+1}=\ell_i-2n$ if $s_i=E$. Now define a collection $X$ of labelled steps as follows. If $\ell_i<0$ then add $(s_i,\ell_i)$ to $X$. If $\ell_i>0$ then add $(s_{i-1},-\ell_i)$. Finally, for the step $s_1$ which is the only step labelled $0$, add $(s_{2n},-n)$. Thus $X$ contains $2n$ labelled steps.

Now draw a path as follows. Choose $(s,\ell)\in X$ such that $\ell$ is the minimal label among all pairs in $X$. Draw the step $s$ and remove $(s,\ell)$ from $X$. Repeat until $X$ is empty. We denote the path obtained in this way by $\op{sw}(\pi)$. See \reff{sweepC}.

\begin{figure}[ht]
\begin{center}
\begin{tikzpicture}[scale=.6]
\begin{scope}
\draw[gray] (0,0) grid (6,6);
\draw[very thick] (0,0)--(0,1)--(4,1)--(4,6)--(6,6);
\draw[xshift=5mm,yshift=-1mm,anchor=south]
	(0,1) node{\small{$13$}}
	(1,1) node{\small{$1$}}
	(2,1) node{\small{$-11$}}
	(3,1) node{\small{$-23$}}
	(4,6) node{\small{$30$}}
	(5,6) node{\small{$18$}};
\draw[xshift=-1mm,yshift=5mm,anchor=west]
	(0,0) node{\small{$0$}}
	(4,1) node{\small{$-35$}}
	(4,2) node{\small{$-22$}}
	(4,3) node{\small{$-9$}}
	(4,4) node{\small{$4$}}
	(4,5) node{\small{$17$}};
\end{scope}
\begin{scope}[xshift=8cm,yshift=5cm]
\begin{scope}
\draw[very thick](0,0)--(1,0);
\draw[xshift=5mm,yshift=-1mm,anchor=south] (0,0) node{\small{$-11$}};
\end{scope}
\begin{scope}[xshift=15mm]
\draw[very thick](0,0)--(1,0);
\draw[xshift=5mm,yshift=-1mm,anchor=south] (0,0) node{\small{$-23$}};
\end{scope}
\begin{scope}[xshift=30mm]
\draw[very thick](0,0)--(0,1);
\draw[xshift=5mm,yshift=-1mm,anchor=south] (0,0) node{\small{$-35$}};
\end{scope}
\begin{scope}[xshift=45mm]
\draw[very thick](0,0)--(0,1);
\draw[xshift=5mm,yshift=-1mm,anchor=south] (0,0) node{\small{$-22$}};
\end{scope}
\begin{scope}[xshift=60mm]
\draw[very thick](0,0)--(0,1);
\draw[xshift=5mm,yshift=-1mm,anchor=south] (0,0) node{\small{$-9$}};
\end{scope}
\end{scope}
\begin{scope}[xshift=8cm,yshift=1cm]
\begin{scope}[yshift=12mm]
\draw[very thick](0,0)--(0,1);
\draw[xshift=5mm,yshift=-1mm,anchor=south] (0,0) node{\small{$-13$}};
\end{scope}
\begin{scope}[xshift=15mm,yshift=12mm]
\draw[very thick](0,0)--(1,0);
\draw[xshift=5mm,yshift=-1mm,anchor=south] (0,0) node{\small{$-1$}};
\end{scope}
\begin{scope}[xshift=30mm,yshift=12mm]
\draw[very thick](0,0)--(0,1);
\draw[xshift=5mm,yshift=-1mm,anchor=south] (0,0) node{\small{$-4$}};
\end{scope}
\begin{scope}
\draw[very thick](0,0)--(0,1);
\draw[xshift=5mm,yshift=-1mm,anchor=south] (0,0) node{\small{$-17$}};
\end{scope}
\begin{scope}[xshift=15mm]
\draw[very thick](0,0)--(0,1);
\draw[xshift=5mm,yshift=-1mm,anchor=south] (0,0) node{\small{$-30$}};
\end{scope}
\begin{scope}[xshift=30mm]
\draw[very thick](0,0)--(1,0);
\draw[xshift=5mm,yshift=-1mm,anchor=south] (0,0) node{\small{$-18$}};
\end{scope}
\end{scope}
\begin{scope}[xshift=13.5cm,yshift=17mm]
\draw[very thick](0,0)--(1,0);
\draw[xshift=5mm,yshift=-1mm,anchor=south] (0,0) node{\small{$-6$}};
\end{scope}
\begin{scope}[xshift=17cm,yshift=-3cm]
\ballot{6}
\draw[very thick] (0,0)--(0,2)--(1,2)--(1,3)--(2,3)--(2,5)--(3,5)--(3,6)--(4,6)--(4,7)--(5,7);
\draw[xshift=-1mm,yshift=5mm,anchor=west]
	(0,0) node{\small{$-35$}}
	(0,1) node{\small{$-30$}}
	(1,2) node{\small{$-22$}}
	(2,3) node{\small{$-17$}}
	(2,4) node{\small{$-13$}}
	(3,5) node{\small{$-9$}}
	(4,6) node{\small{$-4$}};
\draw[xshift=5mm,yshift=-1mm,anchor=south]
	(0,2) node{\small{$-23$}}
	(1,3) node{\small{$-18$}}
	(2,5) node{\small{$-11$}}
	(3,6) node{\small{$-6$}}
	(4,7) node{\small{$-1$}};
\end{scope}
\end{tikzpicture}
\caption{The labelling of $\pi$ (left), the set $X$ of labelled steps (middle), and the path $\op{sw}(\pi)$ of steps in increasing order (right).}
\label{Figure:sweepC}
\end{center}
\end{figure}
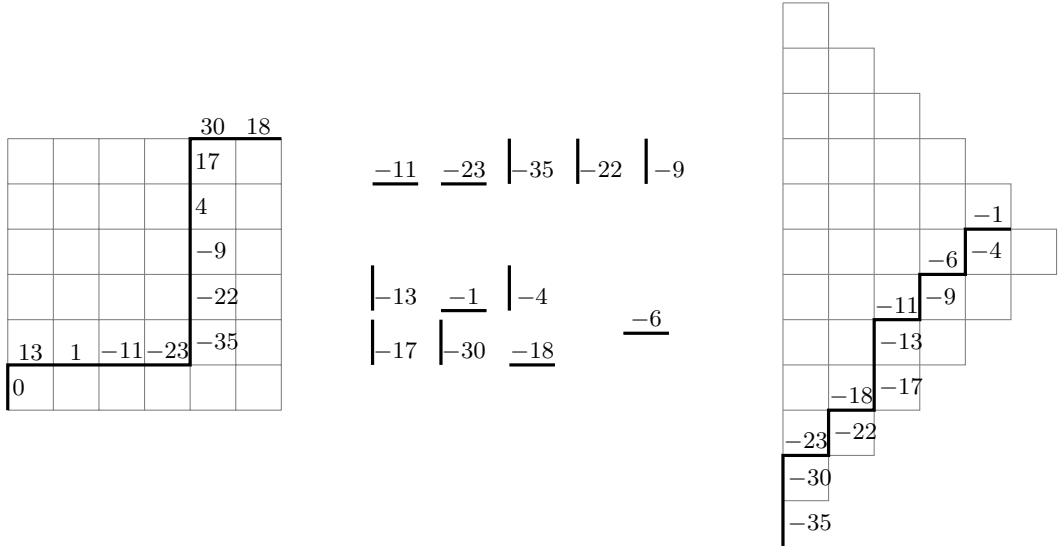

\begin{mythrm}{sweepC} For each lattice path $\pi\in\mathcal L_{n,n}$ we have $\op{sw}(\pi)=\zeta_C(\pi)$. In particular, the sweep map $\op{sw}:\mathcal L_{n,n}\to\mathcal B_{2n}$ is a bijection.
\end{mythrm}

\begin{myproof} The proof consists of a straightforward but rather tedious case by case analysis of the involved labels. Let $\mu$ be the type $C_n$ area vector of a path $\pi\in\mathcal L_{n,n}$. We use the following notation. The label of the $i$-th North step of $\pi$ is denoted by $\ell_i^N$. The corresponding labelled step which is added to $X$ is denoted by $(s_i^N,x_i^N)$.

We pair each North step with an East step. If the North step has a non-negative label, this is the next East step in the same diagonal. If the North step has a negative label, this is the previous East step in the same diagonal. We denote by $\ell_i^E$ the label of the East step corresponding to the $i$-th North step, and by $(s_i^E,x_i^E)$ the associated labelled step in $X$.

For example in \reff{sweepC} we have $\ell_4^N=-9$, $(s_4^N,x_4^N)=(N,-9)$, $\ell_4^E=1$, and $(s_4^E,x_4^E)=(E,-1)$. Also $\ell_6^N=17$, $(s_6^N,x_6^N)=(N,-17)$, $\ell_6^E=30$, and $(s_6^E,x_6^E)=(N,-30)$.

\sk
The label of the $i$-th North step is
\begin{align}
\label{eq:liN} \ell_i^N=(i-1)(2n+1)-(i-\mu_{n-i+1})(2n)=2n(\mu_{n-i+1}-1)+i-1.
\end{align}
First consider the case $\mu_{n-i+1}>0$. Then
\begin{align*}
2n|\mu_{n-i+1}|-2n\leq\ell_i^N<2n|\mu_{n-i+1}|-n.
\end{align*}
If $i>1$ then $x_i^N=-\ell_i^N$. If $i=1$ then $x_i^N=-n$. Hence
\begin{align*}
-2n|\mu_{n-i+1}|+n\leq x_i^N<-2n|\mu_{n-i+1}|+2n.
\end{align*}
On the other hand, if $\mu_{n-i+1}\leq0$ then
\begin{align*}
-2n|\mu_{n-i+1}|-2n\leq\ell_i^N=x_i^N<-2n|\mu_{n-i+1}|-n.
\end{align*}
Now, let us treat the East steps. We start with the case $\mu_{n-i+1}>0$. Then
\begin{align*}
\ell_i^E=\ell_i^N+2n+k_i=2n|\mu_{n-i+1}|+i-1+k_i,
\end{align*}
for some $k_i\in[n-i+1]$. Since $x_i^E=-\ell_i^E$, we obtain
\begin{align}
\label{eq:xiEpos}
-2n|\mu_{n-i+1}|-n\leq x_i^E<-2n|\mu_{n-i+1}|.
\end{align}
If $\mu_{n-i+1}\leq0$ then
\begin{align*}
\ell_i^E=\ell_i^N+2n-k_i=-2n|\mu_{n-i+1}|+i-1-k_i,
\end{align*}
for some $k_i\in\{0,\dots,i-1\}$. If $\mu_{n+i-1}<0$ then
\begin{align*}
-2n|\mu_{n-i+1}|\leq\ell_i^E=x_i^E<-2n|\mu_{n-i+1}|+n.
\end{align*}
Finally assume $\mu_{n-i+1}=0$.
If $\ell_i^E=0$ then $x_i^E=-n$. Otherwise $x_i^E=-\ell_i^E$. Combined this yields
\begin{align*}
-n\leq x_i^E<0,
\end{align*}
which is a special case of \refq{xiEpos}.

\sk
We make the following observation. If $-2nk\leq x_i^s<-2nk+n$, where $k=1,2,\dots$, then either $s=E$ and $\mu_{n-i+1}=-k$, or $s=N$ and $\mu_{n-i+1}=-k+1$. Similarly, if $-2nk+n\leq x_i^s<-2nk+2n$, where $k=1,2,\dots$, then either $s=N$ and $\mu_{n-i+1}=k$ or $s=E$ and $\mu_{n-i+1}=k-1$.

Thus by definition the path $\op{sw}(\pi)$ is composed of segments $T(-k,-k+1)$, $T(k-1,k)$, where $k=1,2,\dots$, such that each step of $T(-k,-k+1)$ corresponds to an entry of the area vector $\mu_{n-i+1}\in\{-k,-k+1\}$, and each step of $T(k-1,k)$ corresponds an entry $\mu_{n-i+1}\in\{k-1,k\}$.

This is a good sign because the path $\zeta_C(\pi)$ is also composed of segments with the same property. Indeed we will carry out the proof by showing that $T(-k,-k+1)=\LS_{k-1}^-(\mu)$ and $T(k-1,k)=\RS_{k-1}^+(\mu)$.

\sk
We first prove that $T(-k,-k+1)=\LS_{k-1}^-(\mu)$ for $k\geq1$.
As $\mu_{n-i+1}=-k<0$ implies $x_i^E=\ell_i^E<0$ and $\mu_{n-i+1}=-k+1\leq0$ implies $x_i^N=\ell_i^N<0$, we have $s_i^E=E$ and $s_j^N=N$.
That is, every entry $\mu_{n-i+1}=-k$ will contribute an East step while each entry $\mu_{n-i+1}=-k+1$ contributes a North step.
This is consistent with the definition of the zeta map.
Therefore, it suffices to check that $i<j$ implies
\begin{align}
\label{eq:xiExjNneg}
x_i^E&<x_j^N \qquad\text{if }\mu_{n-j+1}=-k, \mu_{n-i+1}=-k+1,\\
\label{eq:xiNxjEneg}
x_i^N&<x_j^E \qquad\text{if }\mu_{n-j+1}=-k+1, \mu_{n-i+1}=-k.
\end{align}
Inequality \refq{xiExjNneg} is trivial as
\begin{align*}
x_i^E=-2nk+i-1-k_i<2n(-k+1-1)+j-1=x_j^N
\end{align*}
To see \refq{xiNxjEneg} note that $\mu_{n-j+1}=-k$ and $\mu_{n-i+1}=-k+1$ imply that the path $\pi$ has an East step in the same diagonal as its $i$-th North step somewhere between its $i$-th and $j$-th North steps.
That is, $k_j\leq j-i-1$ and therefore
\begin{align*}
x_i^N=2n(-k+1-1)+i-1<-2nk+j-1-k_j=x_j^E.
\end{align*}
Next we show that $T(k-1,k)=\RS_{k-1}^+(\mu)$ for $k\geq2$.
In this case $\mu_{n-i+1}=k>1$ implies $\ell_i^N>0$ and $\mu_{n-i+1}=k-1>0$ implies $\ell_i^E>0$.
This case is more difficult (confusing) because we do not necessarily have $s_i^E=E$ and $s_j^N=N$.
Instead, if $\mu_{n-i+1}=k-1$ and $\mu_{n-j+1}=k$ then
\begin{align*}
s_i^E=
\begin{cases}
E\quad\text{if }\mu_{n-i}=k,\\
N\quad\text{if }\mu_{n-i}\leq k-1,
\end{cases}
&&
s_j^N=
\begin{cases}
N\quad\text{if }\mu_{n-j+2}=k-1,\\
E\quad\text{if }\mu_{n-j+2}\geq k.
\end{cases}
\end{align*}
Thus an entry $\mu_{n-i+1}=k-1$ contributes an East step instead of a North step if and only if the previous entry $\mu_{n-i}=k$ contributes a North step instead of an East step.
We see that the number of East and North steps in $T(k-1,k)$ is consistent with the definition of the zeta map. 

To see that also the relative orders of the steps in $T(k-1,k)$ and $\RS_{k-1}^+(\mu)$ agree, it suffices to prove that $i<j$ implies
\begin{align}
\label{eq:xiExjEpos} x_j^E&<x_i^E\qquad\text{if }\mu_{n-j+1}=\mu_{n-i+1}=k-1,\\
\label{eq:xiNxjNpos} x_j^N&<x_i^N\qquad\text{if }\mu_{n-j+1}=\mu_{n-i+1}=k,\\
\label{eq:xiNxjEpos} x_j^E&<x_i^N\qquad\text{if }\mu_{n-j+1}=k-1,\mu_{n-i+1}=k.
\end{align}
Moreover, let $i<j$, $\mu_{n-j+1}=k$ and $\mu_{n-i+1}=k-1$. Then we require that
\begin{align}
\label{eq:xiExjNpos} x_i^E&<x_j^N
\end{align}
if and only if there exists no $r$ such that $i<r<j$ and $\mu_{n-r+1}=k-1$.

From \refq{xiExjEpos}--\refq{xiExjNpos} it follows that the order of the steps is (almost) obtained by reading the area vector from left to right drawing $s_i^E$ whenever $\mu_{n-i+1}=k-1$, and $s_i^N$ whenever $\mu_{n-i+1}=k$. The only exception to this rule is when $\mu_{n-i+1}=k-1$ and $\mu_{i+1}=k$. In this case one has to draw the step $s_i^E=E$ before the step $s_{i+1}^N=N$. \footnote{In fact, one also draws $s_i^E=E$ before all the East steps coming from entries of the area vector equal to $k$ occurring between the $(n-i+1)$-th and the $(n-t+1)$-th entry, where $t$ is minimal such that $i<t$ and $\mu_{n-t+1}=k-1$. However, permuting East steps clearly has no effect on the resulting path.}

We now prove the claims \refq{xiExjEpos}--\refq{xiExjNpos}.
If $i<j$ and $\mu_{n-i+1}=\mu_j=k-1$ then there must be an East step on the diagonal between the $i$-th and $j$-th North steps of $\pi$.
Hence $k_i\leq j-i$ and we obtain \refq{xiExjEpos}.
\begin{align*}
x_i^E=-2n(k-1)-i+1-k_i\geq-2n(k-1)-j+1>-2n(k-1)-j+1-k_j=x_j^E
\end{align*}
The inequalities \refq{xiNxjNpos} and \refq{xiNxjEpos} are trivial as
\begin{align*}
x_i^N=-2n(k-1)-i+1&>
\begin{cases}
-2n(k-1)-j+1=x_j^N\\
-2n(k-1)-j+1-k_j=x_j^E.
\end{cases}
\end{align*}
To see claim \refq{xiExjNpos}, first assume that there is no $r$ with $i<r<j$ and $\mu_r=k-1$. Then $k_i\geq j-i+1$ and
\begin{align*}
x_i^E=-2n(k-1)-i+1-k_i\leq-2n(k-1)-j<-2n(k-1)-j+1=x_j^N
\end{align*}
On the other hand, if there is such an $r$ then $k_i\leq r-i<j-i$ and we obtain
\begin{align*}
x_i^E=-2n(k-1)-i+1-k_i>-2n(k-1)-j+1=x_j^N.
\end{align*}

Finally, we need to show $T(0,1)=\RS_{0}^+(\mu)$. Let $\mu_{n-j+1}=0$ and assume $\ell_j^E\neq0$. Then $x_j^E=-j+1+k_j$. Choose $i$ maximal such that $i<j$ and $\mu_{n-i+1}\in\{0,1\}$. Note that such an $i$ always exists. Then
\begin{align*}
s_j^E=
\begin{cases}
N\quad\text{if }\mu_{n-i+1}=0,\\
E\quad\text{if }\mu_{n-i+1}=1.
\end{cases}
\end{align*}
If $\ell_j^E=0$ then $x_j^E=-n$. In this case $\mu_{n-i+1}<0$ for all $i<j$. Instead choose $i\leq n$ maximal such that $\mu_{n-i+1}\in\{0,1\}$. Then
\begin{align*}
s_j^E=
\begin{cases}
N\quad\text{if }\mu_{n-i+1}=0,\\
E\quad\text{if }\mu_{n-i+1}=1.
\end{cases}
\end{align*}
Now let $\mu_{n-j+1}=1$ and assume $j>1$. Then $x_j^N=-j+1$. Choose $i$ maximal such that $i<j$ and $\mu_{n-i+1}\in\{0,1\}$. Again such an $i$ always exists. Then
\begin{align*}
s_j^N=
\begin{cases}
N\quad\text{if }\mu_{n-i+1}=0,\\
E\quad\text{if }\mu_{n-i+1}=1.
\end{cases}
\end{align*}
If $j=1$ then $\ell_j^N=0$ and $x_j^N=-n$. In this case choose $i\leq n$ maximal such that $\mu_{n-i+1}\in\{0,1\}$. Then
\begin{align*}
s_j^N=
\begin{cases}
N\quad\text{if }\mu_{n-i+1}=0,\\
E\quad\text{if }\mu_{n-i+1}=1.
\end{cases}
\end{align*}
We see that every entry of the area vector equal to zero contributes a North step and every entry equal to one contributes an East step. Again this is consistent with the zeta map. To see that the relative orders of North and East steps in $T(0,1)$ and $\RS_0^+(\mu)$ are the same, it suffices to show that for all $i<j$
\begin{align*}
-n&<x_i^E\qquad\text{if }\ell_i^E\neq0,\\
-n&<x_i^N\qquad\text{if }i>1,\\
x_j^E&<x_i^E\qquad\text{if }\mu_{n-j+1}=\mu_{n-i+1}=0,\ell_i^E\neq0,\\
x_j^N&<x_i^N\qquad\text{if }\mu_{n-j+1}=\mu_{n-i+1}=1,i>1,\\
x_j^E&<x_i^N\qquad\text{if }\mu_{n-j+1}=0,\mu_{n-i+1}=1,i>1,\\
x_j^N&<x_i^E\qquad\text{if }\mu_{n-j+1}=0,\mu_{n-i+1}=1,\ell_i^E\neq0.
\end{align*}
\end{myproof}

\subsection{The statistics area and dinv}\label{Section:statC}

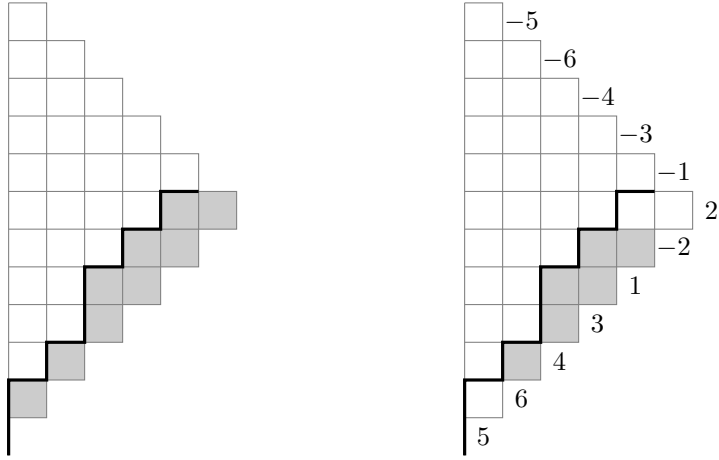
\begin{figure}[t]
\begin{center}
\begin{tikzpicture}[scale=.5]
\begin{scope}
	\fill[black!20] (0,1)--(0,2)--(1,2)--(1,3)--(2,3)--(2,5)--(3,5)--(3,6)--(4,6)--(4,7)--(6,7)--(6,6)--(5,6)--(5,5)--(4,5)--(4,4)--(3,4)--(3,3)--(2,3)--(2,2)--(1,2)--(1,1)--cycle;
	\ballot{6}
	\draw[very thick] (0,0)--(0,2)--(1,2)--(1,3)--(2,3)--(2,5)--(3,5)--(3,6)--(4,6)--(4,7)--(5,7);
\end{scope}
\begin{scope}[xshift=12cm]
	\fill[black!20] (1,2)--(1,3)--(2,3)--(2,5)--(3,5)--(3,6)--(5,6)--(5,5)--(4,5)--(4,4)--(3,4)--(3,3)--(2,3)--(2,2)--cycle;
	\ballot{6}
	\draw[very thick] (0,0)--(0,2)--(1,2)--(1,3)--(2,3)--(2,5)--(3,5)--(3,6)--(4,6)--(4,7)--(5,7);
	\draw[xshift=5mm,yshift=5mm]
		(0,0) node{$5$}
		(1,1) node{$6$}
		(2,2) node{$4$}
		(3,3) node{$3$}
		(4,4) node{$1$}
		(5,5) node{$-2$}
		(6,6) node{$2$}
		(5,7) node{$-1$}
		(4,8) node{$-3$}
		(3,9) node{$-4$}
		(2,10) node{$-6$}
		(1,11) node{$-5$}
;
\end{scope}
\end{tikzpicture}
\caption{A ballot path $\beta\in\mathcal{B}_{12}$ with $\area(\beta)=9$ (left), and a diagonally labelled ballot path $(\beta,w)\in\op{Diag}(C_6)$ with $\area'(\beta,w)=6$ (right). The contributing boxes below the paths are shaded grey.}
\label{Figure:areastatC}
\end{center}
\end{figure}

In this subsection we introduce a natural statistic $\area:\mathcal B_{2n}\to\N$ as well as an extension to labelled paths $\area':\op{Diag}(C_n)\to\N$. We then introduce a less obvious statistic $\dinv_C:\mathcal{L}_{n,n}\to\N$ and its extension $\dinv_C':\op{Vert}(C_n)$. All of these statistics have a known analogue in type $A_{n-1}$. The main result of this subsection is that the zeta map $\zeta_C$ transforms $\dinv$ into $\area$ and $\dinv'$ into $\area'$.

\begin{mydef}{area} Let $\beta\in\mathcal{B}_{2n}$ be a ballot path corresponding to the antichain $A_{\beta}\subseteq\Phi^+$ in the root poset as in \refl{AbetaC}. We define the $\area(\beta)$ as the cardinality of the order ideal $I_{\beta}\subseteq\Phi^+$ given by
\begin{align*}
I_{\beta}
&=\big\{x\in\Phi^+:\text{there exists no }y\in A_{\beta}\text{ with }y\leq x\big\}.
\end{align*}
\end{mydef}

\sk
First note that $\area$ is really a statistic on antichains in the root poset, and as such \refd{area} is completely uniform.
Moreover, its type $A_{n-1}$ instance coincides with the $\area$ statistic on Dyck paths.
Secondly, note that in type $C_n$ the number $\area(\beta)$ equals the number of boxes ``below'' the path $\beta$ in the sense of \reff{areastatC} (left).

\begin{mydef}{area'} Let $(\beta,w)\in\op{Diag}(C_n)$ be a diagonally labelled Dyck path and let $I_{\beta}$ be as in \refd{area} above.
Then we define
\begin{align*}
\area'(\beta,w)
&=\#\big\{x\in I_{\beta}:w\cdot x\in\Phi^+\big\}.
\end{align*}
\end{mydef}

Note that $\area'$ is really a statistic on non-nesting parking functions, and as such \refd{area'} is uniform as well.
The type $A_{n-1}$ instance of $\area'$ was first considered  by Haglund and Loehr~\cite{HagLoe2005}.
Similar uniform descriptions of the statistic $\area'$ using the Shi arrangement were given by Armstrong~\cite{Armstrong2013} and Stump~\cite{Stump2010}.
In type $C_n$ the number $\area'(\beta,w)$ can be visualised easily as the number of boxes below the path $\beta$ such that the label to its right is smaller that the label below it. See \reff{areastatC} (right).

\sk
Next we define the somewhat less obvious statistics $\dinv_C$ and $\dinv_C'$.

\begin{mydef}{dinvC} Let $\pi\in\mathcal{L}_{n,n}$ be a lattice path with type $C_n$ area vector $\mu$. Then we define
\begin{align*}
\dinv_C(\pi)
&=\#\{(i,j):i<j,\mu_{n-i+1}=\mu_{n-j+1}\}
+\#\{(i,j):i<j,\mu_{n-i+1}=\mu_{n-j+1}+1\}\\
&\quad+\#\{(i,j):i<j,\mu_{n-i+1}=-\mu_{n-j+1}\}
+\#\{(i,j):i<j,\mu_{n-i+1}=-\mu_{n-j+1}+1\}\\
&\quad+\#\{i:\mu_{n-i+1}=0\}
\end{align*}
A pair $(i,j)\in[n]\times[n]$ or index $i\in[n]$ contributing to $\dinv_C(\pi)$ is called \emph{diagonal inversion} of $\pi$.
\end{mydef}

\begin{mydef}{dinvC'} Let $(\pi,u)\in\op{Vert}(C_n)$ be a vertically labelled lattice path and let $\mu$ be the type $C_n$ area vector of $\pi$.
Then we define
\begin{align*}
\dinv_C'(\pi)
&=\#\{(i,j):i<j,\mu_{n-i+1}=\mu_{n-j+1},u(i)<u(j)\}\\
&\quad+\#\{(i,j):i<j,\mu_{n-i+1}=\mu_{n-j+1}+1,u(i)>u(j)\}\\
&\quad+\#\{(i,j):i<j,\mu_{n-i+1}=-\mu_{n-j+1},u(i)<-u(j)\}\\
&\quad+\#\{(i,j):i<j,\mu_{n-i+1}=-\mu_{n-j+1}+1,u(i)>-u(j)\}\\
&\quad+\#\{i:\mu_{n-i+1}=0,u(i)<0\}
\end{align*}
We call a pair $(i,j)\in[n]\times[n]$ or index $i\in[n]$ contributing to $\dinv_C'(\pi,u)$ a \emph{diagonal inversion} of $(\pi,u)$.
\end{mydef}

\begin{myex}{dinvC} Consider the vertically labelled lattice path $(\pi,u)$ in \reff{vertC}. 

The path $\pi$ has area vector $\mu=(2,1,0,-1,-2,1)$.
The pair $(1,5)$ is a diagonal inversion of type $\mu_{n-i+1}=\mu_{n-j+1}$, $(1,4)$ is a diagonal inversion of type $\mu_{n-i+1}=\mu_{n-j+1}+1$, $(1,3),(2,6)$ and $(3,5)$ are diagonal inversions of type $\mu_{n-i+1}=-\mu_{n-j+1}$, $(1,4),(3,6)$ and $(4,5)$ are diagonal inversions of type $\mu_{n-i+1}=-\mu_{n-j+1}+1$, and $4$ is a diagonal inversion of type $\mu_{n-i+1}=0$.
Together that makes $\dinv(\pi)=9$.
Note that the pair $(1,4)$ is counted twice!

\sk
We have $u(1)=1<u(5)=3$ thus $(1,5)$ is a diagonal inversion of $(\pi,u)$ of type $\mu_{n-i+1}=\mu_{n-j+1}+1$.
$u(1)=1<u(4)=2$ thus $(1,4)$ is no diagonal inversion of type $\mu_{n-i+1}=\mu_{n-i+1}+1$. However, $u(1)>-u(4)$ thus $(1,4)$ is still a diagonal inversion of type $\mu_{n-i+1}=-\mu_{n-i+1}+1$.
The only other two diagonal inversions of $\pi$ that do not contribute to $\dinv_C'(\pi,u)$ are $(2,6)$ and $4$.
We obtain $\dinv_C'(\pi,u)=6$.
\end{myex}

\sk
The $\dinv$ statistic of type $A_{n-1}$ was suggested by Haiman while the $\dinv'$ statistic of type $A_{n-1}$ is due to Haglund and Loehr.
Both statistics are found in~\cite{HagLoe2005}.
The type $C_n$ version is a natural analogue in the following sense.
If $\pi\in\mathcal{L}_{n,n}$ is a Dyck path then it is actually assigned two numbers $\dinv_A(\pi)$ and $\dinv_C(\pi)$.
In this case \refd{dinvC} agrees with the definition of Haiman, that is, $\dinv_A(\pi)=\dinv_C(\pi)$.

\sk
In the remainder of this section we relate $\area$ and $\dinv$ by means of the zeta map.

\begin{mythrm}{dinvC} Let $\pi\in\mathcal{L}_{n,n}$.
Then $\dinv_C(\pi)=\area(\zeta_C(\pi))$.
\end{mythrm}

\begin{myproof} Let $\pi\in\mathcal{L}_{n,n}$ be a lattice path with type $C_n$ area vector $\mu$ and set $\beta=\zeta_C(\pi)$.
For $k\in\{0,\dots,n\}$ let $\alpha_k$ be the number of indices $i\in[n]$ such that $\abs{\mu_i}=k$. Recall the definition of the bounce path from the proof of \reft{zetaC}.

The area below bounce path of $\zeta_C(\pi)$ is given by
\begin{align*}
\alpha_0^2+\sum_{k=1}^n\binom{\alpha_k}{2}.
\end{align*}
This accounts exactly for the number of diagonal inversions of type $\mu_{n-i+1}=\mu_{n-j+1}$, $\mu_{n-i+1}=-\mu_{n-j+1}$ or $\mu_{n-i+1}=0$.
Let $(r,s)$ be the top right corner of a $1\times1$ square between $\zeta_C(\pi)$ and its bounce path.
Then the $r$-th East step and the $s$-th North step of $\beta$ both exist and belong to the segment $\LS_k^-(\mu)\RS_k^+(\mu)$ for some $k\in\{0,\dots,n\}$.
There are three cases:
If the $r$-East step and the $s$-th North step of $\beta$ both belong to $\LS_k^-(\mu)$ then there exist indices $i,j\in[n]$ with $i<j$, $\mu_{n-j+1}=-k-1$ and $\mu_{n-i+1}=-k$.
Thus $(i,j)$ is a diagonal inversion of type $\mu_{n-i+1}=\mu_{n-j+1}+1$, where $\mu_{n-i+1}\leq0$.
If $r$-th East step and the $s$-th North step of $\beta$ both belong to $\RS_k^+(\mu)$ then there exist indices $i,j\in[n]$ with $i<j$, $\mu_{n-j+1}=k$ and $\mu_{n-i+1}=k+1$.
Here $(i,j)$ is a diagonal inversion of type $\mu_{n-i+1}=\mu_{n-j+1}+1$, where $\mu_{n-i+1}>0$.
Finally, if the $s$-th North step of $\beta$ belongs to $\LS_k^-(\mu)$ and the $r$-th East step of $\beta$ belongs to $\RS_k^+(\mu)$ then there exist indices $i,j\in[n]$ such that $\mu_i=-k$ and $\mu_j=k+1$.
In this case either $(i,j)$ or $(j,i)$ is a diagonal inversion of type $\mu_{n-i+1}=-\mu_{n-j+1}+1$.

We conclude that there is a bijective correspondence between diagonal inversion of $\pi$ and the $1\times1$ squares below $\zeta_C(\pi)$.
\end{myproof}

\begin{mythrm}{dinvC'} Let $(\pi,u)\in\op{Vert}(C_n)$.
Then $\dinv_C'(\pi,w)=\area'(\zeta_C(\pi,w))$.
\end{mythrm}

\begin{myproof} Let $(\pi,u)\in\op{Vert}(C_n)$, let $\mu$ be the type $C_n$ area vector of $\pi$ and set $(\beta,w)=\zeta_C(\pi,u)$
By definition each step of $\beta$ corresponds to an entry of $\mu$.
More precisely, let $i,r\in[n]$. Then the $r$-th East step (the $r$-th North step) of $\beta$ corresponds to the entry $\mu_{n-i+1}$ if and only if
\begin{align*}
r
&=\#\{j\in[n]:|\mu_{n-j+1}|>|\mu_{n-i+1}|\}\\
&\quad+\#\{j\in[n]:\mu_{n-i+1}>0,\mu_{n-j+1}=-\mu_{n-i+1}\}\\
&\quad+\#\{j\in[i]:\mu_{n-i+1}\leq0,\mu_{n-j+1}=\mu_{n-i+1}\}\\
&\quad+\#\{j\in[n]:\mu_{n-i+1}>0,\mu_{n-j+1}=\mu_{n-i+1},j\geq i\}.
\end{align*}
Comparing this to the discussion of the diagonal reading word, see that the label below the $r$-th East step of $\beta$ (to the right of the $r$-th North step of $\beta$) is given by
\begin{align*}
w(n-r+1)
&=
\begin{cases}
u(i)&\quad\text{if }\mu_{n-i+1}>0,\\
-u(i)&\quad\text{if }\mu_{n-i+1}\leq0.
\end{cases}
\end{align*}
If $r>n$ then the $r$-th North step of $\beta$ corresponds to the entry $\mu_{n-i+1}$ if and only if $\mu_{n-i+1}=0$ and
\begin{align*}
r-n
&=\#\{j\in[n]:i\leq j,\mu_{n-j+1}=0\}.
\end{align*}
In this case the label to the right of the $i$-th North step is given by $w(n-r)=u(i)$.

\sk
Recall from \reft{dinvC} that each $1\times1$ square below $\beta$ corresponds to a diagonal inversion of $\pi$.
We demonstrate the claim by showing that such a square contributes to $\area'(\beta,w)$ if and only if the corresponding diagonal inversion contributes to $\dinv_C'(\pi,u)$.

Let $(r,s)$ be a the top right corner of a $1\times1$ square below $\beta$.
Choose $j$ such that the $s$-th North step of $\beta$ corresponds to the entry $\mu_{n-j+1}$.

Suppose the square lies below the bounce path.
Assume further that $\beta$ contains at least $r$ East steps.
In this case $s\leq n$ and $\mu_{n-j+1}\neq0$.
Choose $i$ such that the $r$-th East step of $\beta$ corresponds to the entry $\mu_{n-i+1}$.

Note that $\mu_{n-i+1}\in\{-\mu_{n-j+1},\mu_{n-j+1}\}$.
First assume $\mu_{n-i+1}=\mu_{n-j+1}$.
If $\mu_{n-i+1}<0$ then $i<j$ and $w(n-r+1)=-u(i)$ and $w(n-s+1)=-u(j)$. Thus,
\begin{align*}
w(n-r+1)>w(n-s+1) \Leftrightarrow u(i)<u(j)
\end{align*}
and the square contributes to $\area'(\beta,w)$ if and only if $(i,j)$ is a diagonal inversion of $(\pi,u)$ of type $\mu_{n-i+1}=\mu_{n-j+1}$.

If $\mu_{n-i+1}>0$ then $j<i$ and $w(n-r+1)=u(i)$ and $w(n-s+1)=u(j)$. Thus,
\begin{align*}
w(n-r+1)>w(n-s+1) \Leftrightarrow u(i)>u(j)
\end{align*}
and the square contributes to $\area'(\beta,w)$ if and only if $(i,j)$ is a diagonal inversion of $(\pi,u)$ of type $\mu_{n-j+1}=\mu_{n-i+1}$.

Secondly, assume $\mu_{n-i+1}=-\mu_{n-j+1}$. Then $\mu_{n-i+1}<0<\mu_{n-j+1}$. Thus $w(n-r+1)=-u(i)$ and $w(n-s+1)=u(j)$. We obtain
\begin{align*}
w(n-r+1)>w(n-s+1) \Leftrightarrow u(i)<-u(j) \Leftrightarrow u(j)<-u(i)
\end{align*}
and the square contributes to $\area'(\beta,w)$ if and only if either $i<j$ and $(i,j)$ is a diagonal inversion of $(\pi,u)$ of type $\mu_{n-i+1}=-\mu_{n-j+1}$ or $j<i$ and $(j,i)$ is a diagonal inversion of $(\pi,u)$ of type $\mu_{n-j+1}=-\mu_{n-i+1}$.

Next, consider the case where the square is a cell below the bounce path but $\beta$ does not contain $r$ East steps.
Then $\mu_{n-j+1}=0$.
Set $i=u^{-1}(-w(n-r+1))$, that is, $w(n-r+1)=-u(i)$.
Note that $i\in[n]$ and $\mu_{n-i+1}=0$.

First suppose that $s\leq n$. Then $i<j$ and $w(n-s+1)=-u(j)$. Hence,
\begin{align*}
w(n-r+1)>w(n-s+1) \Leftrightarrow u(i)<u(j)
\end{align*}
and the square contributes to $\area'(\beta,w)$ if and only if $(i,j)$ is a diagonal inversion of $(\pi,u)$ of type $\mu_{n-i+1}=\mu_{n-j+1}$.

Secondly, assume $s>n$. Then $i\leq j$ and $w(n-s)=u(j)$. If $i<j$, we obtain
\begin{align*}
w(n-r+1)>w(n-s) \Leftrightarrow u(i)<-u(j)
\end{align*}
and the square contributes to $\area'(\beta,w)$ if and only if $(i,j)$ is a diagonal inversion of $(\pi,u)$ of type $\mu_{n-i+1}=-\mu_{n-j+1}$.

On the other hand, if $i=j$ then
\begin{align*}
w(n-r+1)>w(n-s) \Leftrightarrow u(i)<-u(i) \Leftrightarrow u(i)<0
\end{align*}
and the square contributes to $\area'(\beta,w)$ if and only if $i$ is a diagonal inversion of $(\pi,u)$ of type $\mu_{n-i+1}=0$.

\sk
Finally we come to the case where the square lies between $\beta$ and its bounce path.
Choose $i$ such that the $r$-th East step of $\beta$ corresponds to $\mu_{n-i+1}$.
Note that $\beta$ contains at least $r$ East steps.

Assume $\mu_{n-i+1}=k+1$ with $k>0$.
Then $w(n-r+1)=u(i)$ and $\mu_{n-j+1}\in\{-k,k\}$, thus in particular $\mu_{n-j+1}\neq0$. If $\mu_{n-j+1}=k$ then $i<j$ and $w(n-s+1)=u(j)$.
Hence,
\begin{align*}
w(n-r+1)>w(n-s+1) \Leftrightarrow u(i)>u(j)
\end{align*}
and the square contributes to $\area'(\beta,w)$ if and only if $(i,j)$ is a diagonal inversion of $(\pi,u)$ of type $\mu_{n-i+1}=\mu_{n-j+1}+1$.

On the other hand if $\mu_{n-j+1}=-k$ then $w(n-s+1)=-u(j)$. Hence,
\begin{align*}
w(n-r+1)>w(n-s+1) \Leftrightarrow u(i)>-u(j) \Leftrightarrow u(j)>-u(i)
\end{align*}
and the square contributes to $\area'(\beta,w)$ if and only if either $i<j$ and $(i,j)$ is a diagonal inversion of $(\pi,u)$ of type $\mu_{n-i+1}=-\mu_{n-j+1}+1$ or $j<i$ and $(j,i)$ is a diagonal inversion of $(\pi,u)$ of type $\mu_{n-j+1}=-\mu_{n-i+1}+1$.

Secondly, assume $\mu_{n-i+1}=1$. Then $w(n-r+1)=u(i)$ and $\mu_{n-j+1}=0$. If $s>n$ then $i<j$ and $w(n-s)=u(j)$. Thus,
\begin{align*}
w(n-r+1)>w(n-s+1) \Leftrightarrow u(i)>u(j)
\end{align*}
and the square contributes to $\area'(\beta,w)$ if and only if $(i,j)$ is a diagonal inversion of $(\pi,u)$ of type $\mu_{n-i+1}=\mu_{n-j+1}+1$.

On the other hand, if $s\leq n$ then $w(n-s+1)=-u(j)$. We obtain
\begin{align*}
w(n-r+1)>w(n-s+1) \Leftrightarrow u(i)>-u(j) \Leftrightarrow u(j)>-u(i)
\end{align*}
and the square contributes to $\area'(\beta,w)$ if and only if either $i<j$ and $(i,j)$ is a diagonal inversion of $(\pi,u)$ of type $\mu_{n-i+1}=-\mu_{n-j+1}+1$ or $j<i$ and $(j,i)$ is a diagonal inversion of $(\pi,u)$ of type $\mu_{n-j+1}=-\mu_{n-i+1}+1$.

Thirdly, assume $\mu_{n-i+1}=-k<0$. Then $w(n-r+1)=-u(i)$ and $\mu_{n-j+1}=-k+1\leq0$ and $j<i$.
Also $s\leq n$ and therefore $w(n-s+1)=-u(j)$.
We obtain
\begin{align*}
w(n-r+1)>w(n-s+1) \Leftrightarrow u(j)>u(i)
\end{align*}
and the square contributes to $\area'(\beta,w)$ if and only if $(j,i)$ is a diagonal inversion of $(\pi,u)$ of type $\mu_{n-j+1}=\mu_{n-i+1}+1$.
\end{myproof}

%% file: typeD.tex
\section{The zeta map of type D}\label{Section:typeD}

In Sections \ref{Section:vertD} and \ref{Section:diagD} we describe combinatorial models for the finite torus of type $D_n$ in terms of vertically labelled signed lattice paths, and for the non-nesting parking functions of type $D_n$ in terms of diagonally labelled signed ballot paths.
Again these models nicely reflect the orbit structure under the action of the Weyl group $\S_n^D$ in the sense that acting by an element of $\S_n^D$ only changes the labels, but not the underlying path.

In \refs{areaD} we recover the dominant Shi region corresponding to a signed lattice path using its area vector.
Moreover we show how the Shi region corresponding to a vertically labelled signed lattice path can be obtained using the diagonal reading word in \refs{readD}.
The chosen models allow for a natural combinatorial description of the zeta map, which is presented in \refs{zetaD}.
We prove that this construction coincides with the uniform zeta map, thus implying that it is bijective.
Moreover we show that the zeta map induces a new bijection between lattice paths and ballot paths of odd length (without signs).

\subsection{The finite torus}\label{Section:vertD}

The Coxeter number of the root system of type $D_n$ is $h=2n-2$ hence the finite torus is $T=\Q/(2n-1)\Q$. Athanasiadis~\cite[Sec~5.4]{Athan2005} demonstrated that a system of representatives for the orbits of $T$ under the action of the Weyl group is given by
\begin{align}\label{eq:T/W}
\Q\cap(2n-1)\overline{\ac}
&=\Big\{(\lambda_1,\lambda_2,\dots,\lambda_n)\in\Q:0\leq\abs{\lambda_1}\leq\lambda_2\leq\dots\leq\lambda_n\text{ and }\lambda_{n-1}+\lambda_n\leq2n-1\Big\},
\end{align}
and that the stabiliser subgroup $\op{Stab}(\lambda)\leq\S_n^D$ is generated by the simple transpositions $s_i$ for each $i\in[n-1]$ such that $\lambda_i=\lambda_{i+1}$, the reflection $s_{\tilde\alpha}$ that exchanges $\lambda_{n-1}$ and $\lambda_n$ and changes the signs of these entries if $\lambda_{n-1}+\lambda_n=2n-1$, and the simple transposition $s_0$ that exchanges $\lambda_1$ and $\lambda_2$ and changes the signs of these entries if $\lambda_1=-\lambda_2$.

\begin{mydef}{signedlatticepath} A \emph{signed lattice path} $\pi\in\ld_{n-1,n}$ is a lattice path in $\mathcal L_{n-1,n}$ except that if it begins with an East step then this East step is replaced by a signed East step from the set $\{E^+,E^-\}$. We also define a sign function on signed lattice paths by setting $\epsilon(\pi)=-1$ if $\pi$ contains $E^-$ and $\epsilon(\pi)=1$ otherwise.
\end{mydef}

\begin{figure}[ht]
\begin{center}
\begin{tikzpicture}[scale=.55]
\begin{scope}
\draw[gray](0,0)grid(1,2);
\draw[very thick](0,0)--(1,0)--(1,2);
\draw[xshift=5mm](0,0)node[circle,fill=white,inner sep=-1.3pt]{$\oplus$};
\end{scope}
\begin{scope}[xshift=3cm]
\draw[gray](0,0)grid(1,2);
\draw[very thick](0,0)--(1,0)--(1,2);
\draw[xshift=5mm](0,0)node[circle,fill=white,inner sep=-1.3pt]{$\ominus$};
\end{scope}
\begin{scope}[xshift=6cm]
\draw[gray](0,0)grid(1,2);
\draw[very thick](0,0)--(0,1)--(1,1)--(1,2);
\end{scope}
\begin{scope}[xshift=9cm]
\draw[gray](0,0)grid(1,2);
\draw[very thick](0,0)--(0,2)--(1,2);
\end{scope}
\end{tikzpicture}
\caption{The signed lattice paths in $\ld_{1,2}$.}
\label{Figure:ld2}
\end{center}
\end{figure}
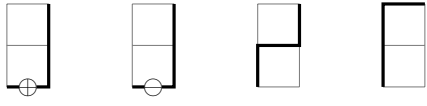

For example the set $\ld_{1,2}=\{E^+NN,E^-NN,NEN,NNE\}$ is displayed in \reff{ld2}.

Given a signed lattice path $\pi\in\ld_{n-1,n}$ let $\pi_i$ denote the number of East steps (with or without sign) of $\pi$ that occur before the $i$-th North step of $\pi$. Define an integer vector $\lambda$ by setting 
\begin{align}\notag
\lambda_1&=\epsilon(\pi)\pi_1\\ \label{eq:lambdaD}
\lambda_i&=\pi_i&&\text{for }1<i<n\text{ and}\\ \notag
\lambda_n&=
\begin{cases}
2\pi_n-\pi_{n-1}\\
2n-1-2\pi_n+\pi_{n-1}
\end{cases}
&&\begin{matrix*}[l]
\text{if }\pi_1+\dots+\pi_{n-2}\text{ is even,}\\
\text{if }\pi_1+\dots+\pi_{n-2}\text{ is odd.}
\end{matrix*}
\end{align}
\begin{myex}{lambdaD} Consider the three signed lattice paths in \reff{vertD} (ignoring the labels for the moment). For the first path we find $\pi=(0,0,4,4,4)$, $\epsilon(\pi)=1$, and $\pi_1+\pi_2+\pi_3=4$ is even. We obtain $\lambda=(0,0,4,4,4)$.

For the second path we find $\pi=(1,1,2,3,5,5)$, $\epsilon(\pi)=1$, and $\pi_1+\pi_2+\pi_3+\pi_4=7$ is odd. We compute $\lambda=(1,1,2,3,5,6)$.

Finally, for the third path we find $\pi=(3,3,3,3,3)$, $\epsilon(\pi)=-1$, and $\pi_1+\pi_2+\pi_3=9$ is odd. Hence $\lambda=(-3,3,3,3,6)$.
\end{myex}

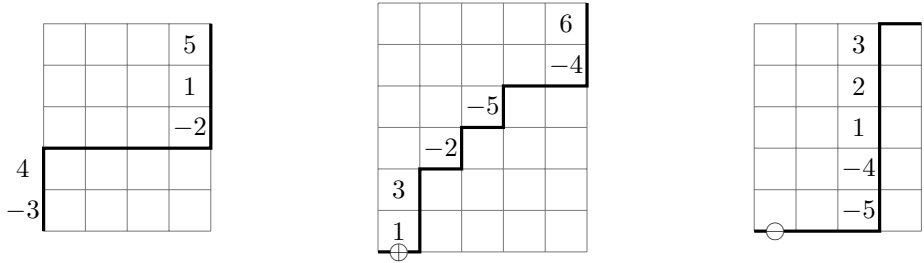
\begin{figure}[ht]
\begin{center}
\begin{tikzpicture}[scale=.55]
\begin{scope}
\draw[gray](0,0)grid(4,5);
\draw[very thick](0,0)--(0,2)--(4,2)--(4,5);
\draw[xshift=5mm,yshift=5mm]
(-1,0)node{$-3$}
(-1,1)node{$4$}
(3,2)node{$-2$}
(3,3)node{$1$}
(3,4)node{$5$};
\end{scope}
\begin{scope}[xshift=8cm,yshift=-5mm]
\draw[gray](0,0)grid(5,6);
\draw[very thick](0,0)--(1,0)--(1,2)--(2,2)--(2,3)--(3,3)--(3,4)--(5,4)--(5,6);
\draw(.5,0)node[circle,fill=white,inner sep=-1.3pt]{$\oplus$};
\draw[xshift=5mm,yshift=5mm]
(0,0)node{$1$}
(0,1)node{$3$}
(1,2)node{$-2$}
(2,3)node{$-5$}
(4,4)node{$-4$}
(4,5)node{$6$};
\end{scope}
\begin{scope}[xshift=17cm]
\draw[gray](0,0)grid(4,5);
\draw[very thick](0,0)--(3,0)--(3,5)--(4,5);
\draw(.5,0)node[circle,fill=white,inner sep=-1.3pt]{$\ominus$};
\draw[xshift=5mm,yshift=5mm]
(2,0)node{$-5$}
(2,1)node{$-4$}
(2,2)node{$1$}
(2,3)node{$2$}
(2,4)node{$3$};
\end{scope}
\end{tikzpicture}
\caption{Three vertically labelled signed lattice paths.}
\label{Figure:vertD}
\end{center}
\end{figure}

It is easy to see that signed lattice paths represent the orbits of the finite torus under the action of the Weyl group.

\begin{myprop}{lambdaD} The map $\psi:\ld_{n-1,n}\to\Q\cap(2n-1)\overline{\ac}$ defined by \refq{lambdaD} is a bijection.
\end{myprop}

\begin{myproof} Suppose $\pi\in\ld_{n-1,n}$ is a signed lattice path. Clearly $\abs{\lambda_1}\leq\lambda_2\leq\dots\leq\lambda_{n-1}$. Moreover, $\pi_{n-1}\leq2\pi_n-\pi_{n-1}<2n-1-\pi_{n-1}$ and $\pi_{n-1}<2n-1-2\pi_n+\pi_{n-1}\leq2n-1-\pi_{n-1}$ hence $\lambda_{n-1}\leq\lambda_n\leq2n-1-\lambda_{n-1}$. Since $\lambda_1+\dots+\lambda_n$ is even by definition, we conclude that $\lambda\in\Q\cap(2n-1)\overline{\ac}$.
\end{myproof}

\sk
By adding suitable labels to the signed lattice paths, we obtain a combinatorial model for the finite torus of type $D_n$. This definition is very much in the spirit of the vertically labelled Dyck paths in type $A_{n-1}$ and the vertically labelled lattice paths in type $C_n$.

\begin{mydef}{vertD} A \emph{vertically labelled signed lattice path} $(\pi,v)$ is a pair of a signed lattice path $\pi\in\ld_{n-1,n}$ and a signed permutation $v\in\S_n^B$ such that $v(i)<v(i+1)$ for each rise $i$ of $\pi$, $\abs{v(1)}<v(2)$ if $\pi$ begins with two North steps, and
\begin{align*}
\prod_{i=1}^n\sgn(v(i))&=\epsilon(\pi)(-1)^{\lambda_{n-1}+\lambda_n}\,,
\end{align*}
where $\lambda$ is defined by \refq{lambdaD}. We denote the set of all vertically labelled signed lattice paths by $\mathrm{Vert}(D_n)$.
\end{mydef}

\sk
Furthermore, given a vertically labelled signed lattice path $(\pi,v)$ define a signed permutation $u\in\S_n^D$, by setting
\begin{align}\notag
u(1)&=\epsilon(\pi)v(1),\\ \label{eq:uD}
u(i)&=v(i)&&\text{for }1<i<n\text{ and}\\ \notag
u(n)&=(-1)^{\lambda_{n-1}+\lambda_n}v(n).
\end{align}
Note that by definition $u$ has an even number of sign changes.

\begin{myex}{vertD} We can confirm that the vertically labelled lattice paths in \reff{vertD} fulfil all requirements. Consider the left path first. Since $\epsilon(\pi)(-1)^{\lambda_4+\lambda_5}=1\cdot(-1)^{4+4}=1$ we must have $v\in\S_5^D$, that is, $v$ must have an even number of sign changes. Since $\pi$ begins with two North steps we require $\abs{v(1)}=3<v(2)=4$. Additionally, $3$ and $4$ are rises of $\pi$ thus $v(3)=-2<v(4)=1<v(5)=5$. The even signed permutation $u$ is given by $u=v=[-3,4,-2,1,5]$.

Secondly, for the middle path we find $\epsilon(\pi)(-1)^{\lambda_5+\lambda_8}=1\cdot(-1)^{5+6}=-1$. Accordingly the signed permutation $v=[1,3,-2,-5,-4,6]$ has an odd number of sign changes, that is, $v\in\S_6^B-\S_6^D$. The path $\pi$ has two rises $1$ and $5$, and $v$ fulfils the conditions $v(1)=1<v(2)=3$ and $v(5)=-4<v(6)=6$. The associated even signed permutation $u=[1,3,-2,-5,-4,-6]\in\S_6^D$ is obtained from $v$ by changing the sign of $v(6)$.

Finally, for the right path we compute $\epsilon(\pi)(-1)^{\lambda_4+\lambda_5}=(-1)\cdot(-1)^{3+6}=1$. Accordingly the signed permutation $v=[-5,-4,1,2,3]\in\S_5^D$ has an even number of sign chances. The rises of $\pi$ are $1,2,3$ and $4$, and we have $v(1)<v(2)<v(3)<v(4)<v(5)$. Moreover the even signed permutation $u=[5,-4,1,2,-3]$ is obtained from $v$ by changing the signs of $v(1)$ and $v(5)$.
\end{myex}

\sk
We now extend the bijection of \refp{lambdaD} to the finite torus.

\begin{myprop}{vertD} The map $\psi:\op{Vert}(D_n)\to\Q/(2n-1)\Q$ given by $(\pi,v)\mapsto u\cdot\lambda+(2n-1)\Q$, where $\lambda$ and $u$ are defined in \refq{lambdaD} and \refq{uD}, is a bijection.
\end{myprop}

\begin{myproof} Let $\pi\in\ld_{n-1,n}$ be a signed lattice path and $u\in\S_n^D$ an even signed permutation. Using \refp{lambdaD} and \refl{canon} it suffices to show that $(\pi,v)\in\op{Vert}(D_n)$ if and only if $u\cdot\walls\subseteq\Phi^+$.

This is not a difficult task and is accomplished by distinguishing a few cases. For example consider the two simple roots in $S=\{e_2-e_1,e_2+e_1\}$. If $(\pi,v)\in\op{Vert}(D_n)$ then
\begin{align*}
S\cap\walls=\{e_2-e_1\}
&\Leftrightarrow\lambda_1=\lambda_2>0
\Leftrightarrow\pi_1=\pi_2\text{ and }\epsilon(\pi)=1\\
&\Rightarrow v(1)<v(2)\text{ and }\epsilon(\pi)=1\\
&\Rightarrow u(1)<u(2)
\Leftrightarrow u\cdot(e_2-e_1)\in\Phi^+,\\
S\cap\walls=\{e_2+e_1\}
&\Leftrightarrow\lambda_1=-\lambda_2<0
\Leftrightarrow\pi_1=\pi_2\text{ and }\epsilon(\pi)=-1\\
&\Rightarrow v(1)<v(2)\text{ and }\epsilon(\pi)=-1\\
&\Rightarrow -u(1)<u(2)
\Leftrightarrow u\cdot(e_2+e_1)\in\Phi^+,\\
S\subseteq\walls
&\Leftrightarrow\lambda_1=\lambda_2=0
\Leftrightarrow\pi\text{ begins with two North steps}\\
&\Rightarrow\abs{v(1)}<v(2)
\Leftrightarrow\abs{u(1)}<u(2)
\Leftrightarrow u\cdot S\subseteq\Phi^+.
\end{align*}
Conversely, assume that $u\cdot\walls\subseteq\Phi^+$. Then
\begin{align*}
\pi_1=\pi_2>0\text{ and }\epsilon(\pi)=1
&\Leftrightarrow\lambda_1=\lambda_2>0
\Leftrightarrow S\cap\walls=\{e_2-e_1\}\\
&\Rightarrow u\cdot(e_2-e_1)\in\Phi^+\text{ and }\epsilon(\pi)=1\\
&\Leftrightarrow u(1)<u(2)\text{ and }\epsilon(\pi)=1\\
&\Rightarrow v(1)<v(2),\\
\pi_1=\pi_2>0\text{ and }\epsilon(\pi)=-1
&\Leftrightarrow\lambda_1=-\lambda_2>0
\Leftrightarrow S\cap\walls=\{e_2+e_1\}\\
&\Rightarrow u\cdot(e_2+e_1)\in\Phi^+\text{ and }\epsilon(\pi)=-1\\
&\Leftrightarrow -u(1)<u(2)\text{ and }\epsilon(\pi)=-1\\
&\Rightarrow v(1)<v(2),\\
\pi\text{ begins with two North steps}
&\Leftrightarrow\lambda_1=\lambda_2=0
\Leftrightarrow S\subseteq\walls\\
&\Rightarrow u\cdot S\subseteq\Phi^+
\Leftrightarrow\abs{u(1)}<u(2)
\Leftrightarrow\abs{v(1)}<v(2).
\end{align*}
All other roots $\alpha\in\Delta\cup\{\tilde\alpha\}$ are treated similarly (See also the proofs of \refp{vertC} and \refp{vertB}).
\end{myproof}

\subsection{Non-nesting parking functions}\label{Section:diagD}
In this section we present an interpretation of the non-nesting parking functions of type $D_n$ in terms of labelled ballot paths of odd length.

One aspect of the root system of type $D_n$ that is different from the other infinite families ($A_{n-1}$, $B_n$ and $C_n$) is the fact that the root poset is no longer planar. For this reason its antichains are seldom associated with lattice paths in the literature. However, we have found a simple way to represent antichains by ballot paths of odd length by adding a sign to a certain East step.

\begin{mydef}{signedballotpath} A \emph{signed ballot path} $\beta\in\bd_{2n-1}$ is a ballot path with $2n-1$ steps except that if its $n$-th North step is followed by an East step, then this East step is replaced with a signed East step from the set $\{E^+,E^-\}$. We define a sign function $\epsilon:\bd_{2n-1}\to\{\pm1\}$ in the same way as for signed lattice paths, that is, $\epsilon(\pi)=-1$ if $\pi$ contains the step $E^-$ and $\epsilon(\pi)=1$ otherwise.

\end{mydef}

\sk
For example the set $\bd_3=\{NEN,NNE^+,NNE^-,NNN\}$ is pictured in \reff{bd}.

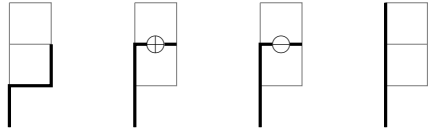
\begin{figure}[ht]
\begin{center}
\begin{tikzpicture}[scale=.55]
\begin{scope}
\ballotodd{2}
\draw[very thick] (0,0)--(0,1)--(1,1)--(1,2);
\end{scope}
\begin{scope}[xshift=3cm]
\ballotodd{2}
\draw[very thick] (0,0)--(0,2)--(1,2);
\draw(.5,2)node[circle,fill=white,inner sep=-1.3pt]{$\oplus$};
\end{scope}
\begin{scope}[xshift=6cm]
\ballotodd{2}
\draw[very thick] (0,0)--(0,2)--(1,2);
\draw(.5,2)node[circle,fill=white,inner sep=-1.3pt]{$\ominus$};
\end{scope}
\begin{scope}[xshift=9cm]
\ballotodd{2}
\draw[very thick] (0,0)--(0,3);
\end{scope}
\end{tikzpicture}
\caption{The set $\bd_{3}$ of signed ballot paths of length three.}
\label{Figure:bd}
\end{center}

\end{figure}

Let $\beta\in\bd_{2n-1}$ be a signed ballot path and $(i,j)$ a valley of $\beta$. We define the corresponding root $\alpha_{i,j}\in\Phi^+$ by
\begin{align*}
\alpha_{i,j}
&=\begin{cases}
e_{n+1-i}-e_{n+1-j}&\quad\text{if }j\leq n-1,\\
e_{n+1-i}-\epsilon(\beta)e_1&\quad\text{if }j=n,\\
e_{n+1-i}+\epsilon(\beta)e_1&\quad\text{if }j=n+1,\\
e_{n+1-i}+e_{j-n}&\quad\text{if }j\geq n+2.
\end{cases}
\end{align*}
Furthermore set
\begin{align*}
A_{\beta}=\big\{\alpha_{i,j}:(i,j)\text{ is a valley of }\beta\big\},
\end{align*}
except if $\beta$ has a valley $(i,n)$ and the $n$-th North step of $\beta$ is not followed by an East step, in which case both $\alpha_{i,n}$ and $\alpha_{i,n+1}$ are added to $A_{\beta}$ due to the valley $(i,n)$. It is easy to check that this correspondence connects signed ballot paths to antichains in the root poset.

\begin{myprop}{AbetaD} The map $\varphi:\bd_{2n-1}\to\op{Antichains}(\Phi^+)$ given by $\beta\mapsto A_{\beta}$ is a bijection between signed ballot paths and antichains in the root poset $\Phi^+$ of type $D_n$.
\end{myprop}
\noproof{}

\begin{myex}{AbetaD} Consider the signed ballot paths in \reff{diagD} (disregarding the labels for the moment). The first path has three valleys $(1,3),(2,5)$ and $(3,7)$. The second valley, marked with a star in \reff{diagD}, is of the special form $(i,n)$ where the $n$-th North step is not followed by an East step, and therefore contributes two roots. We obtain $A_{\beta}=\{e_5-e_3,e_4-e_1,e_4+e_1,e_3+e_2\}$.

The second path has two valleys $(1,7)$ and $(2,10)$. Since $7=n+1$ the former contributes the root $e_6+\epsilon(\beta)e_1=e_6-e_1$. We obtain $A_{\beta}=\{e_6-e_1,e_5+e_4\}$.

For the third path we compute $A_{\beta}=\{e_5-e_4,e_4-e_2,e_3+e_1,e_2-e_1\}$.
\end{myex}

\sk
Next we extend signed ballot paths (introducing diagonal labellings) to obtain a combinatorial model for the non-nesting parking functions of type $D_n$. This construction is akin to diagonally labelled Dyck path in type $A_{n-1}$ and diagonally labelled ballot paths in type $C_n$.

\begin{mydef}{DiagD} A \emph{diagonally labelled signed ballot path} $(\beta,w)$ is a pair of a signed ballot path $\beta\in\bd_{2n-1}$ and a signed permutation with an even number of sign changes $w\in\S_n^D$ such that for each valley $(i,j)$ of $\beta$ we have 
\begin{align*}
w(n+1-i)>
\begin{cases}
w(n+1-j) &\quad\text{if }j\leq n-1,\\
\epsilon(\beta)w(1) &\quad\text{if }j=n,\\
-\epsilon(\beta)w(1) &\quad\text{if }j=n+1,\\
w(n-j) &\quad\text{if }j\geq n+2,
\end{cases}
\end{align*}
and such that $w(n+1-i)>\abs{w(1)}$ if $\beta$ has a valley $(i,n)$ and the $n$-th North step of $\beta$ is not followed by an East step. We denote the set of all diagonally labelled signed ballot paths by $\op{Diag}(D_n)$.
\end{mydef}

\sk
Hence, if we place the labels $w(i)$, where $i=n,n-1,\dots,2,\epsilon(\beta),-\epsilon(\beta),-2,\dots,-n$, in the diagonal then for each valley the label to its right has to be smaller than the label below it.

\begin{figure}[ht]
\begin{center}
\begin{tikzpicture}[scale=.6]
\begin{scope}
\ballotodd{5}
\draw[very thick](0,0)--(0,2)--(1,2)--(1,4)--(2,4)--(2,6)--(3,6);
\draw[xshift=5mm,yshift=5mm]
(0,0)node{$2$}
(1,1)node{$4$}
(2,2)node{$-1$}
(3,3)node{$5$}
(4,4)node{$-3$}
(4,5)node{$3$}
(3,6)node{$-5$}
(2,7)node{$1$}
(1,8)node{$-4$}
(0,9)node{$-2$}
(0,2)node{$\bullet$}
(1,4)node{$*$}
(2,6)node{$\circ$};
\end{scope}
\begin{scope}[xshift=7cm,yshift=-1cm]
\ballotodd{6}
\draw[very thick](0,0)--(0,6)--(1,6)--(1,9)--(2,9);
\ballotsign{0}{6}{\ominus}
\draw[xshift=5mm,yshift=5mm]
(0,0)node{$-1$}
(1,1)node{$4$}
(2,2)node{$6$}
(3,3)node{$-5$}
(4,4)node{$-2$}
(5,5)node{$3$}
(5,6)node{$-3$}
(4,7)node{$2$}
(3,8)node{$5$}
(2,9)node{$-6$}
(1,10)node{$-4$}
(0,11)node{$1$}
(0,6)node{$\bullet$}
(1,9)node{$\circ$};
\end{scope}
\begin{scope}[xshift=15cm]
\ballotodd{5}
\draw[very thick](0,0)--(0,1)--(1,1)--(1,3)--(2,3)--(2,4)--(3,4)--(3,5)--(4,5);
\ballotsign{3}{5}{\ominus}
\draw[xshift=5mm,yshift=5mm]
(0,0)node{$5$}
(1,1)node{$4$}
(2,2)node{$3$}
(3,3)node{$-1$}
(4,4)node{$2$}
(4,5)node{$-2$}
(3,6)node{$1$}
(2,7)node{$-3$}
(1,8)node{$-4$}
(0,9)node{$-5$}
(0,1)node{$\bullet$}
(1,3)node{$\bullet$}
(2,4)node{$\bullet$}
(3,5)node{$\circ$};
\end{scope}
\end{tikzpicture}
\caption{Three diagonally labelled signed ballot paths.}
\label{Figure:diagD}
\end{center}
\end{figure}
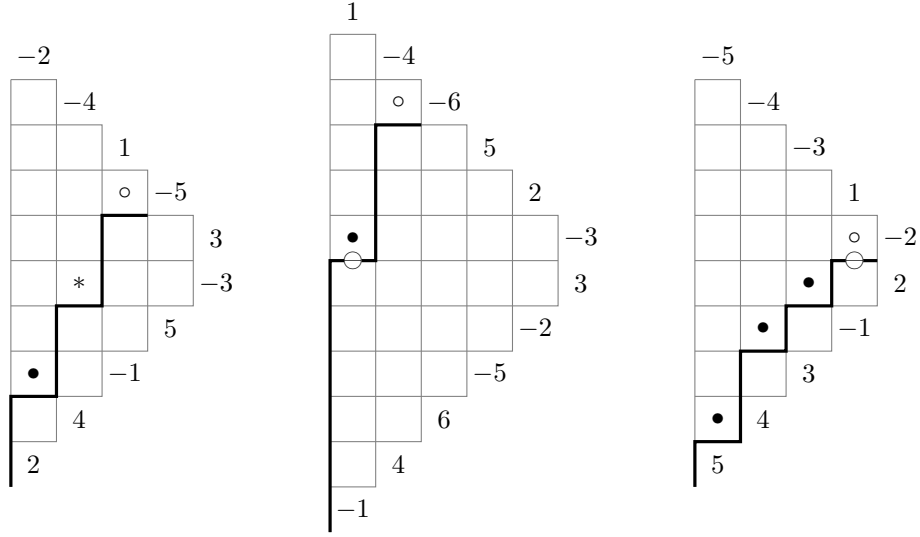

\begin{myex}{diagD} Let us return to the diagonally labelled ballot paths of \reff{diagD}. First consider the leftmost path. From $\epsilon(\beta)=1$ we obtain $w=[-3,5,-1,4,2]\in\S_5^D$. We can confirm that $w$ fulfils the conditions $w(6-2)=2>w(6-3)=-1$ and $w(6-3)=-1>w(5-7)=-5$, which are imposed by the valleys $(1,3)$ and $(3,7)$. Moreover we have $w(6-2)=4>\abs{w(1)}=3$ in accordance with the special valley $(2,5)$.

For the middle path we find $\epsilon(\beta)=-1$ and $w=[-3,-2,-5,6,4,-1]\in\S_6^D$. Note that $w(7-1)=-1>-\epsilon(\beta)w(1)=w(1)=-3$ as it is required by the presence of a valley $(1,7)=(1,n+1)$.

In case of the rightmost path we find $\epsilon(\beta)=-1$ and thus $w=[-2,-1,3,4,5]$. We verify that $w(5)=5>w(4)=4>w(2)=-1>-\epsilon(\beta)w(1)=-2$ and $w(3)=3>\epsilon(\beta)w(1)=2$.
%
%
\end{myex}

\sk
We conclude this section by extending \refp{AbetaD} to a bijection between diagonally labelled signed ballot paths an non-nesting parking functions.

\begin{myprop}{diagD} The map $\varphi:\op{Diag}(D_n)\to\op{Park}(D_n)$ given by $(\beta,w)\mapsto[w,A_{\beta}]$ is a bijection.
\end{myprop}

\begin{myproof} Let $(\beta,w)\in\op{Diag}(D_n)$ be a diagonally labelled signed ballot path. In order to prove that $\varphi$ is injective, it suffices to show that $w$ is the canonical coset representative of $[w,A_{\beta}]$ described in \refl{nnpark}, that is, $w\cdot A_{\beta}\subseteq\Phi^+$. This is a simple consequence of the definitions, for example
\begin{align*}
\alpha_{i,n}=e_{n+1-i}-\epsilon(\beta)e_1\in A_{\beta}
&\Leftrightarrow(i,n)\text{ is a valley of }\beta\\
&\Rightarrow w(n+1-i)>\epsilon(\beta)w(1)
\Leftrightarrow w(\alpha_{i,n})\in\Phi^+.
\end{align*}
Conversely suppose that $[w,A_{\beta}]\in\op{Park}(D_n)$, where $w\in\S_n^D$ is chosen such that $w\cdot A_{\beta}\subseteq\Phi^+$. Then it is not difficult to show that $(w,\beta)$ is a diagonal labelling, and thus $\varphi$ is surjective. For example,
\begin{align*}
(i,n)\text{ is a valley of }\beta
&\Leftrightarrow\alpha_{i,n}\in A_{\beta}\\
&\Rightarrow w(\alpha_{i,n})\in\Phi^+
\Leftrightarrow w(n+1-i)>\epsilon(\beta)w(1).
\end{align*}
\end{myproof}

\subsection{The area vector}\label{Section:areaD}
Let $\pi\in\ld_{n-1,n}$ be a signed lattice path and define $\lambda$ as in \refp{lambdaD}. Recall that $\lambda$ represents an orbit of the finite torus under the action of the Weyl group, and therefore corresponds to a dominant Shi region with minimal alcove $\waf_D\ac$ via the Anderson map of \refs{torus}. The main goal of this section is to recover the element $\waf_D$ of the affine Weyl group from $\pi$.

As before we write $\waf_D^{-1}=t_{\mu}\sigma$ where $\mu\in\Q$ and $\sigma\in\S_n^D$. Since $\waf_D^{-1}$ is Gra{\ss}mannian, the signed permutation $\sigma$ is given by \refl{sigma} once we know $\mu$. We will see shortly that $\mu$ can be interpreted as an area vector of the signed lattice path $\pi$, similar to the area vector of a Dyck path in type $A_{n-1}$ or a lattice path in type $C_n$.

Recall from \refs{uniform} that $\mu$ appears implicitly in the identity $\lambda=\waf_f\waf_D^{-1}\cdot0$. To find $\mu$ we need an explicit description of $w_f$.

\begin{mylem}{wfD}We have $w_f=t_{\nu}\tau$, where $\nu\in\Q$ and $\tau\in\S_n^D$ are given by
\begin{align*}
\nu=
\begin{cases}
(0,1,\dots,n-2,n-1)\\
(0,1,\dots,n-2,n)
\end{cases}
&&
\tau=
\begin{cases}
[1,2,\dots,n-1,n]\\
[-1,2,\dots,n-1,-n]
\end{cases}
&&
\begin{matrix*}[l]
\text{if }n-1\equiv0,3\mod4,\\
\text{if }n-1\equiv1,2\mod4.
\end{matrix*}
\end{align*}
\end{mylem}

\begin{myproof} By \refl{wf} it suffices to show that $t_{\nu}\tau(\widetilde\Delta)=(\Delta-\delta)\cup\{-\tilde\alpha+h\delta\}$.
Suppose $n-1\equiv0,3$ modulo $4$, then
\begin{align*}
t_{\nu}\tau(\alpha_i)
&=\alpha_i-\skal{\tau,\alpha_i}\delta
=\alpha_i-\delta
&&\text{for all }0\leq i\leq n-1,\\ 
t_{\nu}\tau(-\tilde\alpha+\delta)
&=-\tilde\alpha+(1+\skal{\tau,\tilde\alpha})\delta
=-\tilde\alpha+(2n-2)\delta.
\end{align*}
If $n-1\equiv1,2$ modulo $4$, then four cases are different:
\begin{align*}
t_{\nu}\tau(\alpha_0)
&=\alpha_1-\delta,\\
t_{\nu}\tau(\alpha_1)
&=\alpha_0-\delta,\\
t_{\nu}\tau(\alpha_{n-1})
&=-\tilde\alpha+\skal{\tau,\tilde\alpha}\delta
=-\tilde\alpha+(2n-2)\delta
&&\text{and}\\
t_{\nu}\tau(-\tilde\alpha+\delta)
&=\alpha_{n-1}+(1-\skal{\tau,\alpha_{n-1}})\delta
=\alpha_{n-1}-\delta.
\end{align*}
%
%
\end{myproof}

\sk
For example in dimensions $n=5$ and $n=6$ we have
\begin{align*}
w_f=[1,-9,-19,-29,-39]\qquad\text{and}\qquad
w_f=[-1,-11,-23,-35,-47,\,72].
\end{align*}

\begin{mydef}{areaD} Let $\pi\in\ld_{n-1,n}$ be a signed lattice path, and define $\lambda$ as in \refp{lambdaD}, and $\nu\in\Q$ and $\tau\in\S_n^D$ as in \refl{wfD}. We define the type $D_n$ \emph{area vector} of $\pi$ as $\mu=\tau(\lambda-\nu)$.

\end{mydef}

\sk
Indeed note that for $1<i<n$ the entry $\mu_i$ counts the number of boxes in the $i$-th row between the path $\pi$ and the alternating path $\rho=N(EN)^{n-1}\in\mathcal L_{n-1,n}$ (the number being negative while $\pi$ is above $\rho$). Furthermore $\mu_1$ counts the number of such boxes in the first row up to a sign, while $\mu_n$ is a little mysterious if one only looks at the picture of $\pi$.

\begin{figure}[ht]
\begin{center}
\begin{tikzpicture}[scale=.55]
\begin{scope}
\fill[black!20]
(0,1)rectangle(1,2)
(2,2)--(4,2)--(4,4)--(3,4)--(3,3)--(2,3)--cycle;
\draw[gray](0,0)grid(4,5);
\draw[xshift=0pt,yshift=-0pt,very thick](0,0)--(0,2)--(4,2)--(4,5);
\end{scope}
\begin{scope}[xshift=8cm,yshift=-5mm]
\fill[black!20]
(0,0)rectangle(1,1)
(4,4)rectangle(5,5);
\draw[gray](0,0)grid(5,6);
\draw[xshift=0pt,yshift=-0pt,very thick](0,0)--(1,0)--(1,2)--(2,2)--(2,3)--(3,3)--(3,4)--(5,4)--(5,6);
\draw(.5,0)node[circle,fill=white,inner sep=-1.3pt]{$\oplus$};
\end{scope}
\begin{scope}[xshift=17cm]
\fill[black!20]
(0,0)--(3,0)--(3,3)--(2,3)--(2,2)--(1,2)--(1,1)--(0,1)--cycle
(3,4)rectangle(4,5);
\draw[gray](0,0)grid(4,5);
\draw[xshift=0pt,yshift=-0pt,very thick](0,0)--(3,0)--(3,5)--(4,5);
\draw(.5,0)node[circle,fill=white,inner sep=-1.3pt]{$\ominus$};
\end{scope}
\end{tikzpicture}
\caption{Signed lattice paths with area vectors $(0,-1,2,1,0)$, $(-1,0,0,0,1,0)$ and $(-3,2,1,0,2)$.}
\label{Figure:areaD}
\end{center}
\end{figure}
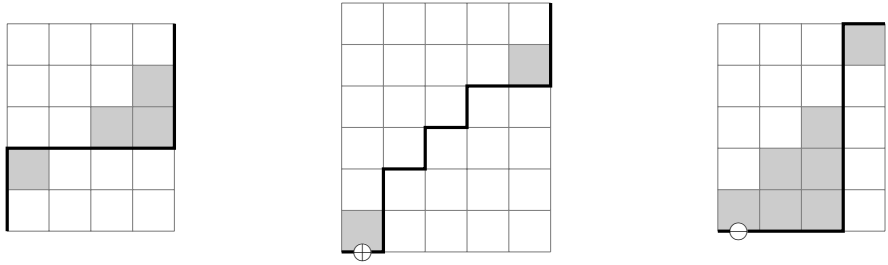

\begin{myex}{areaD} We revisit the signed lattice paths from our running example that are drawn again in \reff{areaD}. The first path begins with a North step, thus the entries $\mu_i$ of the area vector count the number of shaded boxes between $\pi$ and $\rho=N(EN)^4$ for all $i\in[n-1]$. Here $\mu_2=-1$ is negative because $\pi$ is above $\rho$ in the second row. Moreover $n-1=4\equiv0$ modulo $4$. Thus $\mu_n=\lambda_n-(n-1)=4-4=0$. We obtain $\mu=(0,-1,2,1,0)$.

For the second path we find $\epsilon(\pi)=1$ but $n-1=5\equiv1$ modulo $4$. Hence $\tau(1)=-1$ and $\mu_1=-1$. Note that this is the number of boxes between $\pi$ and $\rho=N(EN)^5$ in the first row with a negative sign even though $\pi$ is below $\rho$ in this row. For $1<i<6$ the entries $\mu_i$ give the number of boxes between $\pi$ and $\rho$ in the $i$-th row. Furthermore $\mu_n=-(n-\lambda_n)=-(6-6)=0$. We obtain $\mu=(-1,0,0,0,1,0)$.

In case of the third path we have $\epsilon(\pi)=-1$ and $n-1\equiv0$ modulo $4$. Thus $\mu_1=-3$ is the number of boxes between $\pi$ and the alternating path $\rho$ in the first row with the ``wrong'' sign. Moreover $\mu_n=\lambda_n-(n-1)=6-4=2$. We compute $\mu=(-3,2,1,0,2)$.

\sk
From the area vector we can compute the affine permutation $\waf_D$ that takes the fundamental alcove to the minimal alcove of the dominant Shi region corresponding to the signed path $\pi$. We do so for the middle path of \reff{areaD}. The translation $t_{\mu}\in\affS_n^D$ is given by $t_{\mu}=[14,2,3,4,-8,6]$. The unique even signed permutation $\sigma\in\S_n^D$ such that $(\sigma t_{\mu})^{-1}\ac$ lies in the dominant chamber is $\sigma=[-2,3,4,6,-5,1]$. We obtain $\waf_D^{-1}=t_{\mu}\sigma=[-2,3,4,6,8,14]$ and $\waf_D=[-7,-1,2,3,8,4]$.
%
%

\end{myex}

\sk
We conclude this section by proving a simple lemma that will be useful later on. The description of the representatives for the orbits of the Weyl group action on the finite torus  given in \refq{T/W} imposes certain restrictions on the area vector of a signed lattice path. Our lemma captures some of these properties.

\begin{mylem}{muD} Let $\pi\in\ld_{n-1,n}$ be a signed lattice path with area vector $\mu$.
\myi{(i)} Let $i,j\in[n]$ with $i<j$ such that $\mu_j=\mu_i-1$ and $\mu_{\ell}\notin\{\mu_i-1,\mu_i\}$ for all $\ell$ with $i<\ell<j$. Then $j=i+1$.
\myi{(ii)} Let $i\in[n]$ such that $\mu_i\leq0$ and $\mu_{\ell}\notin\{\mu_i,\mu_i+1\}$ for all $\ell$ with $1\leq\ell<i$. Then $i=1$ or $i=2,\mu_2=0$.
\myi{(iii)} If $\mu_1<0$ then there exists $i\in[n]$ with $\mu_i\in\{-\mu_1-1,-\mu_1\}$.
\myi{(iv)} Assume $\mu_1<0$ and let $i\in[n]$ such that $\mu_i=-\mu_1-1$ and $\mu_{\ell}\notin\{-\mu_1,-\mu_1-1\}$ for all $\ell$ with $1<\ell<i$. Then $i=2$.
\myi{(v)} Let $i\in[n]$ such that $\mu_i>0$ and $\mu_{\ell}\notin\{\mu_i-1,\mu_i\}$ for all $\ell$ with $i<\ell\leq n$. Then $i=n$.
\myi{(vi)} If $\mu_n>0$ then there exists $i\in[n]$ with $\mu_i\in\{-\mu_n-1,-\mu_n\}$.
\myi{(vii)} Assume $\mu_n>0$ and let $i\in[n-1]$ such that $\mu_i=-\mu_n+2$ and $\mu_{\ell}\notin\{-\mu_n+1,-\mu_n+2\}$ for all $\ell$ with $i<\ell<n$. Then $i=n-1$.
\end{mylem}

\sk
\begin{myproof} We first prove claim (i). From $\lambda=\tau\cdot\mu+\nu$ we obtain $\lambda_{\ell}=\mu_{\ell}+\ell-1$ for all $\ell$ with $1<\ell<n$. Thus, in the described situation we have $\lambda_{\ell}\notin\{\mu_i+\ell-2,\mu_i+\ell-1\}$ for all $\ell$ with $i<\ell<j$. If $n-1\equiv0,3$ modulo $4$ or if $j<n$ then $\lambda_j=\mu_i+j-2$. Hence $\lambda_{\ell}\leq\lambda_{\ell+1}$ implies $\lambda_{\ell}<\mu_i+\ell-2$ for all $\ell$ with $i<\ell<j$. On the other hand if $n-1\equiv1,2$ modulo $4$ and $j=n$, then $\lambda_n=-\mu_i+n+1$. Since $\lambda_{n-1}+\lambda_n\leq2n-1$ we obtain that $\lambda_{n-1}\leq\mu_i+n-2$, and again we have $\lambda_{\ell}<\mu_i+\ell-2$ for all $\ell$ with $i<\ell<j$. But for $\ell=i+1$ this yields a contradiction, namely $\abs{\lambda_i}=\abs{\mu_i+i-1}\leq\lambda_{i+1}<\mu_i+i-1$. Therefore we must have $j=i+1$.

To see~(ii) note that by the same argument as in the proof of (i) we obtain $\lambda_{\ell}<\mu_i+\ell-1$ for all $\ell$ with $1<\ell<i$. Thus if $i>2$ then $0\leq\lambda_2<\mu_i+1$ implies $\mu_i=0$, $\lambda_2=0$ and $\mu_1=\lambda_1=0$, which is a contradiction. If $i=2$ and $\mu_2=-1$ then $\lambda_2=0$ and $\mu_1=\lambda_1=0$, which is again a contradiction. Hence either $i=1$ or $i=2$ and $\mu_2=0$ as claimed.

To see claim~(iii) assume that $\mu_i\notin\{-\mu_1-1,-\mu_1\}$ for all $i\in[n]$. From $\lambda_i\geq\abs{\lambda_{i-1}}$ we obtain $\lambda_i>i-1-\mu_1\geq i$ for all $i$ with $1<i<n$. In particular $\lambda_{n-1}\geq n$, which is a contradiction.

Similarly in the situation of~(iv) we have $\lambda_{\ell}>\ell-1-\mu_1$ for all $\ell$ with $1<\ell<i$. If $2<i<n$, then we obtain $\lambda_{i-1}>i-2-\mu_1=\lambda_i$ which is a contradiction. If $i=n>2$ then $\lambda_{n-1}\geq n$ is the same contradiction as in the proof of (iii). Thus $i=2$ by elimination.

In the situation of~(v) we have $\lambda_{\ell}>\ell-1+\mu_i\geq\ell$ for all $\ell$ with $i<\ell<n$ by the same argument as in the proof of~(iii). Hence, if $i<n-1$ then $\lambda_{n-1}\geq n$ is a contradiction. If $i=n-1$ and $\mu_{n-1}=1$, then $\lambda_{n-1}=n-1$. Therefore $\lambda_n\in\{n-1,n\}$ and $\mu_n\in\{0,1\}$, which is a contradiction. The only remaining possibility is $i=n$.

To see~(vi) assume that $\mu_{\ell}\notin\{-\mu_n-1,-\mu_n\}$ for all $\ell\in[n]$. Then, by similar reasoning as in the proof of (i), we have $\lambda_{\ell}<\ell-\mu_n-2$ for all $\ell$ with $1<\ell<n$. This yields a contradiction for $\ell=2$.

Next consider claim~(vii). If $n-1\equiv0,3$ modulo $4$ then $2n-1-\lambda_n=n-\mu_n$. If $n-1\equiv1,2$ modulo $4$ then $\lambda_n=n-\mu_n$. In both cases we have $\lambda_{n-1}\leq n-\mu_n$. Hence $\lambda_{\ell}\notin\{\ell-\mu_n,\ell-\mu_n+1\}$ implies $\lambda_{\ell}<\ell-\mu_n$ for all $\ell$ with $i<\ell<n$. This yields a contradiction for $\ell=i+1$, namely $\lambda_{i+1}<i+1-\mu_n=\abs{\lambda_i}$. Therefore we must have $i=n-1$.
\end{myproof}



\subsection{The diagonal reading word}\label{Section:readD}

Let $(\pi,v)\in\op{Vert}(D_n)$ be a vertically labelled signed lattice path. Recall that $(\pi,v)$ corresponds to an element $u\cdot\lambda+(2n-1)\Q$ of the finite torus by \refp{vertD}, and hence to a region of the Shi arrangement with minimal alcove $\waf_R\ac$ by \reft{anderson}. Decompose $\waf_R=w\waf_D$ such that $w\in \S_n^D$ and $\waf_D\ac$ is the dominant Shi alcove corresponding to $\pi$ as in the previous section. Our aim in this section is to recover the signed permutation $w$, and thereby the element $\waf_R$ of the affine Weyl group, from $(\pi,v)$.

To this end let $u,\tau\in\S_n^D$ be defined as in \refp{vertD} and \refl{wfD} respectively. Recall that $w=u\tau\sigma$ by \refp{utausig}. Consequently one can obtain $w$ by computing $\sigma$ from $\pi$ as in the previous section. We show that $w$ can also be read off the vertically labelled signed lattice path $(\pi,v)$ in similar fashion as in type $C_n$.

\begin{mydef}{readD} Let $(\pi,v)\in\op{Vert}(D_n)$ be a vertically labelled signed lattice path with type $D_n$ area vector $\mu$. Define its type $D_n$ \emph{diagonal reading word} $d_D(\pi,v)$ as follows: For each $i=0,1,2,\dots$ first write down the labels $v(j)$ of the rows with $\mu_j=-i$ from bottom to top, then write down the negative labels $-v(j)$ of rows with $\mu_j=i+1$ from top to bottom. Finally we need to adjust some signs: Multiply the label coming from the top row by $(-1)^{1+\mu_{n-1}+\mu_n}$ and the label coming from the bottom row by $\epsilon(\pi)(-1)^{1+\nu_{n-1}+\nu_n}$. Then change the sign of $d_D(\pi,v)(1)$ if the number of positive entries of $\mu$ is odd.
\end{mydef}

\sk
Except for some necessary twists, the definition of the diagonal reading word of type $D_n$ follows very similar rules as the counterparts in types $A_{n-1}$ and $C_n$. Note that $d_D(\pi,v)\in\S_n^D$, that is, the diagonal reading word has an even number of sign changes.

\begin{myex}{readD} Let us return to the vertically labelled signed lattice paths from \reff{vertD}. 

The first path has labels $v=[-3,4,-2,1,5]$ and area vector $\mu=(0,-1,2,1,0)$. Reading off the labels according to the entries of the area vector yields the signed permutation $[-3,5,-1,4,2]$. As $(-1)^{1+\mu_{n-1}+\mu_n}=1$ and $\epsilon(\pi)(-1)^{1+\nu_{n-1}+\nu_n}=1\cdot1=1$ and since $\mu$ has two positive entries, which is a positive number, all signs remain unchanged. We obtain $d_D(\pi,v)=[-3,5,-1,4,2]$.

In case of the second path we have $v=[1,3,-2,-5,-4,6]$ and $\mu=[-1,0,0,0,1,0]$. Reading off the labels according to the entries of the area vector yields the signed permutation $[3,-2,-5,6,4,1]$. We have $(-1)^{1+\mu_{n-1}+\mu_n}=1$, however, $\epsilon(\pi)(-1)^{1+\nu_{n-1}+\nu_n}=1\cdot(-1)^{1+4+6}=-1$. Thus we need to change the sign of the label coming from the bottom row from $1$ to $-1$. Additionally $\mu$ has an odd number of positive entries, namely a single one. Consequently we change $3$ to $-3$. We obtain $d_D(\pi,v)=[-3,-2,-5,6,4,-1]$.

The third path has labels $v=[-5,-4,1,2,3]$ and area vector $\mu=(-3,2,1,0,2)$. Reading off the labels according to the entries of the area vector yields the signed permutation $[2,1,-3,-4,5]$. As $(-1)^{1+\mu_{n-1}+\mu_n}=1$ and $\epsilon(\pi)(-1)^{1+\nu_{n-1}+\nu_n}=1\cdot1=1$ and since $\mu$ has two positive entries, which is a positive number, all signs remain unchanged. We obtain $d_D(\pi,v)=[-3,5,-1,4,2]$.

\sk
In light of the next following proposition we compute that for the second path we have
\begin{align*}
d_D(\pi,v)
=[1,3,-2,-5,-4,-6]\cdot[-1,2,3,4,5,-6]\cdot[-2,3,4,6,-5,1]
=u\tau\sigma,
\end{align*}
where the respective signed permutations are taken from \refe{vertD} and \refe{areaD}.
%
\end{myex}

\sk
The following proposition asserts that our definition always yields the desired signed permutation.

\begin{myprop}{readD} Let $(\pi,v)\in\op{Vert}(D_n)$ be a vertically labelled signed lattice path with area vector $\mu$ and chose $\sigma\in\S_n^D$ as in \refl{sigma} such that $t_{\mu}\sigma$ is Gra{\ss}mannian. Moreover let $u,\tau\in\S_n^D$ be defined as in \refp{vertD} and \refl{wfD} respectively. Then $d_D(\pi,v)=u\tau\sigma$.
\end{myprop}

\begin{myproof} 
By \refl{sigma} we have $\abs{\sigma(i)}=j$ if and only if
\begin{align*}
i
&=\#\big\{k:1\leq k\leq n,\abs{\mu_k}<\abs{\mu_j}\big\}
+\#\big\{k:j\leq k\leq n,\mu_k=\mu_j>0\big\} \notag\\
&\quad+\#\big\{k:1\leq k\leq j,\mu_k=\mu_j\leq0\big\}
+\#\big\{k:1\leq k\leq n,\mu_k=-\mu_j>0\big\}
\end{align*}
Comparing this to the definition of $d_D(\pi,v)$ we obtain
\begin{align*}
\abs{d_D(\pi,v)(i)}
=\abs{v(j)}
=\abs{v\sigma(i)}
=\abs{u\tau\sigma(i)}.
\end{align*}
Comparing \refl{sigma} to the description of $u$ in \refq{uD}, $\tau$ in \refl{wfD} and the definition of the diagonal reading word, one can check that all signs work out and we may indeed drop the absolute value in the above identity.
\end{myproof}

\subsection{The zeta map}\label{Section:zetaD}

We are now in a position to present the combinatorial definition of the (Haglund--Loehr) zeta map of type $D_n$.

\begin{mydef}{zetaD} Let $\pi\in\ld_{n-1,n}$ be a signed lattice path with area vector $\mu$. Its image $\zeta_D(\pi)$ under the type $D_n$ \emph{zeta map} is obtained from
\begin{align*}
\LS_{2n-1}^-(\mu)\RS_{2n-1}^+(\mu)
\LS_{2n-2}^-(\mu)\RS_{2n-2}^+(\mu)\cdots
\LS_1^-(\mu)\RS_1^+(\mu)
\LS_0^-(\mu)\RS_0^+(\mu),
\end{align*}
by deleting the last letter and, if the $n$-th North step is followed by an East step, adding a sign to this East step such that
\begin{align*}
\epsilon(\zeta_D(\pi))=(-1)^{\#\{r\in[n]:\mu_r>0\}}.
\end{align*}
Moreover, given a vertically labelled signed lattice path $(\pi,v)\in\op{Vert}(D_n)$ we define its image under the \emph{Haglund--Loehr-zeta map} of type $D_n$ as
\begin{align*}
\zeta_D(\pi,v)=(\zeta_D(\pi),d_D(\pi,v)).
\end{align*}
\end{mydef}

\sk
Note that by definition $\zeta_D(\pi)$ is a signed ballot path, that is, $\zeta_D:\ld_{n-1,n}\to\bd_{2n-1}$.

\begin{figure}[ht]
\begin{center}
\begin{tikzpicture}[scale=.6]
\begin{scope}[yshift=2cm]
\draw[gray](0,0)grid(4,5);
\draw[very thick](0,0)--(0,2)--(4,2)--(4,5);
\draw[xshift=5mm,yshift=5mm]
(-1,0)node{$-3$}
(-1,1)node{$4$}
(3,2)node{$-2$}
(3,3)node{$1$}
(3,4)node{$5$};
\end{scope}
\draw[xshift=5cm,yshift=4cm,thick,->](0,0)--node[above]{$\zeta_D$}(1,0);
\begin{scope}[xshift=7cm]
\ballotodd{5}
\draw[very thick](0,0)--(0,2)--(1,2)--(1,4)--(2,4)--(2,6)--(3,6);
\draw[very thick,red]
(0,1)--(0,2)
(1,3)--(1,4)--(2,4)--(2,5);
\draw[xshift=5mm,yshift=5mm]
(0,0)node{$2$}
(1,1)node{$4$}
(2,2)node{$-1$}
(3,3)node{$5$}
(4,4)node{$-3$}
(4,5)node{$3$}
(3,6)node{$-5$}
(2,7)node{$1$}
(1,8)node{$-4$}
(0,9)node{$-2$}
(0,2)node{$\bullet$}
(1,4)node{$*$}
(2,6)node{$\circ$};
\end{scope}
\end{tikzpicture}
\caption{}
\label{Figure:zetaD1}
\end{center}
\end{figure}

\begin{myex}{zetaD} For a first example consider \reff{zetaD1}. Recall that $\mu=(0,-1,2,1,0)$. We obtain
\begin{align*}
\RS_2^+(\mu)=N,
\LS_1^-(\mu)=N,
\RS_1^+(\mu)=EN,
\LS_0^-(\mu)=NEN,
\RS_0^+(\mu)=NEN,
\end{align*}
which combines to $\beta=\zeta_D(\pi)=NNENNENNE$. The last North step of $\RS_0^+(\mu)$ is deleted. Note that in this example $\pi$ begins with two North steps, which translates into a valley of $\beta$ of the special form $(x,n)$, where the $n$-th North step of $\beta$ is followed by a North step. The condition $\abs{v(1)}=3<v(2)=4$ on $(\pi,v)$ corresponds to the condition $w(6-2)=4>\abs{w(1)}=3$ on $(\beta,w)$. We shall later see that this is always the case.

The signed lattice path in \reff{zetaD2} has area vector $\mu=(-1,0,0,0,1,0)$. We compute
\begin{align*}
\LS_1^-(\mu)=N,
\RS_1^+(\mu)=N,
\LS_0^-(\mu)=NNNNE,
\RS_0^+(\mu)=NNNEN.
\end{align*}
The final North step of $\RS_0^+(\mu)$ is deleted. Moreover the sixth North step is followed by an East step. Since $\mu$ has an odd number of positive entries we have $\epsilon(\beta)=-1$, and this East step is replaced by $E^-$. We conclude $\zeta_D(\pi)=NNNNNNE^-NNNE$.

The signed lattice path in \reff{zetaD3} has area vector $\mu=(-3,2,1,0,2)$. We compute
\begin{align*}
\LS_3^-(\mu)=N,
\LS_2^-(\mu)=E,
\RS_2^+(\mu)=NN
\RS_1^+(\mu)=ENE,
\LS_0^-(\mu)=N,
\RS_0^+(\mu)=EN.
\end{align*}
As above the final step of $\RS_0^+(\mu)$ is deleted, and East step following the $n$-th North step is replaced by an East step with a negative sign because $\mu$ has an odd number of positive entries. We obtain $\zeta_D(\pi)=NENNENENE^-$.
\end{myex}

\sk
In the remainder of this section we will prove that both maps defined above are the type $D_n$ instances of the uniform bijections given in \refq{zetauni} and \refq{HLzetauni}. Thereby we establish the following theorem, which is the main result of this section.

\begin{mythrm}{mainD} The Haglund--Loehr-zeta map $\zeta_D:\op{Vert}(D_n)\to\op{Diag}(D_n)$ is a bijection from vertically labelled signed lattice paths to diagonally labelled signed ballot paths.
\end{mythrm}

\sk
From \reft{mainD} we may derive some nice consequences. Since the underlying ballot path of $\zeta_D(\pi,v)$ only depends on $\pi$ and not on $v$, we have the following result.

\begin{mythrm}{zetaD} The zeta map $\zeta_D:\ld_{n-1,n}\to\bd_{2n-1}$ is a bijection from signed lattice paths to signed ballot paths.
\end{mythrm}
\noproof{}

\begin{figure}[ht]
\begin{center}
\begin{tikzpicture}[scale=.6]
\begin{scope}[yshift=2.5cm]
\draw[gray](0,0)grid(5,6);
\draw[very thick](0,0)--(1,0)--(1,2)--(2,2)--(2,3)--(3,3)--(3,4)--(5,4)--(5,6);
\ballotsign{0}{0}{\oplus}
\draw[xshift=5mm,yshift=5mm]
(0,0)node{$1$}
(0,1)node{$3$}
(1,2)node{$-2$}
(2,3)node{$-5$}
(4,4)node{$-4$}
(4,5)node{$6$};
\end{scope}
\draw[xshift=6cm,yshift=5cm,thick,->](0,0)--node[above]{$\zeta_D$}(1,0);
\begin{scope}[xshift=8cm]
\ballotodd{6}
\draw[very thick](0,0)--(0,6)--(1,6)--(1,9)--(2,9);
\draw[very thick,red]
(0,1)--(0,2)
(1,6)--(1,9)--(2,9);
\ballotsign{0}{6}{\ominus}
\draw[xshift=5mm,yshift=5mm]
(0,0)node{$-1$}
(1,1)node{$4$}
(2,2)node{$6$}
(3,3)node{$-5$}
(4,4)node{$-2$}
(5,5)node{$3$}
(5,6)node{$-3$}
(4,7)node{$2$}
(3,8)node{$5$}
(2,9)node{$-6$}
(1,10)node{$-4$}
(0,11)node{$1$}
(0,6)node{$\bullet$}
(1,9)node{$\circ$};
\end{scope}
\end{tikzpicture}
\caption{}
\label{Figure:zetaD2}
\end{center}
\end{figure}

\sk
Moreover, given any signed path $\rho$ define $\rho^*$ as the path obtained from $\rho$ by replacing all signed East steps $E^+,E^-$ by simple East steps $E$. Define $\zeta_D^*(\pi^*)=(\zeta_D(\pi))^*$ for all $\pi\in\ld_{n-1,n}$. The zeta map of type $D_n$ thereby gives rise to a new bijection between lattice paths in an $n-1\times n$ rectangle and ballot paths of odd length.

\begin{mythrm}{oddballotbij} The map $\zeta_D^*:\mathcal{L}_{n-1,n}\to\mathcal{B}_{2n-1}$ is a well-defined bijection.
\end{mythrm}

\begin{myproof} Suppose $\pi,\rho\in\ld_{n-1,n}$ differ only by the sign of the initial East step. Then their respective area vectors differ only by the sign of the first entry. It is easy to see that $\zeta_D(\pi)$ and $\zeta_D(\rho)$ can only differ by the sign of an East step. For example, suppose $a=(k,a_2,\dots,a_2)$ and $b=(-k,a_2,\dots,a_n)$, where $k>0$, then 
\begin{align*}
\LS_k^-(a)N=\LS_k^-(b)\text{ and }\RS_k^+(a)=N\RS_k^+(b).
\end{align*}
Hence
\begin{align*}
\LS_k^-(a)\RS_k^+(a)=\LS_k^-(b)\RS_k^+(b).
\end{align*}
Consequently $\zeta_D^*$ is well-defined and bijectivity follows from \reft{zetaD}.
\end{myproof}

\sk
Before we attack the proof of \reft{mainD} we need another technical lemma. In combination with \refl{sigma} it gives a criterion for when $\abs{\sigma^{-1}(i)}=1$.

\begin{mylem}{sig1D} Let $\pi$ be a signed lattice path with area vector $\mu$, and let $i\in[n]$ be such that
\begin{align*}
1&=\#\big\{r\in[n]:\abs{\mu_r\n-r}\leq\abs{\mu_i\n-i}\big\}.
\end{align*}
Then either $i$ is minimal such that $\mu_i=0$ or $\mu_j\neq0$ for all $j\in[n]$, $\mu_n=1$ and $i=n$.
\end{mylem}

\begin{myproof} Clearly $\mu_i\neq0$ implies that $\mu_{\ell}\neq0$ for all $\ell\in[n]$. Moreover $\mu_{\ell}\neq0$ for all $\ell\in[n]$ and $\mu_n=1$ implies $i=n$. Thus assume that $\mu_{\ell}\neq0$ for all $\ell\in[n]$ and $\mu_n\neq1$. From $\lambda_{n-1}\leq\min\{\lambda_n,2n-1-\lambda_n\}$ we obtain $\mu_{n-1}<0$. Thus $\mu_{\ell}<0$ for all $\ell$ with $1<\ell<n$ which yields a contradiction for $\ell=2$.
\end{myproof}

\begin{figure}[ht]
\begin{center}
\begin{tikzpicture}[scale=.6]
\begin{scope}[yshift=2cm]
\draw[gray](0,0)grid(4,5);
\draw[very thick](0,0)--(3,0)--(3,5)--(4,5);
\ballotsign{0}{0}{\ominus}
\draw[xshift=5mm,yshift=5mm]
(2,0)node{$-5$}
(2,1)node{$-4$}
(2,2)node{$1$}
(2,3)node{$2$}
(2,4)node{$3$};
\end{scope}
\draw[xshift=5cm,yshift=4cm,thick,->](0,0)--node[above]{$\zeta_D$}(1,0);
\begin{scope}[xshift=7cm]
\ballotodd{5}
\draw[very thick](0,0)--(0,1)--(1,1)--(1,3)--(2,3)--(2,4)--(3,4)--(3,5)--(4,5);
\draw[very thick,red]
(0,1)--(1,1)
(1,3)--(2,3)--(2,4)--(3,4)
(3,5)--(4,5);
\ballotsign{3}{5}{\ominus}
\draw[xshift=5mm,yshift=5mm]
(0,0)node{$5$}
(1,1)node{$4$}
(2,2)node{$3$}
(3,3)node{$-1$}
(4,4)node{$2$}
(4,5)node{$-2$}
(3,6)node{$1$}
(2,7)node{$-3$}
(1,8)node{$-4$}
(0,9)node{$-5$}
(0,1)node{$\bullet$}
(1,3)node{$\bullet$}
(2,4)node{$\bullet$}
(3,5)node{$\circ$};
\end{scope}
\end{tikzpicture}
\caption{}
\label{Figure:zetaD3}
\end{center}
\end{figure}

\sk
We proof \reft{mainD} in two steps. First we prove \reft{risevalleyD}, which establishes a strong correspondence between the rises of $(\pi,v)$ and the valleys of $\zeta_D(\pi,v)$.

Let $i$ be a rise of a vertically labelled signed lattice path $(\pi,v)$. We say the rise $i$ is labelled by $(v(i),v(i+1))$. If $\pi$ begins with two North steps and $i=1$, then we say $i$ is labelled by $(\pm\abs{v(i)},v(i+1))$ instead. Let $(i,j)$ be a valley of a diagonally labelled signed ballot path $(\beta,w)$. We say 
\begin{align*}
(i,j)\text{ is labelled by }
\begin{cases}
\big(w(n+1-i),w(n+1-j)\big)&\quad\text{ if }j<n,\\
\big(w(n+1-i),\epsilon(\beta)w(1)\big)&\quad\text{ if }j=n,\\
\big(w(n+1-i),\epsilon(\beta)w(-1)\big)&\quad\text{ if }j=n+1,\\
\big(w(n+1-i),w(n-j)\big)&\quad\text{ if }j>n+1.
\end{cases}
\end{align*}
If $j=n$ and the $n$-th North step of $\beta$ is followed by another North step, then we say $(i,j)$ is labelled by $(w(n+1-i),\pm\abs{w(1)})$ instead.

Note that if we place the labels $w(i)$, where $i=n,\dots,2,\epsilon(\beta),-\epsilon(\beta),-2,\dots,-n$, in the diagonal, then each valley is labelled by the number below it and the number to its right.

\begin{mythrm}{risevalleyD} Let $(\pi,v)\in\op{Vert}(D_n)$ be a vertically labelled signed lattice path and $a,b\in v([n])$.
Then $(\pi,v)$ has a rise labelled $(a,b)$ if and only if $\zeta_D(\pi,v)$ has a valley labelled $(b,a)$ or $(-a,-b)$.
Moreover $(\pi,v)$ has a rise labelled $(\pm\abs{a},b)$ if and only if $\zeta_D(\pi,v)$ has a valley labelled $(b,\pm\abs{a})$.
\end{mythrm}

\begin{myex}{risevalleyD} Consider once more the labelled paths in \reff{zetaD3}. The vertically labelled signed lattice path has four rises that are labelled by $(-5,-4),(-4,1),(1,2)$ and $(2,3)$. This corresponds exactly to the valleys of the diagonally labelled signed ballot path, which are labelled by $(5,4),(4,-1),(3,2)$ and $(-1,-2)$.

Moreover, the vertically labelled signed lattice path in \reff{zetaD1} has a rise labelled by $(\pm3,4)$, and its image under the zeta map has a valley labelled by $(4,\pm3)$.
\end{myex}

\begin{myproof*}{of \reft{risevalleyD} (Part 1)} We start out by demonstrating the backward implication. Therefore assume that $\zeta_D(\pi)$ has a valley $(x,y)$. Recall that $\zeta_D(\pi)$ is the concatenation of sequences $\LS_k^-(\mu)$ and $\RS_k^+(\mu)$, thus there are multiple situations in which a valley can arise: within such a sequence or at the join of two sequences. These cases, while being similar, have to be treated separately.

\myi{(1.1)} The valley $(x,y)$ appears within a sequence $\LS_k^-(\mu)$. Then there exist indices $i,j\in[n]$ with $i<j$ such that $\mu_i=-k$, $\mu_j=-k-1$ and $\mu_{\ell}\notin\{-k,-k-1\}$ for all $\ell$ with $i<\ell<j$. By \refl{muD}~(i) we have $j=i+1$.

If $i=1$ then $\lambda_2=1+\mu_2=-k\leq0$ thus $\lambda_1=\lambda_2=0$. It follows that $\pi_1=\pi_2=0$ and $i$ is a rise of $\pi$. 
If $n-1\equiv1,2$ modulo $4$ and $i=n-1$ then $\lambda_n=n-\mu_n=n+k+1$ and $\lambda_{n-1}=n-2-k$. Hence $\lambda_{n-1}+\lambda_n=2n-1$ and $i$ is a rise of $\pi$.
In all other cases $\lambda_i=i-1-k$ and $\lambda_{i+1}=i-k-1$. Thus $\pi_i=\lambda_i=\lambda_{i+1}=\pi_{i+1}$ and again $i$ is a rise of $\pi$.

The number $x$ equals the number of East steps in the sequence
\begin{align*}
\LS_{2n+1}^-(\mu)\RS_{2n-1}^+(\mu)\cdots\LS_{k+1}^-(\mu)\RS_{k+1}^+(\mu)\LS_k^-(\mu_{i+1},\dots,\mu_n).
\end{align*}
Hence,
\begin{align*}
n+1-x
&=\#\big\{r\in[n]:\abs{\mu_r}>k+1\big\}+\#\big\{r\in[n]:i+1\leq r, \mu_r=-k-1\big\}\\
&=\#\big\{r\in[n]:\abs{\mu_r\n-r}\geq\abs{(-k-1)\n-(i+1)}\big\}.
\end{align*}
Note that $\abs{\sigma^{-1}(i+1)}\neq1$ by \refl{sig1D} because $\mu_{i+1}<0$. \refl{sigma} therefore provides
\begin{align*}
n+1-x
&=\#\big\{r\in[n]:\abs{\mu_r\n-r}\leq\abs{\mu_{i+1}\n-(i+1)}\big\}\\
&=\abs{\sigma^{-1}(i+1)}
=\sigma^{-1}(i+1).
\end{align*}
We conclude that
\begin{align*}
d_D(\pi,v)(n+1-x)
&=u\tau\sigma(n+1-x)
=u\tau(i+1)
=v(i+1)
\end{align*}
because even if $i+1=n$ we have $\mu_{n-1}+\mu_n=-2k-1$, which is odd.

On the other hand, $y$ equals the number of North steps in the sequence
\begin{align*}
\LS_{2n-1}^-(\mu)\RS_{2n-1}^+(\mu)\cdots\LS_{k+1}^-(\mu)\RS_{k+1}^+(\mu)\LS_k^-(\mu_i,\dots,\mu_n).
\end{align*}
By \refl{sigma} we therefore have
\begin{align*}
n+1-y
&=\#\big\{r\in[n]:\abs{\mu_r\n-r}\leq\abs{\mu_i\n-i}\big\}
=\abs{\sigma^{-1}(i)}.
\end{align*}
We first treat the case where $\abs{\sigma^{-1}(i)}\neq1$. Note that this implies $i>1$. Then by \refl{sigma} we have $n+1-y=\sigma^{-1}(i)$ because $\mu_i\leq0$. We conclude that
\begin{align*}
d_D(\pi,v)(n+1-y)
&=u\tau\sigma(n+1-y)
=u\tau(i)
=v(i),
\end{align*}
and the valley is labelled by $(v(i+1),v(i))$.

Next assume that $\abs{\sigma^{-1}(i)}=1$, that is, $y=n$. Then $k=0$. If $i=1$ then we have already seen that $\pi$ begins with two North steps. Since $\LS_0^-(\mu)$ ends with a valley and $\RS_0^+(\mu)$ begins with a North step, we are in the special situation that the valley $(x,n)$ is labelled by $(v(2),\pm\abs{v(1)})$.

Finally if $i>1$ then $\mu_{\ell}\neq0$ for all $\ell$ with $1\leq\ell<i$. It follows that there has to be an index $\ell\in[i-1]$ such that $\mu_{\ell}\in\{1,-1\}$ and thus the $n$-th North step of $\zeta_D(\pi)$ is followed by an East step. Consequently $\zeta_D(\pi)$ contains a signed East step whose sign is determined by the number of positive entries of $\mu$. We conclude that
\begin{align*}
d_D(\pi,v)(1)
&=u\tau\sigma(1)
=(-1)^{\#\{r\in[n]:\mu_r>0\}}u\tau(i)
=\epsilon(\zeta_D(\pi))v(i),
\end{align*}
in which case the valley is labelled by $(v(i+1),v(i))$.

\myi{(1.2)} Secondly, assume that the valley $(x,y)$ arises within a sequence $\RS_k^+(\mu)$ for some $k>0$. Then there exist indices $i,j\in[n]$ with $i<j$ such that $\mu_i=k+1$, $\mu_j=k$ and $\mu_{\ell}\notin\{k,k+1\}$ for all $\ell$ with $i<\ell<j$. From \refl{muD}~(i) we obtain $j=i+1$.

If $n-1\equiv1,2$ modulo $4$ and $i=1$, then $\lambda_1=-k-1<0$ and $\lambda_2=1+k$. Hence $\pi_1=\pi_2=k+1$, $\epsilon(\pi)=-1$ and $i$ is a rise of $\pi$. Note that $i=n-1$ yields a contradiction: Either $n-1\equiv0,3$ modulo $4$, then $\lambda_{n-1}=(n-2)+(k+1)$ and $\lambda_n=n-1+k$ thus $\lambda_{n-1}+\lambda_n\geq2n$, or $n-1\equiv1,2$ modulo $4$, then $\lambda_n=n-k$ and $\lambda_{n-1}>\lambda_n$.
In all other cases
we have $\lambda_i=(i-1)+(k+1)=i+k=\lambda_{i+1}$. Thus $\pi_i=\pi_{i+1}$ and $i$ is a rise of $\pi$. Note that if $n-1\equiv0,3$ modulo $4$ and $i=1$ then $\lambda_1=k+1>0$ and $\epsilon(\pi)=1$.

Similar to the case above, $x$ equals the number of East steps in the sequence
\begin{align*}
\LS_{2n-1}^-(\mu)\RS_{2n-1}^+(\mu)\cdots\RS_{k+1}^+(\mu)\LS_k^-(\mu)\RS_k^+(\mu_1,\dots,\mu_i).
\end{align*}
Using \refl{sigma} we conclude
\begin{align*}
n+1-x
&=\#\big\{r\in[n]:\abs{\mu_r\n-r}\leq\abs{\mu_i\n-i}\big\}
=\abs{\sigma^{-1}(i)}.
\end{align*}
By \refl{sig1D} we have $\abs{\sigma^{-1}(i)}\neq1$ because $\mu_i>1$, thus \refl{sigma} yields $-\sigma^{-1}(i)=n+1-x$. We obtain
\begin{align*}
d_D(\pi,v)(n+1-x)
=u\tau\sigma(n+1-x)
=-u\tau(i)=-v(i),
\end{align*}
where we use that $\epsilon(\pi)=-1$ if and only if $n-1\equiv1,2$ modulo $4$ in the case where $i=1$.

On the other hand, $y$ equals the number of North steps in the sequence
\begin{align*}
\LS_{2n-1}^-(\mu)\RS_{2n-1}^+(\mu)\cdots\RS_{k+1}^+(\mu)\LS_k^-(\mu)\RS_k^+(\mu_1,\dots,\mu_{i+1}).
\end{align*}
It follows that
\begin{align*}
n+1-y
&=\#\big\{r\in[n]:\abs{\mu_r\n-r}\leq\abs{\mu_{i+1}\n-(i+1)}\big\}
=\abs{\sigma^{-1}(i+1)}.
\end{align*}
Suppose $y=n$, then by \refl{sig1D} we see that $\mu_{\ell}\neq0$ for all $\ell\in[n]$, $k=1$ and $i+1=n$, which is a contradiction as mentioned above. Thus $y<n$. From $\mu_{i+1}>0$ and \refl{sigma} we obtain $-\sigma^{-1}(i+1)=n+1-y$. Therefore
\begin{align*}
d_D(\pi,v)(n+1-y)
=u\tau\sigma(n+1-y)
=-u\tau(i+1)
=-v(i+1),
\end{align*}
and the valley is labelled by $(-v(i),-v(i+1))$.

\myi{(1.3)} The sequence $\LS_k^-(\mu)$ ends with an East step and the next non-empty sequence begins with a North step. If $\LS_k^-(\mu)$ ends with an East step then there exists an index $i\in[n]$ such that $\mu_i=-k-1$ and $\mu_{\ell}\notin\{-k-1,-k\}$ for all $\ell$ with $1\leq\ell<i$. By \refl{muD}~(ii) we have $i=1$. Consequently by \refl{muD}~(iii) the sequence $\RS_k^+(\mu)$ is non-empty. By assumption this means that there exists an index $j\in[n]$ such that $\mu_j=k=-\mu_1-1$ and $\mu_{\ell}\notin\{k,k+1\}$ for all $\ell$ with $1<\ell<j$. Now \refl{muD}~(iv) implies that $j=2$.

If $n-1\equiv0,3$ modulo $4$ then $\lambda_1=-k-1$ and $\lambda_2=k+1$. Hence $\pi_1=\pi_2=k+1$ and $\epsilon(\pi)=-1$. On the other hand if $n-1\equiv1,2$ modulo $4$ then $\pi_1=\lambda_1=k+1=\lambda_2=\pi_2$ and $\epsilon(\pi)=1$. In both cases $i$ is a rise of $\pi$.

\myi{(1.3.1)} Assume that $k>0$, and let $(x,y)$ be the valley under consideration. Then $x$ equals the number of East steps in the sequence
\begin{align*}
\LS_{2n-1}^-(\mu)\RS_{2n-1}^+(\mu)\cdots\RS_{k+1}^+(\mu)\LS_k^-(\mu).
\end{align*}
We deduce that 
\begin{align*}
n+1-x
=\#\big\{r\in[n]:\abs{\mu_r\n-r}\leq\abs{(-k-1)\n-1}\big\}
=\abs{\sigma^{-1}(1)},
\end{align*}
and since $\mu_1<0$, implying $\abs{\sigma^{-1}(1)}\neq1$, that $n+1-x=\sigma^{-1}(1)$. Moreover
\begin{align*}
d_D(\pi,v)(n+1-x)
&=u\tau\sigma(n+1-x)
=u\tau(1)=-v(1),
\end{align*}
because $\epsilon(\pi)=-1$ is and only if $n-1\equiv0,3$ modulo $4$.

On the other hand, $y$ equals the number of North steps in the sequence
\begin{align*}
\LS_{2n-1}^-(\mu)\RS_{2n-1}^+(\mu)\cdots\RS_{k+1}^+(\mu)\LS_k^-(\mu)\RS_k^+(\mu_1,\mu_2).
\end{align*}
As before we obtain
\begin{align*}
n+1-y
=\#\big\{r\in[n]:\abs{\mu_r\n-r}\leq\abs{k\n-2}\big\}
=\abs{\sigma^{-1}(2)}.
\end{align*}
We know that $\abs{\sigma^{-1}(2)}\neq1$ because of \refl{sig1D}. Since $\mu_2=k>0$ we conclude
\begin{align*}
d_D(\pi,v)(n+1-y)
&=u\tau\sigma(n+1-y)
=-u\tau(2)
=-v(2).
\end{align*}
\myi{(1.3.2)} Next assume that $k=0$. In this case the valley is of the form $(x,n+1)$, where $n+1-x=\sigma^{-1}(1)$ just as in (1.3.1) above. On the other hand by \refl{sig1D} we see that $\abs{\sigma^{-1}(2)}=1$. Note that $\LS_0^+(\mu)$ ending with an East step implies that the $n$-th North step of $\zeta_D(\pi)$ is followed by an East step. We obtain
\begin{align*}
d_D(\pi,v)(1)
&=u\tau\sigma(1)
=(-1)^{\#\{r\in[n]:\mu_r>0\}}u\tau(2)
=\epsilon(\zeta_D(\pi))v(2),
\end{align*}
and the valley $(x,n+1)$ is labelled by
\begin{align*}
\big(d_D(\pi,v)(n+1-x),-\epsilon(\zeta_D(\pi))d_D(\pi,v)(1)\big)
&=\big(-v(1),-v(2)\big).
\end{align*}

\myi{(1.4)} For some $k>0$ the sequence $\RS_k^+(\mu)$ ends with an East step and the next non-empty sequence begins with a North step. Then there exists an index $j\in[n]$ such that $\mu_j=k+1$ and $\mu_{\ell}\notin\{k,k+1\}$ for all $\ell$ with $j<\ell\leq n$. By \refl{muD}~(v) and~(vi) we know that $j=n$ and that $\LS_{k-1}^-(\mu)$ is non-empty. By assumption $\LS_{k-1}^-(\mu)$ begins with a North step, thus there exists an index $i\in[n]$ such that $\mu_i=-k+1$ and $\mu_{\ell}\notin\{-k,-k+1\}$ for all $\ell$ with $i<\ell\leq n$. Using \refl{muD}~(vii) we see that $i=n-1$.

If $n-1\equiv0,3$ modulo $4$ then $\lambda_{n-1}+\lambda_n=(n-k-1)+(n+k)=2n-1$. If $n-1\equiv1,2$ modulo $4$ then $\lambda_{n-1}=n-k-1=\lambda_n$. In both cases $i$ is a rise of $\pi$.

Once more let $(x,y)$ be the valley under consideration. Then $x$ is the number of East steps in the sequence
\begin{align*}
\LS_{2n-1}^-(\mu)\RS_{2n-1}^+(\mu)\cdots\LS_k^-(\mu)\RS_k^+(\mu),
\end{align*}
and therefore
\begin{align*}
n+1-x
&=\#\big\{r\in[n]:\abs{\mu_r\n-r}\leq\abs{(k+1)\n-n}\big\}
=\abs{\sigma^{-1}(n)}
\end{align*}
From $\mu_n>1$ we obtain $n+1-x=-\sigma^{-1}(n)$ and
\begin{align*}
d_D(\pi,v)(n+1-x)
=u\tau\sigma(n+1-x)
=-u\tau(n)
=v(n)
\end{align*}
because $\mu_{n-1}+\mu_n=(k+1)+(-k+1)=2$ is even.

Moreover $y$ equals the number of North steps in the sequence
\begin{align*}
\LS_{2n-1}^-(\mu)\RS_{2n-1}^+(\mu)\cdots\LS_k^-(\mu)\RS_k^+(\mu)\LS_{k-1}^-(\mu_{n-1},\mu_n).
\end{align*}
Thus,
\begin{align*}
n+1-y
&=\#\big\{r\in[n]:\abs{\mu_r\n-r}\leq\abs{(-k+1)\n-(n-1)}\big\}
=\abs{\sigma^{-1}(n-1)}.
\end{align*}
First suppose that $\abs{\sigma^{-1}(n-1)}\neq1$. Then $\mu_{n-1}\leq0$ implies $n+1-y=\sigma^{-1}(n-1)$, and we obtain
\begin{align*}
d_D(\pi,v)(n+1-y)
&=u\tau\sigma(n+1-y)
=u\tau(n-1)
=v(n-1).
\end{align*}
On the other hand if $\abs{\sigma^{-1}(n-1)}=1$ then $y=n$, $\mu_{n-1}=0$, $k=1$ and $\mu_{\ell}\neq0$ for all $\ell\in[n-2]$. We claim that the $n$-th North step of $\zeta_D(\pi)$, which belongs to $\LS_0^-(\mu)$ and corresponds to $\mu_{n-1}=0$, is followed by an East step. To see this assume that $\mu_{\ell}\notin\{0,-1\}$ for all $\ell\in[n-2]$. Then there has to be an index $\ell\in[n-2]$ with $\mu_{\ell}=1$. Thus $\RS_0^+(\mu)$ begins with an East step and this East step is not the last step of $\RS_0^+(\mu)$. Consequently, this East step is replaced by a signed East step such that $\epsilon(\zeta_D(\pi))=(-1)^{\#\{r\in[n]:\mu_r>0\}}$. From
\begin{align*}
d_D(\pi,v)(1)
&=u\tau\sigma(1)
=(-1)^{\#\{r\in[n]:\mu_r>0\}}u\tau(n-1)
=\epsilon(\zeta_D(\pi))v(n-1).
\end{align*}
we obtain that the valley is labelled by $(v(n),v(n-1))$.

\myi{(1.5)} The the valley arises within (or at the end of) the sequence $\RS_0^+(\mu)$. Note that $\RS_0^+(\mu)$ is non-empty by \refl{sig1D}, and recall that the last letter of $\RS_0^+(\mu)$ does not contribute to $\zeta_D(\pi)$.

\myi{(1.5.1)} There exist indices $i,j\in[n]$ with $i<j$ such that $\mu_i=1$, $\mu_j=0$ and $\mu_{\ell}\notin\{0,1\}$ for all $\ell$ with $i<\ell<j$. Note that it does not make a difference if the North step corresponding to $\mu_j=0$ is deleted. In this case $\zeta_D(\pi)$ ends with an East step, which is still counted as a valley. By \refl{muD}~(i) we have $j=i+1$.

If $i=1$ and $n-1\equiv1,2$ modulo $4$, then $\lambda_1=-1$ and $\lambda_2=1$. Hence $\pi_1=\pi_2=1$ and $\epsilon(\pi)=-1$. If $i=n-1$ and $n-1\equiv1,2$ modulo $4$, then $\lambda_{n-1}=n-1$, $\lambda_n=n$ and $\lambda_{n-1}+\lambda_n=2n-1$. Otherwise $\pi_i=\lambda_i=i=\lambda_{i+1}=\pi_{i+1}$. In all cases $i$ is a rise of $\pi$. Also note that $\epsilon(\pi)=1$ if $i=1$ and $n-1\equiv0,3$ modulo $4$.

Denote the present valley by $(x,y)$. Then $x$ equals the number of East steps in the sequence
\begin{align*}
\LS_{2n-1}^-(\mu)\RS_{2n-1}^+(\mu)\cdots\LS_0^-(\mu)\RS_0^+(\mu_1,\dots,\mu_i),
\end{align*}
We obtain
\begin{align*}
n+1-x
&=\#\big\{r\in[n]:\abs{\mu_r\n-r}\leq\abs{\n-i}\big\}
=\abs{\sigma^{-1}(i)}.
\end{align*}
Since $\mu_i=1$ and $\mu_{i+1}=0$ we have $\abs{\sigma^{-1}(i)}\neq1$ and $-\sigma^{-1}(i)=n+1-x$. Thus
\begin{align*}
d_D(\pi,v)(n+1-x)
&=u\tau\sigma(n+1-x)
=-u\tau(i)
=-v(i)
\end{align*}
as even if $i=1$ we have $\epsilon(\pi)=-1$ if and only if $n-1\equiv1,2$ modulo $4$.

On the other hand $y$ is the number of North steps in the sequence
\begin{align*}
\LS_{2n-1}^-(\mu)\RS_{2n-1}^+(\mu)\cdots\LS_0^-(\mu)\RS_0^+(\mu_1,\dots,\mu_{i+1}).
\end{align*}
Hence
\begin{align*}
y-n
&=\#\big\{r\in[i+1]:\mu_r=0\big\}
=\big\{r\in[n]:\abs{\mu_r\n-r}\leq i+1\big\}
=\abs{\sigma^{-1}(i+1)}.
\end{align*}
First assume that $y>n+1$, then by \refl{sigma} we have $\sigma^{-1}(i+1)=y-n$ since $\mu_{i+1}=0$. It follows that
\begin{align*}
d_D(\pi,v)(n-y)
&=-u\tau\sigma(y-n)
=-u\tau(i+1)
=-v(i+1)
\end{align*}
because even if $i+1=n$ then $\mu_{n-1}+\mu_n=1$, which is odd. Thus the valley under consideration is labelled by $(-v(i),-v(i+1))$.

If $y=n+1$ then there are exactly $n$ North steps in $\zeta_D(\pi)$ that occur before the East step corresponding to $\mu_i=1$. Hence the $n$-th North step of $\zeta_D(\pi)$ is followed by an East step, and $\zeta_D(\pi)$ contains a signed East step. Thus
\begin{align*}
d_D(\pi,v)(y)
&=u\tau\sigma(1)
=(-1)^{\#\{r\in[n]:\mu_r>0\}}u\tau(i+1)
=\epsilon(\zeta_D(\pi))v(i+1)
\end{align*}
because even if $i+1=n$ then $\mu_{n-1}+\mu_n=1$, which is odd. Thus the valley $(x,n+1)$ is labelled by $(-v(i),-v(i+1))$ in this case as well.

\myi{(1.5.2)} There exist indices $i,j\in[n]$ with $i<j$ such that $\mu_i=\mu_j=1$ and $\mu_{\ell}\notin\{0,1\}$ for all $\ell$ with $i<\ell<j$ or $j<\ell\leq n$. In this case the final two steps of $\RS_0^+(\mu)$ are East steps. The latter one is deleted and $\zeta_D(\pi)$ ends with an East step.

From \refl{muD}~(v) we obtain $j=n$, hence by \refl{muD}~(vii) we have $i=n-1$. If $n-1\equiv0,3$ modulo $4$ then $\lambda_n=n$ and $\lambda_{n-1}=n-1$ hence $\lambda_{n-1}+\lambda_n=2n-1$. On the other hand, if $n-1\equiv1,2$ modulo $4$ then $\lambda_{n-1}=\lambda_n=n-1$. In both cases $i$ is a rise of $\pi$.

Let $(x,y)$ be the valley under consideration. Then $x$ equals the number of East steps in the sequence
\begin{align*}
\LS_{2n-1}^-(\mu)\RS_{2n-1}^+(\mu)\cdots\LS_0^-(\mu)\RS_0^+(\mu_1,\dots,\mu_{n-1}).
\end{align*}
Therefore
\begin{align*}
n+1-x
&=\#\big\{r\in[n]:\abs{\mu_r\n-r}\leq n+2\big\}
=\abs{\sigma^{-1}(n-1)}.
\end{align*}
From \refl{sigma} and \refl{sig1D} we obtain $n+1-x=-\sigma^{-1}(n-1)$ because $\mu_{n-1}=1$. As above we conclude that
\begin{align*}
d_D(\pi,v)(n+1-x)
&=u\tau\sigma(n+1-x)
=-u\tau(n-1)
=-v(n-1).
\end{align*}

On the other hand $y$ equals one plus the number of North steps in $\zeta_D(\pi)$, that is,
\begin{align*}
y-n
&=1+\#\big\{r\in[n]:\mu_r=0\big\}
=\#\big\{r\in[n]:\abs{\mu_r\n-r}\leq n+1\big\}
=\abs{\sigma^{-1}(n)}.
\end{align*}
First assume that $\mu_{\ell}=0$ for some $\ell\in[n]$. Then $\abs{\sigma^{-1}(n)}\neq1$. \refl{sigma} and \refl{sig1D} therefore imply $y-n=-\sigma^{-1}(n)$, and we obtain
\begin{align*}
d_D(\pi,v)(n-y)
&=-u\tau\sigma(y-n)
=u\tau(n)
=-v(i)
\end{align*}
because $\mu_{n-1}+\mu_n=2$, which is even. The valley is thus labelled by
$(-v(n-1),-v(n))$.

On the other hand if $\mu_{\ell}\neq0$ for all $\ell\in[n]$, then $y-n=\abs{\sigma^{-1}(n)}=1$ by \refl{sig1D}. Note that the $n$-th North step of $\zeta_D(\pi)$ is followed by an East step, because $\zeta_D(\pi)$ only has $n$ North steps and ends with an East step. Hence $\zeta_D(\pi)$ contains a signed East step.
Since $\mu_n=1$, \refl{sigma} then implies that
\begin{align*}
d_D(\pi,v)(1)
&=u\tau\sigma(1)
=-(-1)^{\#\{r\in[n]:\mu_r=0\}}u\tau(n)
=\epsilon(\zeta_D(\pi))v(n),
\end{align*}
where we again use that $\mu_{n-1}+\mu_n$ is even. The valley under consideration is labelled by
\begin{align*}
\big(d_D(\pi,v)(n+1-x),-\epsilon(\zeta_D(\pi))d_D(\pi,v)(1)\big)
&=\big(-v(n-1),-v(n)\big).
\end{align*}

\myi{(Part 2)} In the second part of the proof we demonstrate the forward implication. Therefore let $i\in[n-1]$ be a rise of $\pi$. We have to show that we are in one of the five cases of Part~1 of the proof.

\myi{(2.1)} Assume $1<i<n-1$. Then $\pi_i=\lambda_i=i-1+\mu_i$ and $\pi_{i+1}=\lambda_{i+1}=i+\mu_{i+1}$, hence $\pi_i=\pi_{i+1}$ implies
\begin{align*}
\mu_i=\mu_{i+1}+1.
\end{align*}
If $\mu_i=-k\leq0$ then there is a valley in the sequence $\LS_k^-(\mu)$ and we are in case~(1.1). If $\mu_i=k+1>0$ then there is a valley in the sequence $\RS_k^+(\mu)$ and we are in case~(1.2) or~(1.5.1).

\myi{(2.2)} Next assume $i=1$. Then $\pi_1=\abs{\lambda_1}=\abs{\mu_1}$ and $\pi_2=\lambda_2=1+\mu_2$, hence $\pi_1=\pi_2$ implies $\abs{\mu_1}=\mu_2+1$. If $\mu_1\geq0$ we are in the same situation as in~(2.1). On the other hand if $\mu_1<0$ then
\begin{align*}
-\mu_1=\mu_2+1.
\end{align*}
Set $-k-1=\mu_1$ then the sequence $\LS_k^-(\mu)$ ends with an East step and the sequence $\RS_k^+(\mu)$ begins with a North step. We are therefore in case~(1.3).

\myi{(2.3)} Finally assume that $i=n-1$. Then $\pi_{n-1}=\lambda_{n-1}=n-2+\mu_{n-1}$.
\myi{(2.3.1)} Suppose $n-1\equiv0,3$ modulo $4$, thus $\lambda_n=n-1+\mu_n$. If $\pi_1+\dots+\pi_{n-2}$ is even, then $\pi_{n-1}=\pi_n$ implies $\lambda_n=2\pi_n-\pi_{n-1}=\pi_{n-1}$ and $\mu_{n-1}=\mu_n+1$ as in~(2.1). Otherwise $\pi_1+\dots+\pi_{n-2}$ is odd and $\pi_{n-1}=\pi_n$ implies $\lambda_n=2n-1-2\pi_n+\pi_{n-1}=2n-1-\pi_{n-1}$ and
\begin{align*}
\mu_{n-1}=-\mu_n+2.
\end{align*}
If $\mu_{n-1}=-k\leq0$ then the sequence $\RS_{k+1}^+(\mu)$ ends with an East step and the sequence $\LS_k^-(\mu)$ begins with a North step. We are therefore in case~(1.4). If $\mu_{n-1}=k+1>1$ then the sequence $\RS_k^+(\mu)$ ends with an East step while the sequence $\LS_{k-1}^-(\mu)$ begins with a North step. This again corresponds to case~(1.4). If $\mu_{n-1}=1$ then also $\mu_n=1$ which puts us into case (1.5.2).
\myi{(2.3.2)} Suppose $n-1\equiv1,2$ modulo $4$, thus $\lambda_n=n-\mu_n$. If $\pi_1+\dots+\pi_{n-2}$ is even, then $\pi_{n-1}=\pi_n$ implies $\lambda_n=\pi_{n-1}$ and $-\mu_{n-1}=\mu_n+2$. If $\pi_1+\dots+\pi_{n-2}$ is odd, then $\pi_{n-1}=\pi_n$ implies $\lambda_n=2n-1-\pi_{n-1}$ and $\mu_{n-1}=\mu_n+1$. Hence we are in the same situation as in~(2.3.1). 

\end{myproof*}

\sk
We now use \reft{risevalleyD} to show that the Haglund--Loehr-zeta map of type $D_n$ does indeed correspond to the specialisation of the uniform zeta map to this case. This completes the proof of \reft{mainD}.

\begin{mythrm}{HLzetaD} Let $\Phi$ be the root system of type $D_n$ with coroot lattice $\Q$ and zeta map $\zeta$, and let $\psi$ and $\varphi$ be defined as in \refp{vertD} and \refp{diagD}. Then the following diagram commutes.
\begin{center}
\begin{tikzpicture}
\draw (0,0) node(q){$\Q/(2n-1)\Q$};
\draw (4,0) node(park){$\mathrm{Park}(\Phi)$};
\draw (0,1.4) node(v){$\op{Vert}(D_n)$};
\draw (4,1.4) node(d){$\op{Diag}(D_n)$};
\draw[->] (q)--node[below]{$\zeta$}(park);
\draw[->] (v)--node[above]{$\zeta_D$}(d);
\draw[->] (v)--node[left]{$\psi$}(q);
\draw[->] (d)--node[right]{$\varphi$}(park);
\end{tikzpicture}
\end{center}
\end{mythrm}

\begin{myproof} Let $(\pi,v)\in\op{Vert}(D_n)$ and $(\beta,w)=\zeta_D(\pi,v)\in\op{Diag}(D_n)$. Furthermore, let $\lambda$ and $u$ be defined as in $\refq{lambdaD}$ and $\refq{uD}$ such that $\psi(\pi,v)=u\cdot\lambda+(2n-1)\Q$. Finally, let $\mu$ be the area vector of $\pi$, chose $\sigma\in\S_n^D$ as in \refl{sigma} such that $t_{\mu}\sigma$ is a Gra{\ss}mannian affine permutation and fix $\tau$ as in \refl{wfD}. Recall from \refp{antichain} and \refp{utausig} that it suffices to show $w=u\tau\sigma$ and $A_\beta=A$, where
\begin{align*}
A=\big\{(\tau\sigma)^{-1}\cdot\alpha:\alpha\in\Delta\cup\{-\tilde\alpha\}\text{ and }s_{\alpha}\cdot\lambda=\lambda\big\}.
\end{align*}
The first claim is immediate from \refp{readD}. In order to show the second claim, we first prove $A\subseteq A_{\beta}$. Assume $(\tau\sigma)^{-1}\cdot\alpha_i\in A$ for some $1<i<n-1$. Then $s_i\cdot\lambda=\lambda$ and $\lambda_i=\lambda_{i+1}$. Thus $\pi_i=\pi_{i+1}$ and $i$ is a rise of $\pi$. By \reft{risevalleyD} $\zeta_D(\pi,v)$ has a valley $(x,y)$ labelled either $(v(i+1),v(i))$ or $(-v(i),-v(i+1))$.
That is, either $u(i+1)=w(n+1-x)$ and
\begin{align*}
u(i)=
\begin{cases}
w(n+1-y)&\quad\text{if }y<n,\\
w(\epsilon(\beta))&\quad\text{if }y=n,\\
w(-\epsilon(\beta))&\quad\text{if }y=n+1,\\
w(n-y)&\quad\text{if }y>n+1,
\end{cases}
\end{align*}
or $-u(i)=w(n+1-x)$ and
\begin{align*}
-u(i+1)=
\begin{cases}
w(n+1-y)&\quad\text{if }y<n,\\
w(\epsilon(\beta))&\quad\text{if }y=n,\\
w(-\epsilon(\beta))&\quad\text{if }y=n+1,\\
w(n-y)&\quad\text{if }y>n+1.
\end{cases}
\end{align*}
Note that we may replace $v(i),v(i+1)$ by $u(i),u(i+1)$ by the choice of $i$.
Applying $u^{-1}$ to the above identities and using $w=u\tau\sigma$, we see that $\tau\sigma\cdot\alpha_{x,y}=\alpha_i$.
Hence $(\tau\sigma)^{-1}\cdot\alpha_i=\alpha_{x,y}\in A_{\beta}$.

Next set $S=(\tau\sigma)^{-1}\cdot\alpha_0,(\tau\sigma)^{-1}\cdot\alpha_1$. If $S\cap A=\{(\tau\sigma)^{-1}\cdot\alpha_1\}$ then $\lambda_1=\lambda_2>0$. We have $\epsilon(\pi)=1$ and $1$ is a rise of $\pi$. In particular $u(1)=v(1)$. By \reft{risevalleyD} $(\beta,w)$ has a valley $(x,y)$ labelled by $(u(2),u(1))$ or $(-u(1),-u(2))$. We obtain $\tau\sigma\cdot\alpha_{x,y}=e_2-e_1$. Consequently $(\tau\sigma)^{-1}\cdot\alpha_1=\alpha_{x,y}\in A_{\beta}$.

If $S\cap A=\{(\tau\sigma)^{-1}\cdot\alpha_0\}$ then $\lambda_1=-\lambda_2>0$. We have $\epsilon(\pi)=-1$ and $1$ is a rise of $\pi$. In particular $u(1)=-v(1)$. By \reft{risevalleyD} $(\beta,w)$ has a valley $(x,y)$ labelled by $(u(2),-u(1))$ or $(u(1),-u(2))$. We obtain $\tau\sigma\cdot\alpha_{x,y}=e_2+e_1$ and thus $(\tau\sigma)^{-1}\cdot\alpha_1=\alpha_{x,y}\in A_{\beta}$.

If $S\subseteq A$ then $\lambda_1=\lambda_2=0$. Hence $\pi$ begins with two North steps.
By \reft{risevalleyD} $(\beta,w)$ has a valley $(x,n)$ labelled by $(u(2),\pm\abs{u(1)})$.
As above we see that $\tau\sigma\cdot\{\alpha_{x,n},\alpha_{x,n+1}\}=\{\alpha_0,\alpha_1\}$.
Hence $S\subseteq A_{\beta}$.

If $(\tau\sigma)^{-1}\cdot\alpha_{n-1}\in A$ then $\lambda_{n-1}=\lambda_n$. It follows that $\pi_{n-1}=\pi_n$, that is, $n-1$ is a rise of $\pi$, and $\pi_1+\dots+\pi_{n-2}$ is even. Moreover $u(n)=v(n)$. By \reft{risevalleyD} $(\beta,w)$ has a valley $(x,y)$ labelled $(u(n),u(n-1))$ or $(-u(n-1),-u(n))$. We obtain $\tau\sigma\cdot\alpha_{x,y}=e_n-e_{n-1}$ and thus $(\tau\sigma)^{-1}\cdot\alpha_{n-1}\in A_{\beta}$.

If $(\tau\sigma)^{-1}\cdot(-\tilde\alpha)\in A$ then $\lambda_{n-1}+\lambda_n=2n-1$. It follows that $\pi_{n-1}=\pi_n$ and $\pi_1+\dots+\pi_{n-2}$ is odd. Moreover $u(n)=-v(n)$. By \reft{risevalleyD} $(\beta,w)$ has a valley $(x,y)$ labelled $(-u(n),u(n-1))$ or $(-u(n-1),u(n))$. We obtain $\tau\sigma\cdot\alpha_{x,y}=-e_{n-1}-e_n$ and thus $(\tau\sigma)^{-1}\cdot(-\tilde\alpha)\in A_{\beta}$.

\sk
To complete the proof we need to demonstrate $A_{\beta}\subseteq A$. Therefore suppose $\alpha_{x,y}\in A_{\beta}$ for some valley $(x,y)$ of $(\beta,w)$. Then by \reft{risevalleyD} $(x,y)$ is labelled either $(v(i+1),v(i))$ or $(-v(i),-v(i+1))$ for some rise $i$ of $\pi$, or by $(v(2),\pm\abs{v(1)})$ if $y=n$ and the valley is not followed by an East step.

If $u(i)=v(i)$ and $u(i+1)=v(i+1)$ then the valley $(x,y)$ is labelled by $(u(i+1),u(i))$ or $(-u(i),-u(i+1))$, and $\tau\sigma\cdot\alpha_{x,y}=\alpha_i$ for a rise $i$ of $\pi$. By similar arguments as above we see that $\alpha_{x,y}\in A$.

If $i=1$ and $u(1)=-v(1)$ then $\tau\sigma\cdot\alpha_{x,y}=\alpha_0$, and $\lambda_1=-\lambda_2$. Again we obtain $\alpha_{x,y}\in A$.

If $i=n-1$ and $u(n)=-v(n)$ then $\tau\sigma\cdot\alpha_{x,y}=-\tilde\alpha$ and $\lambda_{n-1}+\lambda_n=2n-1$. Again we obtain $\alpha_{x,y}\in A$.

Finally, if the valley is of the special form $(x,n)$, not followed by an East step and labelled by $(v(2),\pm\abs{v(1)})$ then $\tau\sigma\cdot\{\alpha_{x,n},\alpha_{x,n+1}\}=\{\alpha_0,\alpha_1\}$, $\pi$ begins with two North steps and $\lambda_1=\lambda_2=0$. We see that $\alpha_{x,n},\alpha_{x,n+1}\in A$.
Thus $A_{\beta}\subseteq A$ and the proof is complete.
\end{myproof}

%% file: typeB.tex
\section{The zeta map of type B}\label{Section:typeB}

In Sections~\ref{Section:vertB} and~\ref{Section:diagB} we present combinatorial models for the finite torus of type $B_n$ in terms of vertically labelled lattice paths and for the non-nesting parking functions in terms of diagonally labelled ballot paths.
These models are almost identical to the respective objects in type $C_n$ and are compatible with the action of the Weyl group $\S_n^B$ in the sense that acting by an element of $\S_n^B$ only changes the labels and leaves the underlying path intact.
As a nice side result we obtain an explicit $\S_n^B$ isomorphism between the finite torus of type $B_n$ and the finite torus of type $C_n$.

In \refs{areaB} we compute the dominant Shi region corresponding to a lattice path using its area vector.
In \refs{readB} we obtain the Shi region corresponding to a vertically labelled lattice path by means of the diagonal reading word.
The combinatorial zeta map, which we present in \refs{zetaB}, is closely related to the zeta map of type $D_{n+1}$.
The proof that our construction coincides with the uniform zeta map relies on the respective results in \refs{zetaD}.
This also accounts for our unusual choice of order: C-D-B.

\subsection{The finite torus}\label{Section:vertB}

The Coxeter number of type $B_n$ is $h=2n$. As always we first recall the analysis of Athanasiadis~\cite[Sec 5.3]{Athan2005}. A system of representatives for the orbits of the finite torus $\Q/(2n+1)\Q$ under the action of $\S_n^B$ is given by
\begin{align*}
\Q\cap(2n+1)\overline{A_{\circ}}
=\big\{(\lambda_1,\lambda_2,\dots,\lambda_n)\in\Q: 0\leq\lambda_1\leq\lambda_2\leq\dots\leq\lambda_n\text{ and }\lambda_{n-1}+\lambda_n\leq 2n+1\big\}.
\end{align*}
Moreover, the stabiliser $\op{Stab}(\lambda)\leq\S_n^B$ of such a $\lambda$ with respect to the action of the Weyl group on the finite torus is generated by the simple transpositions $s_i$ for each $i\in[n-1]$ such that $\lambda_i=\lambda_{i+1}$, the simple transposition $s_0$ if $\lambda_1=0$ and the reflection $s_{\tilde\alpha}$ that exchanges the last two entries and changes their signs if $\lambda_{n-1}+\lambda_n=2n+1$.

Our first aim is to find a set of lattice paths representing the orbits described above. Thus let $\pi\in\mathcal L_{n,n}$ and, as in the previous section, define $\pi_i$ to be the number of East steps preceding the $i$-th North step of $\pi$. We define an integer vector $\lambda$ setting
\begin{align}\notag
\lambda_i&=\pi_i&&\text{for }1\leq i<n\text{ and}\\ \label{eq:lambdaB}
\lambda_n&=
\begin{cases}
2\pi_n-\pi_{n-1}\\
2n+1-2\pi_n+\pi_{n-1}
\end{cases}
&&\begin{matrix*}[l]
\text{if }\pi_1+\dots+\pi_{n-2}\text{ is even,}\\
\text{if }\pi_1+\dots+\pi_{n-2}\text{ is odd.}
\end{matrix*}
\end{align} 
The following result was already used in \cite{Athan2005}.

\begin{myprop}{lambdaB} The map $\psi:\mathcal L_{n,n}\to\Q\cap(2n+1)\overline{\ac}$ given by \refq{lambdaB} is a bijection.
\end{myprop}\noproof{}

\sk
The next definition, which extends our model to the finite torus by introducing vertical labellings, highlights the similarity to type $C_n$.

\begin{mydef}{vertB} Let $\op{Vert}(B_n)=\op{Vert}(C_n)$.
\end{mydef}

\sk
Given a vertically labelled lattice path $(\pi,v)\in\op{Vert}(B_n)$ we define a signed permutation $u\in\S_n^B$ by
\begin{align}\label{eq:uB}
u(i)&=v(i)&&\text{for }1\leq i<n\text{ and}\\ \notag
u(n)&=(-1)^{\lambda_{n-1}+\lambda_n}v(n).
\end{align}

\begin{figure}[ht]
\begin{center}
\begin{tikzpicture}[scale=.6]
\begin{scope}
\draw[gray](0,0)grid(6,6);
\draw[very thick] (0,0)--(0,1)--(4,1)--(4,6)--(6,6);
\draw[xshift=5mm,yshift=5mm]
(-1,0)node{$1$}
(3,1)node{$-5$}
(3,2)node{$-4$}
(3,3)node{$2$}
(3,4)node{$3$}
(3,5)node{$6$};
\end{scope}
\begin{scope}[xshift=10cm,yshift=1cm]
\draw[gray](0,0)grid(4,4);
\draw[very thick](0,0)--(1,0)--(1,1)--(2,1)--(2,4)--(4,4);
\draw[xshift=5mm,yshift=5mm]
(0,0)node{$-1$}
(1,1)node{$-4$}
(1,2)node{$-3$}
(1,3)node{$-2$};
\end{scope}
\draw[white](16,0)--(16,1);
\end{tikzpicture}
\caption{Two vertically labelled lattice paths.}
\label{Figure:vertB}
\end{center}
\end{figure}
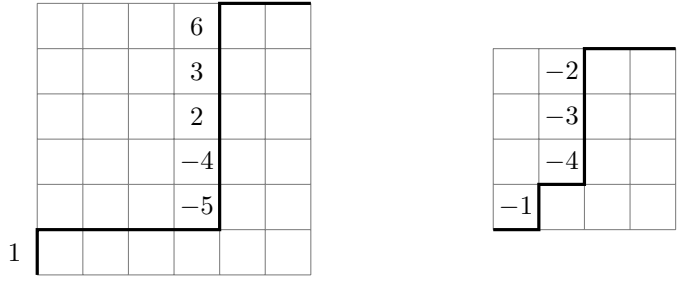

The following proposition asserts that our chosen model corresponds to the finite torus of type $B_n$ nicely.

\begin{myprop}{vertB} The map $\psi:\op{Vert}(B_n)\to\Q/(2n+1)\Q$ given by $(\pi,v)\mapsto u\cdot\lambda+(2n+1)\Q$, where $\lambda$ and $u$ are defined as in \refq{lambdaB} and \refq{uB}, is a bijection.
\end{myprop}

\begin{myproof} Let $\pi\in\mathcal L_{n,n}$ be a lattice path and $v\in\S_n^B$ a signed permutation. Using \refp{lambdaB} and \refl{canon} it suffices to show that $(\pi,v)\in\op{Vert}(B_n)$ if and only if $u\cdot\walls\subseteq\Phi^+$.

For example consider the two roots $\alpha_{n-1}$ and $\tilde\alpha$.
If $(\pi,v)\in\op{Vert}(B_n)$, then
\begin{align*}
e_n-e_{n-1}\in\walls
&\Leftrightarrow\lambda_{n-1}=\lambda_n
\Leftrightarrow\pi_{n-1}=\pi_n\text{ and }\pi_1+\dots+\pi_{n-2}\text{ is even}\\
&\Rightarrow v(n-1)<v(n)\text{ and }u(n)=v(n) \\
&\Rightarrow u(n-1)<u(n)
\Leftrightarrow u\cdot(e_n-e_{n-1})\in\Phi^+,\\
e_{n-1}+e_n\in\walls
&\Leftrightarrow\lambda_{n-1}+\lambda_n=2n+1
\Leftrightarrow\pi_{n-1}=\pi_n\text{ and }\pi_1+\dots+\pi_{n-2}\text{ is odd}\\
&\Rightarrow v(n-1)<v(n)\text{ and }u(n)=-v(n)\\
&\Rightarrow u(n-1)<-u(n)
\Leftrightarrow u\cdot(e_{n-1}+e_n)\in\Phi^+.
\end{align*}
Conversely, if $u\cdot\walls\subseteq\Phi^+$ then
\begin{align*}
\pi_{n-1}=\pi_n&\text{ and }\pi_1+\dots+\pi_{n-2}\text{ is even}
\Leftrightarrow \lambda_{n-1}=\lambda_n
\Leftrightarrow e_n-e_{n-1}\in\walls\\
&\Rightarrow u\cdot(e_n-e_{n-1})\in\Phi^+\text{ and }v(n)=u(n)
\Leftrightarrow u(n-1)<u(n)\text{ and }v(n)=u(n)\\
&\Rightarrow v(n-1)<v(n),\\
\pi_{n-1}=\pi_n&\text{ and }\pi_1+\dots+\pi_{n-2}\text{ is odd}
\Leftrightarrow \lambda_{n-1}+\lambda_n=2n+1
\Leftrightarrow e_{n-1}+e_n\in\walls\\
&\Rightarrow u\cdot(e_{n-1}+e_n)\in\Phi^+\text{ and }v(n)=-u(n)
\Leftrightarrow u(n-1)<-u(n)\text{ and }v(n)=-u(n)\\
&\Rightarrow v(n-1)<v(n).
\end{align*}
Other roots are treated similarly (See also the proof of \refp{vertC}).
\end{myproof}

\begin{myrem}{BCQ} Note that together \refp{vertC} and \refp{vertB} yield an explicit $\S_n^B$-set isomorphism between the finite torus of type $B_n$ and the finite torus of type $C_n$.
\end{myrem}

\begin{myex}{vertB} Consider the vertically labelled lattice paths of \reff{vertB}. The left path is the example we also use in \refs{typeC}. Indeed we have $\pi=(0,4,4,4,4,4)=\lambda$ just as in \refe{vertC}. Moreover, $u=v=[1,-5,-4,2,3,6]$, and the corresponding element of the finite torus $u\cdot\lambda+13\Q=(0,4,4,-4,-4,4)+13\Q$ looks the same as in type $C_6$ too.

On the other hand for the right path we compute $\pi=(1,2,2,2)$ and $\lambda=(1,2,2,7)$ which is certainly different from the $\lambda$ we use in type $C_n$. Moreover $u$ is obtained from $v=[-1,-4,-3,-2]$ by changing the sign of $v(4)$, that is, $u=[-1,-4,-3,2]$. The corresponding element of the finite torus of type $B_4$ is $u\cdot\lambda+9\Q=(-1,7,-2,-2)+9\Q$. However, note that $(-1,7,-2,-2)$ is not the canonic representative of its class in the finite torus of type $C_4$.
\end{myex}

\subsection{Non-nesting parking functions}\label{Section:diagB}

As in type $C_n$ ballot paths are the correct lattice paths to consider in this context. Let $\beta\in\mathcal B_{2n}$ be a ballot path with valley $(i,j)$ then we define a root $\alpha_{i,j}\in\Phi^+$ by
\begin{align*}
\alpha_{i,j}=
\begin{cases}
e_{n+1-i}-e_{n+1-j}&\quad\text{if }j<n+1,\\
e_{n+1-i}&\quad\text{if }j=n+1,\\
e_{n+1-i}+e_{j-n-1}&\quad\text{if }j>n+1.
\end{cases}
\end{align*}
Furthermore define
\begin{align*}
A_{\beta}=\big\{\alpha_{i,j}:(i,j)\text{ is a valley of }\beta\big\}.
\end{align*}
Then it is well-known that this correspondence identifies ballot paths with antichains in the root poset.

\begin{mylem}{AbetaB} The map $\varphi:\mathcal B_{2n}\to\op{Antichains}(\Phi^+)$ given by $\beta\mapsto A_{\beta}$ is a bijection between ballot paths of length $2n$ and the antichains in the root poset $\Phi^+$ of type $B_n$.
\end{mylem}\noproof{}

\sk
Next we define a combinatorial model for the non-nesting parking functions in terms of diagonally labelled ballot paths. In contrast to the previous section, the type $B_n$ objects are very similar to the analogous objects in type $C_n$, but not quite the same.

\begin{mydef}{diagB} A \emph{diagonally labelled ballot path} of type $B_n$ is a pair $(\beta,w)$ of a ballot path $\beta\in\mathcal B_{2n}$ and a signed permutation $w\in\S_n^B$ such that for each valley $(i,j)$ of $\beta$ we have
\begin{align*}
w(n+1-i)>w(n+1-j).
\end{align*}
We denote the set of all diagonally labelled ballot paths of type $B_n$ by $\op{Diag}(B_n)$.
\end{mydef}

\sk
These conditions are even simpler to written down than those we have encountered in types $C_n$ (or $D_n$), and can be visualised as follows: Given a ballot path $\beta$ and a signed permutation $w\in\S_n^B$ we determine if the pair $(\beta,w)$ lies in the set $\op{Diag}(B_n)$ by writing the numbers $w(i)$, where $i=n,n-1\dots,1,0,-1,\dots,-n$, in the main diagonal. Then for each valley the number below has to be greater than the number to its right.

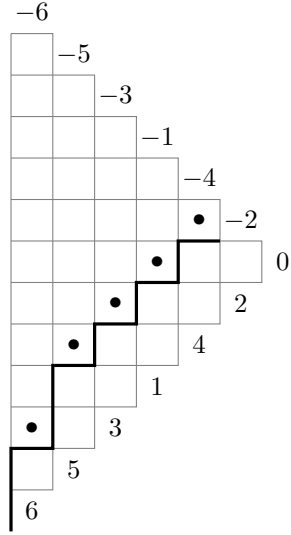
\begin{figure}[ht]
\begin{center}
\begin{tikzpicture}[scale=.55]
\begin{scope}
\ballot{6}
\draw[very thick](0,0)--(0,2)--(1,2)--(1,4)--(2,4)--(2,5)--(3,5)--(3,6)--(4,6)--(4,7)--(5,7);
\draw[xshift=5mm,yshift=5mm]
(0,0)node{$6$}
(1,1)node{$5$}
(2,2)node{$3$}
(3,3)node{$1$}
(4,4)node{$4$}
(5,5)node{$2$}
(6,6)node{$0$}
(5,7)node{$-2$}
(4,8)node{$-4$}
(3,9)node{$-1$}
(2,10)node{$-3$}
(1,11)node{$-5$}
(0,12)node{$-6$}
(0,2)node{$\bullet$}
(1,4)node{$\bullet$}
(2,5)node{$\bullet$}
(3,6)node{$\bullet$}
(4,7)node{$\bullet$};
\end{scope}
%
%
\end{tikzpicture}
\caption{A diagonally labelled ballot path $(\beta,w)\in\op{Diag}(B_6)$.}
\label{Figure:diagB}
\end{center}
\end{figure}

\begin{myex}{diagB} Consider the diagonally labelled ballot path in \reff{diagB}. The path $\beta$ has five valleys and corresponds to the antichain $A_{\beta}=\{e_6-e_4,e_5-e_2,e_4-e_1,e_3,e_2+e_1\}$. We may check that all required conditions on the signed permutation $w=[2,4,1,3,5,6]$ are fulfilled. For example $w(6+1-1)=6>w(6+1-3)=3$ and $w(6+1-4)=1>w(6+1-7)=0$ are the conditions imposed by the valleys $(1,3)$ and $(4,7)$.
\end{myex}

\begin{myprop}{diagB} The map $\varphi:\op{Diag}(B_n)\to\op{Park}(B_n)$ given by $(\beta,w)\mapsto[w,A_{\beta}]$ is a bijection.
\end{myprop}

\begin{myproof} Let $\beta\in\mathcal B_{2n}$ be a ballot path and $w\in\S_n^B$ be a signed permutation. By \refl{AbetaB} and \refl{nnpark} it suffices to show that $(\beta,w)\in\op{Diag}(B_n)$ if and only if $w\cdot A_{\beta}\subseteq\Phi^+$.

Suppose $w\cdot A_{\beta}\subseteq\Phi^+$. Then an easy case by case check reveals
\begin{align*}
(i,j)\text{ is a valley of }\beta
&\Leftrightarrow\alpha_{i,j}\in A_{\beta}\\
&\Rightarrow w\cdot\alpha_{i,j}\in\Phi^+
\Leftrightarrow w(n+1-i)>w(n+1-j).
\end{align*}
Conversely if $(\beta,w)\in\op{Diag}(B_n)$ then
\begin{align*}
\alpha_{i,j}\in A_{\beta}
&\Leftrightarrow(i,j)\text{ is a valley of }\beta\\
&\Rightarrow w(n+1-i)>w(n+1-j)
\Leftrightarrow w\cdot\alpha_{i,j}\in\Phi^+.
\end{align*}
\end{myproof}

\subsection{The area vector}\label{Section:areaB}

Let $\pi\in\mathcal L_{n,n}$ be a lattice path. By \refp{lambdaB} $\pi$ corresponds bijectively to an element $\lambda\in\Q\cap(2n+1)\overline{\ac}$. The orbit of the finite torus under the action of the Weyl group represented by $\lambda$ corresponds to a dominant region of the Shi arrangement via the Anderson map of \refs{torus} as follows: Let $\waf_D\ac$ be the minimal alcove of this region, then we have $\A(\waf_D)\in W\cdot\lambda+(2n+1)\Q$ or equivalently $\waf_f\waf_D^{-1}\cdot0=\lambda$. In this section we describe how the affine permutation $\waf_D$ can be recovered from the lattice path $\pi$.

In order to do so, write $\waf_D^{-1}=t_{\mu}\sigma$ where $\mu\in\Q$ and $\sigma\in\S_n^B$. Since $\waf_D^{-1}$ is Gra{\ss}mannian the signed permutation $\sigma$ is determined by \refl{sigma} once $\mu$ is known. On the other hand $\mu$ is related to $\lambda$ implicitly by $\lambda=\waf_f\waf_D^{-1}\cdot0$. Hence an explicit description of $\waf_f$ is needed. It is provided by the following lemma.

\begin{mylem}{wfB} Let $\waf_f=t_{\nu}\tau$ where $\nu\in\Q$ and $\tau\in\S_n^B$. Then
\begin{align*}
\nu=
\begin{cases}
(1,\dots,n-1,n) \\
(1,\dots,n-1,n+1)
\end{cases}
\qquad\tau=
\begin{cases}
[1,\dots,n-1,n] &\qquad\text{if }n\equiv0,3\mod4,\\
[1,\dots,n-1,-n] &\qquad\text{if }n\equiv1,2\mod4.
\end{cases}
\end{align*}
\end{mylem}

\begin{myproof} By \refl{wf} it suffices to show that $t_{\nu}\tau(\widetilde\Delta)=(\Delta-\delta)\cup\{-\tilde\alpha+h\delta\}$.
Suppose $n\equiv0,3$ modulo $4$, then
\begin{align*}
t_{\nu}\tau(\alpha_i)
&=\alpha_i-\skal{\tau,\alpha_i}\delta
=\alpha_i-\delta
&&\text{for all }0\leq i\leq n-1,\\ 
t_{\nu}\tau(
-\tilde\alpha+\delta)
&=-\tilde\alpha+(1+\skal{\tau,\tilde\alpha})\delta
=-\tilde\alpha+2n\delta.
\end{align*}
If $n-1\equiv1,2$ modulo $4$, then the only difference is
\begin{align*}
t_{\nu}\tau(\alpha_{n-1})
&=-\tilde\alpha+\skal{\tau,\tilde\alpha}\delta
=-\tilde\alpha+2n\delta
&&\text{and}\\
t_{\nu}\tau(-\tilde\alpha+\delta)
&=\alpha_{n-1}+(1-\skal{\tau,\alpha_{n-1}})\delta
=\alpha_{n-1}-\delta.
\end{align*}
%
\end{myproof}

\sk
For example if $n=4$, respectively $n=5$, then
\begin{align*}
w_f=[-8,-16,-24,-36],&&
w_f=[-10,-20,-30,-40,\,61].
\end{align*}

\begin{mydef}{avB} Let $\pi\in\mathcal L_{n,n}$ be a lattice path, $\lambda$ be defined as in \refp{lambdaB} and $\nu$ and $\tau$ as in \refl{wfB} above. Then we define the \emph{area vector} of type $B_n$ of $\pi$ as
\begin{align*}
\mu
&=\tau\cdot(\lambda-\nu)
=\begin{cases}
(\lambda_1-1,\dots,\lambda_{n-1}-n+1,\lambda_n-n)&\quad\text{if }n\equiv0,3\mod4,\\
(\lambda_1-1,\dots,\lambda_{n-1}-n+1,n+1-\lambda_n)&\quad\text{if }n\equiv1,2\mod4.
\end{cases}
\end{align*}
\end{mydef}

Note that similarly to the other types, the entry $\mu_i$ of the area vector counts the number of boxes in the $i$-th row that lie between the path $\pi$ and the alternating path $(EN)^n\in\mathcal L_{n,n}$, where $\mu_i$ is negative as long as $\pi$ is above $(EN)^n$. The only exception to this rule is the top row, where $\mu_n$ does not have as nice of an interpretation.

\begin{figure}[ht]
\begin{center}
\begin{tikzpicture}[scale=.5]
\begin{scope}
\fill[black!20]
(0,0)rectangle(1,1)
(2,1)--(4,1)--(4,3)--(3,3)--(3,2)--(2,2)--cycle
(4,4)--(5,4)--(5,5)--(6,5)--(6,6)--(4,6)--cycle;
\draw[gray](0,0)grid(6,6);
\draw[very thick] (0,0)--(0,1)--(4,1)--(4,6)--(6,6);
\end{scope}
\begin{scope}[xshift=10cm,yshift=1cm]
\fill[black!20](2,2)--(3,2)--(3,3)--(4,3)--(4,4)--(2,4)--cycle;
\draw[gray](0,0)grid(4,4);
\draw[very thick](0,0)--(1,0)--(1,1)--(2,1)--(2,4)--(4,4);
\end{scope}
\draw[white](16,0)--(16,1);
\end{tikzpicture}
\caption{The lattice paths with type $B$ area vectors $(-1,2,1,0,-1,3)$ and $(0,0,-1,3)$.}
\label{Figure:areaB}
\end{center}
\end{figure}
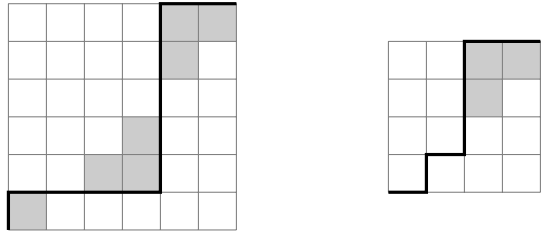

\begin{myex}{areaB} Consider the two lattice paths in \reff{areaB}. In case of the first path we have $\lambda=(0,4,4,4,4,4)$ and $n=6\equiv2$ modulo $4$. Hence we compute $\mu=(-1,2,1,0,-1,3)$. The corresponding translation in $\affS_6^B$ is given by $t_{\mu}=[14,-24,-10,4,18,-33]$. The unique signed permutation $\sigma\in\S_6^B$ such that $t_{\mu}\sigma$ is Gra{\ss}mannian is given by $\sigma=[4,-3,1,5,-2,-6]$. Thus we obtain $\waf_D^{-1}=t_{\mu}\sigma=[4,10,14,18,24,33]$ and $\waf_D=[-10,21,11,1,-9,33]$.
In particular note the difference to \refe{areaC}.

In case of the second path we find $n=4\equiv0$ modulo $4$, $\lambda=(1,2,2,7)$ and therefore $\mu=(0,0,-1,3)$. The translation $t_{\mu}=[1,2,12,-23]\in\affS_4^B$ determines a unique signed permutation $\sigma=[1,2,3,-4]\in\S_4^B$ such that $t_{\mu}\sigma$ is Gra{\ss}mannian. We obtain $\waf_D^{-1}=[1,2,12,23]$ and $\waf_D=[1,2,-6,23]$.
\end{myex}

\subsection{The diagonal reading word}\label{Section:readB}

Let $(\pi,v)\in\op{Vert}(B_n)$ be a vertically labelled lattice path. By \refp{vertB} the pair $(\pi,v)$ corresponds to an element $u\cdot\lambda+(2n+1)\Q$ of the finite torus, and thus via the Anderson map (\reft{anderson}) to a region $R$ of the Shi arrangement with minimal alcove $\waf_R\ac$. Write $\waf_R=w\waf_D$ where $w\in\S_n^B$ is a signed permutation and $\waf_D\ac$ is the minimal alcove of a dominant Shi region. Our aim for this section is to recover the signed permutation $w$, and thereby $\waf_R$, from the labelled path $(\pi,v)$.

The affine permutation $\waf_D=t_{\mu}\sigma$ can be obtained from $\lambda$ as discussed in the previous section. Thus in theory we are already able to compute $w=u\tau\sigma$ using \refp{utausig}. As in types $C_n$ and $D_n$ above, we demonstrate a convenient way to read $w$ off the vertical labelling.

\begin{mydef}{dB} Let $(\pi,v)\in\op{Vert}(B_n)$ be a vertically labelled lattice path with type $B_n$ area vector $\mu$. Define the type $B_n$ \emph{diagonal reading word} $d_B(\pi,v)$ as follows: For each $i=0,1,2,\dots$ first write down the labels $v(j)$ of the rows with $\mu_j=-i$ from bottom to top, then write down the negative labels $-v(j)$ of rows with $\mu_j=i+1$ from top to bottom. Finally, if $\mu_{n-1}+\mu_n$ is even, change the sign of the label coming from the top row.
\end{mydef}

\sk
Note that the diagonal reading word of type $B_n$ is closely related to that of type $D_n$. Indeed the above definition is the same as \refd{readD} except that some technical details are less complicated in type $B_n$.

\begin{myex}{readB} Let us recall the vertically labelled lattice paths from \reff{vertB}. The area vector of the first path is $\mu=(-1,2,1,0,-1,3)$ and its labels are given by $v=[1,-5,-4,2,3,6]$. Reading off the vertical labels according to the entries of $\mu$ yields $[2,4,1,3,5,-6]$. Since $\mu_5+\mu_6=2$ is even we need to change the sign of the label coming from the top row, that is, $v(6)=6$. We obtain $d_B(\pi,v)=[2,4,1,3,5,6]$.

In light of \refp{readB} below the reader may wish to verify that
\begin{align*}
d_B(\pi,v)
&=u\tau\sigma
=[1,-5,-4,2,3,6]\cdot[1,2,3,4,5,-6]\cdot[4,-3,1,5,-2,-6],
\end{align*}
where $\sigma\in\S_6^B$ was already computed in \refe{areaB}.

We are also able to compute the affine permutation $\waf_R=d_B(\pi,v)\waf_D=[-12,21,9,2,-10,33]$ that takes the fundamental alcove to the minimal alcove of the Shi region $R$ corresponding to the vertically labelled lattice path $(\pi,v)$. Indeed, we have $\waf_R\waf_f^{-1}=[1,47,48,54,55,58]$ and
\begin{align*}
\mathcal{A}(\waf_R)
&=-\waf_R\waf_f^{-1}\cdot0+13\Q
=(0,4,4,-4,-4,4)+13\Q,
\end{align*}
which agrees with the results in \refe{vertB}.
Again it is interesting to compare this to \refe{readC} to see the differences and similarities between type $B_n$ and $C_n$.

\sk
The second path in \reff{vertB} has area vector $\mu=(0,0,-1,3)$ and its labels are given by $v=[-1,-4,-3,-2]$.
Reading off the labels according to the entries of $\mu$ yields $[-1,-4,-3,2]$. As above $\mu_2+\mu_3=2$ is even, thus we need to change the sign of the label $v(4)=-2$. We obtain $d_B(\pi,v)=[-1,-4,-3,-2]$.
Using once more the results from \refe{areaB}, we compute $\waf_R=d_B(\pi,v)\waf_D=[-1,-4,-12,29]$, $\waf_R\waf_f^{-1}=[8,14,15,65]$ and
\begin{align*}
\mathcal{A}(\waf_R)
&=-\waf_R\waf_f^{-1}\cdot0+9\Q
=(-1,7,-2,-2)+9\Q.
\end{align*}
\end{myex}

\sk
The following proposition asserts that the diagonal reading word is the correct signed permutation.

\begin{myprop}{readB} Let $(\pi,v)\in\op{Vert}(B_n)$ be a vertically labelled lattice path with area vector $\mu$, define $u$ as in \refq{uB}, $\tau$ as in \refl{wfB} and $\sigma$ as in \refl{sigma} such that $t_{\mu}\sigma$ is Gra{\ss}mannian. Then $d_B(\pi,v)=u\tau\sigma$.
\end{myprop}

\begin{myproof} Let $i\in[n]$ and chose $j\in[n]$ such that $\abs{d_B(\pi,v)(i)}=\abs{v(j)}$. Then
\begin{align*}
i
&=\#\big\{r\in[n]:\abs{\mu_r}<\abs{\mu_j}\big\}
+\#\big\{r\in[n]:j\leq r,\mu_r=\mu_j>0\big\} \\
&\quad+\#\big\{r\in[j],\mu_r=\mu_j\leq0\big\}
+\#\big\{r\in[n],\mu_r=-\mu_j>0\big\} \\
&=\#\big\{r\in[r]:\abs{\mu_r\n+r}\leq\abs{\mu_j\n-j}\big\}
=\abs{\sigma^{-1}(j)}
\end{align*}
Thus $\abs{d_B(\pi,v)(i)}=\abs{v\sigma(i)}=\abs{u\tau\sigma(i)}$. If $j<n$ then $d_B(\pi,v)=v(j)=u\tau(j)$ if and only if $\mu_j\leq0$, which is the case if and only if $\sigma(i)=j$. If $j<n$ and $\mu_j>0$ then $d_B(\pi,v)=-u\tau(j)$ and $\sigma(i)=-j$. If $j=n$ and $\mu_n\leq0$ then
\begin{align*}
d_B(\pi,v)
&=(-1)^{1+\mu_{n-1}+\mu_n}v(n)\\
&=(-1)^{1+\mu_{n-1}+\mu_n}(-1)^{\lambda_{n-1}+\lambda_n}u(n)\\
&=(-1)^{1+\nu_{n-1}+\nu_n}u(n)\\
&=u\tau(n)
=u\tau(\sigma(j))
\end{align*}
and analogously one treats the case where $j=n$ and $\mu_n>0$.
\end{myproof}

\subsection{The zeta map}\label{Section:zetaB}

We can now define the (Haglund--Loehr) zeta map of type $B_n$, which is closely related to type $D_{n+1}$.

\begin{mydef}{zetaB} Given a lattice path $\pi\in\mathcal L_{n,n}$ with type $B_n$ area vector $\mu$ we define its image under the type $B_n$ \emph{zeta map} as
\begin{align*}
\zeta_B(\pi)
=\LS_{2n}^-(\mu)\RS_{2n}^+(\mu)
\LS_{2n-1}^-(\mu)\RS_{2n-1}^+(\mu)\cdots
\LS_1^-(\mu)\RS_1^+(\mu)
\LS_0^-(\mu)\big(N\RS_0^+(\mu)\big)^{\circ},
\end{align*}
where $(N\RS_0^+(\mu))^{\circ}$ is obtained from $N\RS_0^+(\mu)$ by deleting the last letter.

Moreover, given a vertically labelled lattice paths $(\pi,v)\in\op{Vert}(B_n)$ we define its image under the \emph{Haglund--Loehr-zeta map} of type $B_n$ as
\begin{align*}
\zeta_B(\pi,v)=(\zeta_B(\pi),d_B(\pi,v)).
\end{align*}
\end{mydef}

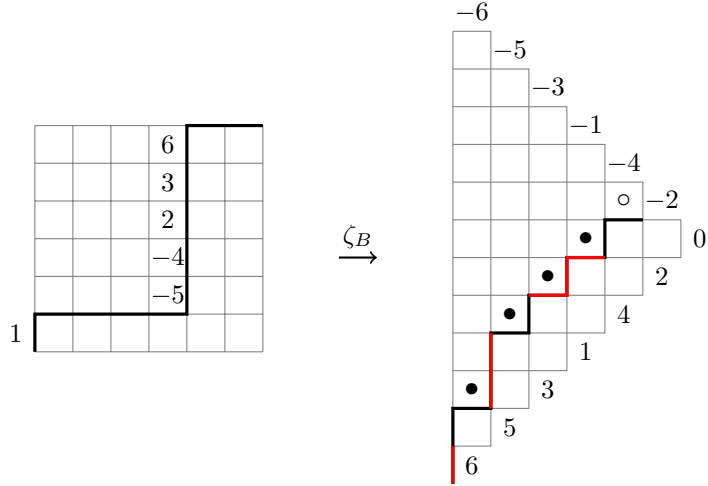
\begin{figure}[ht]
\begin{center}
\begin{tikzpicture}[scale=.5]
\begin{scope}
\draw[gray](0,0)grid(6,6);
\draw[very thick](0,0)--(0,1)--(4,1)--(4,6)--(6,6);
\draw[xshift=5mm,yshift=5mm]
(-1,0)node{$1$}
(3,1)node{$-5$}
(3,2)node{$-4$}
(3,3)node{$2$}
(3,4)node{$3$}
(3,5)node{$6$};
\end{scope}
\draw[xshift=8cm,yshift=2.5cm,->,thick](0,0)--node[above]{$\zeta_B$}(1,0);
\begin{scope}[xshift=11cm,yshift=-35mm]
\ballot{6}
\draw[very thick](0,0)--(0,2)--(1,2)--(1,4)--(2,4)--(2,5)--(3,5)--(3,6)--(4,6)--(4,7)--(5,7);
\draw[very thick,red]
(0,0)--(0,1)
(1,2)--(1,4)
(2,5)--(3,5)--(3,6)--(4,6);
\draw[xshift=5mm,yshift=5mm]
(0,0)node{$6$}
(1,1)node{$5$}
(2,2)node{$3$}
(3,3)node{$1$}
(4,4)node{$4$}
(5,5)node{$2$}
(6,6)node{$0$}
(5,7)node{$-2$}
(4,8)node{$-4$}
(3,9)node{$-1$}
(2,10)node{$-3$}
(1,11)node{$-5$}
(0,12)node{$-6$}
(0,2)node{\large{$\bullet$}}
(1,4)node{\large{$\bullet$}}
(2,5)node{\large{$\bullet$}}
(3,6)node{\large{$\bullet$}}
(4,7)node{\large{$\circ$}};
\end{scope}
\end{tikzpicture}
\caption{A vertically labelled lattice path and its image under the Haglund--Loehr-zeta map.}
\label{Figure:zetaB1}
\end{center}
\end{figure}

\begin{myex}{zetaB} Consider the vertically labelled lattice path $(\pi,v)\in\op{Vert}(B_6)$ drawn in \reff{zetaB1}, which has area vector $\mu=(-1,2,1,0,-1,3)$. We compute
\begin{align*}
\RS_3^+(\mu)=N,
\LS_2^-(\mu)&=\emptyset,
\RS_2^+(\mu)=NE,
\LS_1^-(\mu)=NN,\\
\RS_1^+(\mu)=EN,
\LS_0^-(\mu)&=ENE,
\RS_0^+(\mu)=EN,
\end{align*}
Note that the last segment is replaced by $\big(N\RS_0^+(\mu)\big)^{\circ}=NE$.
In light of \reft{risevalleyB} below we remark that the initial North step of $\pi$ labelled by $v(1)=1$ corresponds to a valley labelled $(1,0)$.
Moreover the rises of $(\pi,v)$, which are labelled $(-5,-4),(-4,2),(2,3)$ and $(3,6)$, correspond exactly to the remaining valleys of $(\beta,w)$.
\end{myex}

\sk
Note that $\zeta_B(\pi)\in\mathcal B_{2n}$ is by definition a ballot path.

\begin{mythrm}{zetaB} The zeta map $\zeta_B:\mathcal L_{n,n}\to\mathcal B_{2n}$ is a bijection.
\end{mythrm}

\begin{myproof}
Let $\pi\in\mathcal L_{n,n}$ be a lattice path with type $B_n$ area vector $\mu$. Consider the path $N\pi\in\mathcal L_{n,n+1}$. Recalling the bijection $\zeta_D^*:\mathcal L_{n,n+1}\to\mathcal B_{2n+1}$ from \reft{oddballotbij} we have
\begin{align*}
\zeta_D^*(N\pi)
=\LS_{2n}^-(\mu)\RS_{2n}^+(\mu)
\LS_{2n-1}^-(\mu)\RS_{2n-1}^+(\mu)\cdots
\LS_1^-(\mu)\RS_1^+(\mu)
\LS_0^-(\mu)N\big(N\RS_0^+(\mu)\big)^{\circ}.
\end{align*}
It follows from the proof of \reft{oddballotbij} that $\zeta_D^*$ restricts to a bijection from the set of lattice paths in $\mathcal L_{n,n+1}$ that begin with a North step to the set of ballot paths in $\mathcal B_{2n+1}$ whose $(n+1)$-st North step is not followed by an East step. Since $\zeta_D^*(N\pi)$ is easily obtained from $\zeta_B(\pi)$ by inserting a North step, we conclude that $\zeta_B$ is also a bijection.
\end{myproof}

\sk
The next theorem further exploits the relation between $\zeta_B$ and $\zeta_D$ to connect the rises of $(\pi,v)$ to the valleys of $\zeta_B(\pi,v)$.

Let $(\pi,v)$ be a vertically labelled lattice path and $i$ a rise of $\pi$. Recall that we say the rise $i$ is labelled by $(v(i),v(i+1))$. On the other hand let $(\beta,w)$ be a pair of a ballot path $\beta\in\mathcal B_{2n}$ and a signed permutation $w\in\S_n^B$, and let $(i,j)$ be a valley of $\beta$. Then we say $(\beta,w)$ has a valley $(i,j)$ that is labelled by $(w(n+1-i),w(n+1-j))$.

\begin{mythrm}{risevalleyB} Let $(\pi,v)\in\op{Vert}(B_n)$ be a vertically labelled lattice path and $a,b\in v([n])$. Then $(\pi,v)$ has a rise labelled $(a,b)$ if and only if $\zeta_B(\pi,v)$ has a valley labelled $(b,a)$ or $(-a,-b)$. Furthermore $\pi$ begins with a North step if and only if $\zeta_B(\pi,v)$ has a valley labelled $(v(1),0)$.
\end{mythrm}

\begin{myproof} Let $(\pi,v)\in\op{Vert}(B_n)$ and consider the vertically labelled path $(N\pi,\bar{v})$ where we label the initial North step by $0$ and $\pi$ retains the labelling $v$. This is basically an element of $\op{Vert}(D_{n+1})$ except that the absolute values of all labels have been decreased by one.

We will prove the claim using \reft{risevalleyD}. First notice that $\zeta_B(\pi,v)$ is obtained from $\zeta_D(N\pi,\bar{v})$ simply by deleting the $(n+1)$-st North step and its label, which is always $0$. Furthermore let $a,b\in\Z-\{0\}$, then $(\pi,v)$ has a rise $i$ labelled $(a,b)$ if and only if $(N\pi,\bar{v})$ has a rise $i+1$ labelled $(a,b)$. This is the case if and only if $\zeta_D(N\pi,\bar{v})$ has a valley labelled $(b,a)$ or $(-a,-b)$ and equivalently $\zeta_B(\pi,v)$ has a valley labelled $(a,b)$ or $(-a,-b)$. On the other hand $\pi$ begins with a North step if and only if $(N\pi,\bar{v})$ has rise $i=1$ labelled by $(\pm0,v(1))$. This is equivalent to $\zeta_D(N\pi,\bar{v})$ having a valley labelled $(v(1),\pm0)$, which is the case if and only if $\zeta_B(\pi,v)$ has a valley labelled $(v(1),0)$.
%
%
%
%
\end{myproof}

\sk
From \reft{zetaB} and \reft{risevalleyB} we obtain the following.

\begin{mythrm}{HLzetaB} The Haglund--Loehr-zeta map $\zeta_B:\op{Vert}(B_n)\to\op{Diag}(B_n)$ is a bijection.
\end{mythrm}

\sk
We conclude this section by proving that the combinatorial zeta $\zeta_B$ map is indeed the type $B_n$ instance of the uniform zeta map.

\begin{mythrm}{uniformB} Let $\Phi$ be the root system of type $B_n$ with coroot lattice $\Q$ and zeta map $\zeta$. Moreover define $\psi$ and $\varphi$ as in \refp{vertB} respectively \refp{diagB}. Then the following diagram commutes.
\begin{center}
\begin{tikzpicture}
\draw (0,0) node(q){$\Q/(2n+1)\Q$};
\draw (4,0) node(park){$\mathrm{Park}(\Phi)$};
\draw (0,1.4) node(v){$\op{Vert}(B_n)$};
\draw (4,1.4) node(d){$\op{Diag}(B_n)$};
\draw[->] (q)--node[below]{$\zeta$}(park);
\draw[->] (v)--node[above]{$\zeta_B$}(d);
\draw[->] (v)--node[left]{$\psi$}(q);
\draw[->] (d)--node[right]{$\varphi$}(park);
\end{tikzpicture}
\end{center}
\end{mythrm}

\begin{myproof} Let $(\pi,v)\in\op{Vert}(B_n)$ and set $(\beta,w)=\zeta_B(\pi,v)\in\op{Diag}(B_n)$. 
Define $\lambda$ and $u$ as in \refq{lambdaB} respectively \refq{uB} such that $\psi(\pi,v)=u\cdot\lambda+(2n+1)\Q$.
Let $\tau$ be as in \refl{wfB} and $\mu$ be the type $B_n$ area vector of $\pi$.
Chose $\sigma\in\S_n^B$ such that $t_{\mu}\sigma\in\affS_n^B$ is a Gra{\ss}mannian affine permutation.

Recall that by \refp{antichain} and \refp{utausig} it suffices to show that $w=u\tau\sigma$ and $A_{\beta}=A$, where
\begin{align*}
A
&=\big\{(\tau\sigma)^{-1}\cdot\alpha:\alpha\in\Delta\cup\{-\tilde\alpha\}\text{ and }s_{\alpha}\cdot\lambda=\lambda\big\}.
\end{align*}
The first claim is taken care of by \refp{readB}.
In order to demonstrate the second claim we first show $A\subseteq A_{\beta}$.
Therefore let $i\in[n-1]$ and suppose that $(\tau\sigma)^{-1}\cdot\alpha_i\in A$.
Then $\lambda_i=\lambda_{i+1}$, hence $\pi_i=\pi_{i+1}$ and $i$ is a rise of $\pi$.
By \reft{risevalleyB} $(\beta,w)$ has a valley $(x,y)$ labelled either $(v(i+1),v(i))$ or $(-v(i),-v(i+1))$.
Note that in particular $y\neq n+1$.
Moreover $u(i)=v(i)$ and $u(i+1)=v(i+1)$ for this choice of $i$.
We obtain
\begin{align*}
(w(n+1-x),w(n+1-y))
\in\big\{(u(i+1),u(i)),(-u(i),-u(i+1))\big\}.
\end{align*}
Applying $u^{-1}$ and using the fact that $w=u\tau\sigma$, yields
\begin{align*}
(\tau\sigma(n+1-x),\tau\sigma(n+1-y)
\in\big\{(i+1,i),(-i,-i-1)\big\}.
\end{align*}
Recalling that
\begin{align*}
\alpha_{x,y}=
\begin{cases}
e_{n+1-x}+e_{n+1-y}&\quad\text{if }y<n+1,\\
e_{n+1-x}-e_{y-n-1}&\quad\text{if }y<n+1,
\end{cases}
\end{align*}
we compute $(\tau\sigma)^{-1}\cdot\alpha_i=\alpha_{x,y}\in A_{\beta}$ in all cases.

Next suppose $(\tau\sigma)^{-1}\cdot\alpha_0\in A$.
Then $\lambda_1=0$ and thus $\pi_1=0$.
By \reft{risevalleyB} $(\beta,w)$ has a valley $(x,n+1)$ labelled $(v(1),0)=(u(1),0)$.
From $u\tau\sigma(n+1-x)=w(n+1-x)=u(1)$ we obtain $\tau\sigma(n+1-x)=1$, and compute $(\tau\sigma)^{-1}\cdot\alpha_0=e_{n+1-x}=\alpha_{x,n+1}\in A_{\beta}$.

Similarly suppose $(\tau\sigma)^{-1}\cdot\alpha_n\in A$.
Then $\lambda_{n-1}=\lambda_n$. 
Consequently $\pi_{n-1}=\pi_n$, that is, $n-1$ is a rise of $\pi$ and $\pi_1+\dots+\pi_{n-2}$ is even.
In particular $u(n)=v(n)$.
By \reft{risevalleyB} $(\beta,w)$ has a valley $(x,y)$ labelled $(u(n),u(n-1))$ or $(-u(n-1),-u(n))$.
As in the cases above we conclude that $(\tau\sigma)^{-1}\cdot\alpha_{n-1}=\alpha_{x,y}\in A_{\beta}$.

Finally suppose $(\tau\sigma)^{-1}\cdot(-\tilde\alpha)\in A$.
Then $\lambda_{n-1}+\lambda_n=2n+1$.
Here again $\pi_{n-1}=\pi_n$, but contrary to the previous case $\pi_1+\dots+\pi_{n-2}$ is now odd.
Therefore $u(n)=-v(n)$.
By \reft{risevalleyB} $(\beta,w)$ has a valley labelled $(u(n-1),-u(n))$ or $(u(n),-u(n-1))$.
As before we compute $(\tau\sigma)^{-1}\cdot(-\tilde\alpha)=\alpha_{x,y}\in A_{\beta}$.

\sk
It remains to prove the reverse inclusion $A_{\beta}\subseteq A$.
Therefore assume $\alpha_{x,y}\in A_{\beta}$ for some valley $(x,y)$ of $(\beta,w)$.

If $y<n+1$ then $\alpha_{x,y}=e_{n+1-x}-e_{n+1-y}$.
Furthermore by \reft{risevalleyB} the valley $(x,y)$ is labelled
\begin{align*}
(w(n+1-x),w(n+1-y))
\in\big\{(v(i+1),v(i)),(-v(i),-v(i+1))\big\}
\end{align*}
for some rise $i$ of $\pi$.
If $i<n-1$ or if $\pi_1+\dots+\pi_{n-2}$ is even, then $\lambda_i=\pi_i=\pi_{i+1}=\lambda_{i+1}$ and $v(i)=u(i)$ and $v(i+1)=u(i+1)$.
Thus $s_{\alpha_i}\cdot\lambda=\lambda$ and we obtain
\begin{align*}
\alpha_{x,y}
&=(\tau\sigma)^{-1}\cdot\alpha_i
\in A.
\end{align*}
If $i=n-1$ and $\pi_1+\dots+\pi_{n-2}$ is odd, then $\lambda_{\n-1}+\lambda_n=2n+1$ and $u(n)=-v(n)$.
We obtain $s_{\tilde\alpha}\cdot\lambda=\lambda$ and $\alpha_{x,y}=(\tau\sigma)^{-1}\cdot(-\tilde\alpha)\in A$. 

If $y=n+1$ then $\alpha_{x,y}=e_{n+1-x}$ and the valley $(x,y)$ is labelled $(w(n+1-x,0)=(v(1),0)$ by \reft{risevalleyB}.
Moreover $\pi$ begins with a North step, that is, $\lambda_1=\pi_1=0$, and $s_{\alpha_0}\cdot\lambda=\lambda$.
We conclude $\alpha_{x,y}=(\tau\sigma)^{-1}\cdot\alpha_0\in A$.

Finally the case $y>n+1$ can be treated in a similar fashion as the case $y<n+1$ above, which completes the proof.
\end{myproof}

%% file: conclusion.tex
\section{Open problems}\label{Section:end}

In this final section we discuss some open problems and perspectives for further research.

\begin{myprob}[ (The Fu{\ss}--Catalan case)]{fuss}
In this paper we discuss the combinatorial interpretation of the bijection $\zeta$ between the finite torus $\Q/(h+1)\Q$ and the set of non-nesting parking functions $\op{Park}(\Phi)$.
However, the zeta map can be defined more generally.
Rhoades~\cite{Rhoades2014} defined a Fu{\ss} analogue of the non-nesting parking functions $\op{Park}^{(m)}(\Phi)$ and proved it is isomorphic to the generalised finite torus $\Q/(mh+1)\Q$ as a $W$-set.
The second named author~\cite{Thiel2015} defined an explicit $W$-equivariant bijection $\zeta:\op{Park}^{(m)}(\Phi)\to\Q/(mh+1)\Q$ generalising the zeta map treated in this work.

In type $A_{n-1}$ the combinatorial objects allowing for a convenient treatment of the Fu{\ss}--Catalan level are $m$-Dyck paths, and the combinatorial zeta map was studied by Loehr~\cite{Loehr2005}.
One possible next step would be to develop the combinatorial framework for the Fu{\ss}--Catalan case of the zeta map for the other infinite families of crystallographic root systems $B_n,C_n$ and $D_n$.

\end{myprob}

\begin{myprob}[ (The inverse zeta map)]{inverse} \refs{typeC} contains an explicit way to invert the zeta map of type $C_n$, namely \reft{zetaC}.
As the attentive reader may have noticed, such a result is missing from the respective chapters on types $B_n$ and $D_n$.
Indeed our only way of establishing the bijectivity of the zeta maps in these types is to prove that they are special cases of the uniform zeta map which is known to be bijective.
It would be interesting to have an explicit combinatorial description of the inverse of the zeta map in types $B_n$ and $D_n$.
One might expect such a description to be connected to the construction of a suitable bounce path, which also appears in types $A_{n-1}$ and $C_n$.
\end{myprob}

\begin{myprob}[ (The $\dinv$ statistic)]{dinv} As pointed out in \refs{statC} the definition of the $\area$ statistic is uniform, while the definition of the $\dinv$ statistic is particular to the types $A_{n-1}$ respectively $C_n$.
One (rather uninspired) way to give a uniform definition of the $\dinv$ statistic is as the pullback of the $\area$ statistic by the zeta map, that is,
\begin{align*}
\dinv=\area\circ\zeta.
\end{align*}
Thereby an analogous result to \reft{dinvC} is fulfilled automatically.
A next logical step would be to ask for a description of the statistics $\dinv'$ and $\dinv$ that does not use the zeta map but can be seen directly from the combinatorial models for the finite torus respectively its $W$-orbits in types $B_n$ and $D_n$.
We suspect this will not prove too difficult.
As a starting point we suggest the following definition, which is similar to the type $C_n$ case:
\begin{align*}
\dinv_B(\pi)
&=\#\big\{(i,j):1\leq i<j\leq n,\mu_i=\mu_j\big\}
+\#\big\{(i,j):1\leq i<j\leq n,\mu_i=\mu_j-1\big\}\\
&\qquad+\#\big\{(i,j):1\leq i<j\leq n,-\mu_i=\mu_j\big\}
+\#\big\{(i,j):1\leq i<j\leq n,-\mu_i=\mu_j-1\big\}\\
&\qquad+\#\big\{i:1\leq i\leq n,\mu_i\in\{0,1\}\big\}.
\end{align*}
It might be a more challenging task to find a uniform definition of the statistics $\dinv'$ and $\dinv$ that can be seen directly on the finite torus respectively its $W$-orbits, which does not rely on the zeta map.
\end{myprob}

%